\documentclass[11pt,reqno]{amsart}
\usepackage[left=33mm,right=33mm,top=30mm,bottom=32mm]{geometry}
\usepackage{mathtools,amssymb,amsthm,mathrsfs,color,float}
\usepackage{paralist}
\usepackage{stackengine}
\usepackage{centernot}
\usepackage{txfonts}
\usepackage{tabularx}
\usepackage{longtable}
\usepackage[colorlinks,
linkcolor=red,
anchorcolor=green,
citecolor=blue, 
]{hyperref}

\usepackage[T1]{fontenc}
\usepackage[utf8]{inputenc}
\usepackage{booktabs}
\usepackage{calc}
\definecolor{bleu1}{RGB}{0,57,128}
\def\bleu1{\color{bleu1}}

\usepackage{etoolbox}
\patchcmd{\section}{\normalfont}{\normalfont \bleu1}{}{}
\patchcmd{\subsection}{\normalfont}{\normalfont \bleu1}{}{}
\patchcmd{\subsubsection}{\normalfont}{\normalfont \bleu1}{}{}

\newtheorem{The}{\bleu1 Theorem}[section]
\newtheorem{Pro}{\bleu1 Proposition}[section]
\newtheorem{Lem}{\bleu1 Lemma}[section]
\newtheorem{Cor}{\bleu1 Corollary}[section]

\theoremstyle{definition}
\newtheorem{defn}{\bleu1 Definition}[section]

\newtheorem{Rem}{\bleu1 Remark}[section]

\newcommand{\T}{\mathbb{T}}
\newcommand{\R}{\mathbb{R}}

\newcommand{\E}{\mathbb{E}}

\newcommand{\A}{\mathbb{A}}

\newcommand{\Cut}{\mathrm{Cut}\,}

\newcommand{\SCL}{\mathrm{SCL\,}}
\newcommand{\Lip}{\mathrm{Lip\,}}
\newcommand{\Liploc}{\mathrm{Lip\,}_{\mathrm{loc}}}
\newcommand{\Sing}{\mathrm{Sing\,}}
\newcommand{\law}{\mathrm{law\,}}

\newcommand{\UC}{\mathrm{UC}\,}
\newcommand{\supp}{\mathrm{supp\,}}

\makeatletter
\@namedef{subjclassname@2020}{\textup{2020} Mathematics Subject Classification}
\makeatother

\title[Singularities and their propagation in optimal transport]{Singularities and their propagation in optimal transport}
\author{Piermarco Cannarsa, Wei Cheng, Tianqi Shi \and Wenxue Wei}
\address[Piermarco Cannarsa]{Dipartimento di Matematica, Universit\`a di Roma ``Tor Vergata'', Via della Ricerca Scientifica 1, 00133 Roma, Italy}
\email{cannarsa@mat.uniroma2.it}
\address[Wei Cheng]{School of Mathematics, Nanjing University, Nanjing 210093, China}
\email{chengwei@nju.edu.cn}
\address[Tianqi Shi]{Department of Mathematical Sciences, Tsinghua University, Beijing 100084, China}
\email{tqshi.math@gmail.com}
\address[Wenxue Wei]{School of Mathematics, Nanjing University, Nanjing 210093, China}
\email{wwx3708@gmail.com}
\keywords{optimal transport, potential energy functional, cut locus, propagation of singularities, Hamilton-Jacobi equation}
\subjclass[2020]{35D40,49Q22,	37Jxx,37Kxx}

\begin{document}
\maketitle

%\tableofcontents

\begin{abstract}
In this paper, we investigate the singularities of potential energy functionals \(\phi(\cdot)\) associated with semiconcave functions \(\phi\) in the Borel probability measure space and their propagation properties. Our study covers two cases: when \(\phi\) is a semiconcave function and when \(u\) is a weak KAM solution of the Hamilton-Jacobi equation \(H(x, Du(x)) = c[0]\) on a smooth closed manifold. By applying previous work on Hamilton-Jacobi equations in the Wasserstein space, we prove that the singularities of \(u(\cdot)\) will propagate globally when \(u\) is a weak KAM solution, and the dynamical cost function \(C^t\) is the associated fundamental solution. We also demonstrate the existence of solutions evolving along the cut locus, governed by an irregular Lagrangian semiflow on the cut locus of \(u\).
\end{abstract}

\section{Introduction}

In principle, the formation and evolution of singularities in optimal control, partial differential equations (PDEs), and geometry arise from the crossing and focusing of underlying characteristics or geodesics. Once formed, these singularities will propagate. This phenomenon reflects a certain irreversibility in the evolution of Hamiltonian systems. In a series of papers \cite{Albano_Cannarsa2002,Cannarsa_Yu2009,ACNS2013,Cannarsa_Cheng_Fathi2017,Cannarsa_Cheng_Fathi2021,Khanin_Sobolevski2016,Stromberg2013}, the authors developed a theory to study the propagation of singularities using the concept of \emph{generalized characteristics} (see also \cite{Dafermos_book2005} for this concept in the context of hyperbolic conservation laws). Recently, various notions of generalized characteristics have been clarified, and some refined regularity properties have been established in \cite{CCHW2024}, using ideas from abstract gradient flow theory. In this paper, we aim to study the problem of singularities propagation in optimal transport by applying the ideas from \cite{CCHW2024}, along with the concept of singularities of functionals in probability measure space.

Over the past 30 years, the vibrant development of optimal transport theory has been driven by a series of geometric and dynamical perspectives. A pivotal advancement was made through the work of Gangbo and McCann (\cite{Gangbo_McCann1996}), and later by McCann (\cite{McCann1997}), who introduced the concept of displacement interpolation under a quadratic cost function in Euclidean space. Shortly thereafter, Benamou and Brenier reformulated these findings as an action-minimization problem within the space of measures (\cite{Benamou_Brenier1999,Brenier2003}). Otto, in a series of pioneering papers (\cite{Otto1998, Otto1999, Otto2001}), introduced the notion of gradient flows in Wasserstein space and explicitly formulated geodesic equations in a Riemannian geometric setting. Building on this foundation, Ambrosio, Gigli, and Savar\'e systematically developed the theory of gradient flows in metric spaces (\cite{Ambrosio_GigliNicola_Savare_book2008}). At this stage, the close connection between Hamiltonian systems and optimal transport became evident. Bernard and Buffoni constructed displacement interpolation with cost functions determined by Tonelli Lagrangians, establishing a clear link between optimal transport and Aubry-Mather theory in Hamiltonian dynamical systems (\cite{Bernard_Buffoni2006, Bernard_Buffoni2007a, Bernard_Buffoni2007b}). Additionally, Ambrosio, Gangbo, and Tudorascu, utilizing the geometry of the 2-Wasserstein space, considered Hamiltonians and flows within its phase space, deriving continuity equations, Euler-Poisson equations, and energy conservation results (\cite{Ambrosio_Gangbo2008, Gangbo_Nguyen_Tudorascu2009}). This work further contributed to the development of weak KAM theory in Wasserstein space by Gangbo and Tudorascu (\cite{Gangbo_Tudorascu2010a, Gangbo_Tudorascu2010b}).

With the continued development of calculus theory on Wasserstein and metric spaces, there has been growing interest in studying viscosity solutions to Hamilton-Jacobi equations on these spaces. Gangbo, Nguyen, and Tudorascu were the first to define viscosity solutions for a class of Hamilton-Jacobi equations in Wasserstein space using subdifferentials and superdifferentials, and to investigate their well-posedness (\cite{Gangbo_Nguyen_Tudorascu2008, Gangbo_Tudorascu2019}). Additionally, using the concept of metric derivatives, Ambrosio and Feng, Gangbo and \'{S}wi\c{e}ch, as well as other researchers, have explored the well-posedness of viscosity solutions to Hamilton-Jacobi equations in general metric spaces (\cite{Ambrosio_Feng2014, Gangbo_Swiech2015a, Gangbo_Swiech2015b, Giga_Hamamuki_Nakayasu2015}). Building on these foundational results, some related dynamical aspects, such as the asymptotic behavior of viscosity solutions and Hamilton-Jacobi equations with cohomological terms, have been examined (\cite{Nakayasu_Namba2018, Bessi2020}). In the areas of calculus of variations and optimal (or mean-field) control, we refer to \cite{Cardaliaguet_Souganidis2023, Bonnet_Frankowska2021b, Bonnet_Frankowska2021a, Bonnet_Frankowska2022, Badreddine_Frankowska2022, Hynd_Kim2015b, Hynd_Kim2015a, CGKPR2024, Feng_Swiech2013}, where some works provide variational representations of viscosity solutions.%, due to our limited knowledge on these topic.

The study of the geometric aspects of optimal transport, with an emphasis on the cut locus, has been the subject of numerous previous works. Under the assumption that the cost function is the square of the distance, the existence of the optimal transport map from an absolutely continuous (with respect to Lebesgue measure) measure was obtained in \cite{Brenier1991} in Euclidean case and \cite{McCann2001} in Riemannian case. The first-order properties of the cut locus, such as Lipschitz continuity and $(n-1)$-rectifiability, have been investigated in \cite{Castelpietra_Rifford2010, Itoh_Tanaka2001, Li_Nirenberg2005, Mantegazza_Mennucci2003}. The authors in \cite{Loeper_Villani2010} study the second-order characterization of the cut locus acting as a barrier to the regularity of optimal transport maps, including the semiconvexity of the injectivity domain. By virtue of the Ma-Trudinger-Wang (MTW) condition, for a Riemannian manifold \( M \) with a nonfocal cut locus satisfying the strong MTW condition, the following results hold:
\begin{itemize}
    \item For any two \( C^\infty \) positive probability densities \( f \) and \( g \) on \( M \), the optimal transport map from \( \mu(dx) = f(x) \operatorname{vol}(dx) \) to \( \nu(dy) = g(y) \operatorname{vol}(dy) \) with cost function \( c = d^2 \) is \( C^\infty \);
    \item There exists a \( \kappa > 0 \) such that all injectivity domains are \( \kappa \)-uniformly convex.
\end{itemize}
Furthermore, by employing a slightly modified version of the MTW condition, a more delicate characterization of the injectivity domains and the regularity of optimal transport maps is obtained (see \cite[Theorem 1.8]{Loeper_Villani2010}). If \((M, g)\) is a surface, the authors in \cite{Figalli_Rifford_Villani2011} prove that the injectivity domains of \(M\) are semiconvex if \(M\) satisfies an appropriate signed curvature condition. In \cite{Figalli_Rifford2009}, the authors demonstrated that any \( C^4 \)-deformation of the round sphere \((\mathbb{S}^2, g^{\text{can}})\) satisfies the transport continuity property. In their subsequent work \cite{Figalli_Rifford_Villani2012}, they further showed that if \((M, g)\) is a \(C^4\)-deformation of the round sphere \((\mathbb{S}^n, g^{\text{can}})\), then all injectivity domains of \(M\) are uniformly convex. For the verification of MTW condition on Riemannian surfaces, see also \cite{Du_Li2014}.

To address singularities and cut locus of functionals on Wasserstein spaces, we follow the framework in \cite{Bonnet_Frankowska2022}. For a metric space \((X, d)\), we denote by \(\mathscr{P}(X)\) the set of all probability Borel measures on \(X\), and by \(\mathscr{P}_c(\mathbb{R}^m)\) the subspace of \(\mathscr{P}(\R^m)\) consisting of measures with compact support. For any functional \(U: \mathscr{P}_c(\mathbb{R}^m) \to \mathbb{R}\), we introduce \(\partial^{\pm} U(\mu)\), the localized Fréchet superdifferential and subdifferential of \(U\) at \(\mu\) (see Definition \ref{def:local frechet}). A measure \(\mu \in \mathscr{P}_c(\mathbb{R}^m)\) is called a regular point of \(U\) if \(U\) is locally differentiable at \(\mu\), i.e., both \(\partial^+ U(\mu)\) and \(\partial^- U(\mu)\) are non-empty. A measure \(\mu \in \mathscr{P}_c(\mathbb{R}^m)\) is called a singular point of \(U\) if \(U\) is not locally differentiable at \(\mu\). We denote by \(\mathscr{S}(U)\) the set of all singular points of \(U\). We also introduce the notions of semiconcavity for a functional \(U: \mathscr{P}_c(\mathbb{R}^m) \to \mathbb{R}\) (see Definition \ref{def:local-scl}).

We concentrate to the potential energy functional (see \cite{Ambrosio_GigliNicola_Savare_book2008}) induced by a function $\phi:\R^m\to\R$
\begin{align*}
	\Phi(\mu):=\int_{\R^m}\phi(x)\,d\mu,\,\,\,\,\forall\mu\in\mathscr P(\R^m).
\end{align*}
We denote $\Phi(\cdot)$ by $\phi(\cdot)$ for brevity. That is $\phi$ is real-valued function on $\R^m$ and $\phi(\cdot)$ is also a functional on $\mathscr P(\R^m)$ if there is no confusion. If $\phi$ is a locally semiconcave function on $\R^m$, its associated potential energy functional $\phi(\cdot)$ has the following properties:
\begin{enumerate}
	\item If $\phi$ is locally semiconcave, then potential energy functional $\phi(\cdot)$ is locally strongly semiconcave. (Proposition \ref{pro:potential-scl})
	\item If $\phi$ is locally semiconcave and $\mu\in\mathscr P_c(\R^m)$, then $\alpha\in\partial^+\phi(\mu)$ if and only if $\alpha(x)\in D^+\phi(x)$, the superdifferential of $\phi$ at $x$, holds true for $\mu$-a.e. $x\in\R^m$. (Theorem \ref{thm:partial+D+}) Moreover, $\partial^+\phi(\mu)\neq\varnothing$ and $\phi(\cdot)$ is locally differentiable at $\mu$ if and only if $\partial^+\phi(\mu)$ is a singleton, and $\alpha=\partial\phi(\mu)$ if and only if $\alpha(x)=D\phi(x)$ holds for $\mu$-a.e. $x\in\R^m$. (Corollary \ref{cor:partial+singleton})
	\item If $\phi$ is locally semiconcave and $\mu\in\mathscr P_c(\R^m)$, then $\mu\in \mathscr S(\phi(\cdot))$ if and only if $\mu(\Sing(\phi))>0$. (Theorem \ref{thm:singular measure}) Therefore, when $\mu\ll\mathscr L^m$, the Lebesgue measure on $\R^m$, $\mu\notin\mathscr S(\phi(\cdot))$.
\end{enumerate}

We aim to study the problem of propagation of singularities of the potential energy functional $u(\cdot)$ with $u$ a weak KAM solution of the Hamilton-Jacobi equation
\begin{equation}\label{eq:wkam_intro}
	H(x,Du(x))=c[0],\qquad x\in M,
\end{equation}
where $M$ is a smooth closed manifold, $H:T^*M\to\R$ is a Tonelli Hamiltonian and $c[0]$ is the Ma\~n\'e's critical value. 

Let \( L: TM \to \mathbb{R} \) be the Tonelli Lagrangian associated with \( H \), and let \( A_t(x,y) = \inf_{\eta} \int_0^t L(\eta, \dot{\eta}) \, ds \) be the fundamental solution, where the infimum is taken over all absolutely continuous curves \( \eta: [0,t] \to M \) connecting \( x = \eta(0) \) to \( y = \eta(t) \). Set \( c^t(x,y) = A_t(x,y) \) and recall the dynamical cost functional
\[
C^t(\mu, \nu) := \inf_{\gamma \in \Gamma(\mu, \nu)} \int_{M \times M} c^t(x,y) \, d\gamma. 
%= \inf_{\substack{\mathrm{law}(X) = \mu \\ \mathrm{law}(Y) = \nu}} \mathbb{E}(c^t(X,Y)),
\]
%where \( X(\omega) \) and \( Y(\omega) \) are \( M \)-valued random variables defined on the probability space \( (\Omega, \mathscr{F}, \mathbb{P}) \). For brevity, we denote by \( \Gamma_o^t(\mu, \nu) \) the set of minimizers of \( C^t(\mu, \nu) \). 
%\textcolor{red}{We say that a random curve \( \xi \) is a dynamical coupling of \( \mu \) and \( \nu \) for \( t > 0 \) if \( \xi: \Omega \to AC([0,t]; M) \) is measurable, \( \mathrm{law}(\xi(0,\cdot)) = \mu \), and \( \mathrm{law}(\xi(t,\cdot)) = \nu \).} %The set of all dynamical couplings of \( \mu \) and \( \nu \) for \( t > 0 \) is denoted by \( L_{\mu, \nu}^t \). 

Based on the regularity properties of the fundamental solution \( A_t(x, y) \), we have the following regularity results for the dynamical cost functional \( C^t(\mu, \nu) \) (Theorem \ref{thm:dyn-cost}): Suppose \( t, a, b > 0 \) with \( a < b \), and let \( K \subset \mathbb{R}^m \) be a convex and compact subset. Then, there exists \( \tau_1 > 0 \) such that:
\begin{enumerate}
    \item \( C^t(\mu, \cdot) \) and \( C^t(\cdot, \mu) \) are superlinear on \( \mathscr{P}_c(\mathbb{R}^m) \) endowed with the \( W_1 \) metric, where \( \mu \in \mathscr{P}_c(\mathbb{R}^m) \);
    \item \( (t, \nu) \mapsto C^t(\mu, \nu) \in \text{Lip}([a,b] \times \mathscr{P}(K)) \), where \( \mu \in \mathscr{P}(K) \);
    \item For \( t \in (0, \tau_1] \), there exists \( C_K > 0 \) depending on \( K \), such that for any \( \mu \in \mathscr{P}(K) \), \( C^t(\mu, \cdot) \) is locally strongly semiconcave on \( \mathscr{P}(K) \) with constant \( \frac{C_K}{t} \).
\end{enumerate}
Moreover, if \( (t, \nu) \mapsto C^t(\mu, \nu) \) is locally differentiable, the dynamical cost functional satisfies the Hamilton-Jacobi inequality
\[
D_t C^t(\mu, \nu) + \int_{\mathbb{R}^m} H(y, \partial_\nu C^t(\mu, \nu)) \, d\nu \leqslant 0.
\]
In general, a strict inequality can be obtained in the Hamilton-Jacobi inequality above.

By addressing the Rubinstein-Kantorovitch duality in terms of abstract Lax-Oleinik operators (\cite{Cheng_Hong_Shi2024}), we introduce the random Lax-Oleinik operators:
\[
P_t^+ \phi(\mu) := \sup_{\nu \in \mathscr{P}(\mathbb{R}^m)} \left\{ \phi(\nu) - C^t(\mu, \nu) \right\},
\]
\[
P_t^- \phi(\mu) := \inf_{\nu \in \mathscr{P}(\mathbb{R}^m)} \left\{ \phi(\nu) + C^t(\nu, \mu) \right\}.
\]
We consider the following Hamilton-Jacobi equation:
\[
\partial_t U(t, \mu) + \int_{\mathbb{R}^m} H(x, \partial_\mu U(t, \mu)(x)) \, d\mu = 0, \tag{PHJ$_{e}$} \label{phje}
\]
where \( (t, \mu) \in \mathbb{R}^+ \times \mathscr{P}_c(\mathbb{R}^m) \). Then, if \( \phi \) is lower semicontinuous and \( (\kappa_1, \kappa_2) \)-Lipschitz in the large, then \( U(t, \mu) := P_t^- \phi(\mu) \) satisfies the following Cauchy problem for equation \eqref{phje} in the viscosity sense (see Definition \ref{defn:viscosity}):
\[
\left\{
\begin{array}{ll}
\partial_t U(t, \mu) + \int_{\mathbb{R}^m} H\left( x, \partial_\mu U(t, \mu) \right) \, d\mu = 0, & (t, \mu) \in \mathbb{R}^+ \times \mathscr{P}_c(\mathbb{R}^m), \\
U(0, \mu) = \phi(\mu), & \mu \in \mathscr{P}_c(\mathbb{R}^m).
\end{array}
\right.
\]
Moreover, $u$ is a weak KAM solution of \eqref{eq:wkam_intro} if and only if $P_t^-u+c[0]t=u(\cdot)$ for all $t\geqslant0$.

Analogous to classic weak KAM theory (\cite{Fathi1997_1,Fathi_book}), an absolutely continuous curve \( \{\mu_s\}_{s \in I} \subset \mathscr{P}_c(\mathbb{R}^m) \) defined on some interval \( I \subset \mathbb{R} \) is called a \((u(\cdot), C, c[0])\)-calibrated curve if for any \( a, b \in I \), with \( a \leqslant b \), the following condition holds:
\[
u(\mu_b) - u(\mu_a) = C^{b-a}(\mu_a, \mu_b) + c[0](b-a).
\]
We say that \( \mu \in \mathscr{P}_c(\mathbb{R}^m) \) is a cut point of \(u(\cdot) \), the potential energy functional induced by the weak KAM solution \( u \), if all \((u(\cdot), C, c[0])\)-calibrated curves ending at \( \mu \) cannot be extended forward as \((u(\cdot), C, c[0])\)-calibrated curves. The set of all cut points of \( u(\cdot) \) is denoted by \( \mathscr{C}(u(\cdot)) \), which is called the cut locus of the functional \( u(\cdot) \) (with respect to the dynamical cost functional \( C^t \)).

In general, it holds that \( \mathscr{S}(u(\cdot)) \subset \mathscr{C}(u(\cdot))\), and \(\mu \in \mathscr{S}(u(\cdot))\) if and only if \(\mu(\Sing(u)) > 0\), where \(\Sing(u)\) denotes the set of singularities of \(u\). %Consequently, \(\mu \notin \mathscr{S}(u(\cdot))\) if \(\mu \ll \mathscr{L}^m\), meaning that \(\mu\) is absolutely continuous with respect to the Lebesgue measure \(\mathscr{L}^m\) (Theorem \ref{thm:singular measure}, Corollary \ref{cor:cut_sing}). 

Building on the work of the first two authors on intrinsic singular characteristics in \cite{Cannarsa_Cheng3} (see also \cite{Cannarsa_Cheng_Fathi2017, CCMW2019, Cannarsa_Cheng_Fathi2021}), we consider the following result: for any \( \phi \in \text{Lip}(\mathbb{R}^m) \cap \text{SCL}_{\mathrm{loc}}(\mathbb{R}^m) \) and any \( \mu \in \mathscr{P}_c(\mathbb{R}^m) \), there exists \( t_{\phi, \mu} > 0 \) depending on \( \phi \) and \( \mu \), such that the problem
\[
\mathop{\arg\max}_{\nu \in \mathscr{P}(\mathbb{R}^m)} \left\{ \phi(\nu) - C^t(\mu, \nu) \right\}, \quad t \in (0, t_{\phi, \mu}]
\]
has a unique solution (singleton). We define the curve \( \nu_{\phi, \mu}(t) \) as follows:
\[
\nu_{\phi, \mu}(t) := 
\begin{cases}
\mu, & t = 0, \\
\mathop{\arg\max}_{\nu \in \mathscr{P}(\mathbb{R}^m)} \left\{ \phi(\nu) - C^t(\mu, \nu) \right\}, & t \in (0, t_{\phi, \mu}].
\end{cases}
\]
Then, we have the following convergence result:
\[
\lim_{t \to 0^+} \nu_{\phi, \mu}(t) = \mu \quad \text{in the sense of the metric} \ W_p, \quad \text{for any fixed} \ p \geqslant 1. \tag{Theorem \ref{thm:wass-argmax}}
\]
However, the curve \( \nu_{\phi, \mu}(\cdot) \) does not generally possess a semi-group property. To address this, we invoke the notion of a maximal slope curve from abstract gradient flow theory (see, for instance, \cite{Ambrosio_GigliNicola_Savare_book2008}). Using an EDI (energy dissipation inequality) framework, we derive a minimizing movement for any given pair \( (\phi, C^t) \) (see also \cite{CCHW2024} for the case in the calculus of variation).

Now, we formulate the following main results of this paper:
\begin{enumerate}[(1)]
	\item Let $u(\cdot)$ be the potential energy functional induced by a weak KAM solution $u$ of \eqref{eq:wkam_intro}. If $\mu\in\mathscr C(u(\cdot))$, $\nu_{u,\mu}(t)$ and $t_{u,\mu}>0$ is mentioned above, then $\nu_{u,\mu}(t)\in\mathscr S(u(\cdot))$ for all $t\in(0,t_{u,\mu}]$. (Theorem \ref{thm:sing-propogation})
	\item Suppose $\phi\in\SCL(\T^m)$, $T>0$ and $\mu_0\in\mathscr P(\T^m)$. Then there exists a Lipschitz curve $\mu(\cdot):[0,T]\to\mathscr P(\T^m)$, which satisfies the following continuity equation:
	\begin{equation}\label{eq:CE_intro}
		\left\{
	\begin{array}{ll}
		\frac{d}{dt}\mu+\mathrm{div}(H_p(\cdot,\mathbf{p}_\phi^\#(\cdot))\cdot\mu)=0,\\
		\mu(0)=\mu_0,
	\end{array}\right.
	\end{equation}
	where $\mathbf{p}_\phi^\#(x):=\arg\min\{H(x,p):p\in D^+\phi(x)\}$ for any $x\in\T^m$ (Theorem \ref{thm:conti-eq}). Moreover, if $\phi=u$ in \eqref{eq:CE_intro} is a weak KAM solution of \eqref{eq:wkam_intro},
	\begin{itemize}
		\item \(\mu_0(\overline{\Sing(\phi)}) > 0\) implies \([\mu(t)](\overline{\Sing(\phi)}) > 0\) for all \(t \in [0, T]\). (Corollary \ref{cor:cut-prop}) 
		\item \(\mu_0\in\mathscr C(u(\cdot))\) implies \(\mu(t)\in\mathscr C(u(\cdot))\) for all \(t \in [0, T]\). (Corollary \ref{cor:cut-prop})
	\end{itemize}
		% In addition, if \( H(x,p) = \frac{1}{2} \langle A(x) p, p \rangle + V(x) \) is a mechanical Hamiltonian, then \(\mu(t) \in \mathscr{S}(\phi(\cdot))\) for all \(t \in [0,T]\) provided \(\mu_0 \in \mathscr{S}(\phi(\cdot))\). (Corollary \ref{cor:conti-eq-2})
	\item For such a Lipschitz curve satisfying \eqref{eq:CE_intro} we have that for any $f\in C^\infty(\T^m)$, the map $t\mapsto f(\mu(t))$ is in $\text{Lip}\,([0,T];\R)$ and its right derivative exists such that
	\begin{align*}
		\frac{d^+}{dt}f(\mu(t))=\int_{\T^m}\langle Df(x),H_p(x,\mathbf p_\phi^\#(x))\rangle\,d\mu(t),\qquad \forall t\in(0,T].
	\end{align*}
	(Theorem \ref{thm:right_derivative})
\end{enumerate}

We emphasize the following novel aspects of this work:

\begin{itemize}
    \item {Dynamics after the Formation of Singularities}: In this paper, we study the dynamics that arise after the formation of singularities in optimal transport. In this setting, the forward dynamics is no longer governed by a regular Lagrangian flow, as discussed in \cite{Ambrosio_Gigli2013}. This makes the classical DiPerna-Lions theory (which is typically used in the regular case) inapplicable to our case.
    
    \item {Challenges from Non-Compactness}: Due to the lack of local compactness in the space of probability measures, our problem cannot be treated as a straightforward extension of results from finite-dimensional settings. This lack of compactness introduces significant technical challenges that distinguish our approach from previous treatments of similar problems.
\end{itemize}

\medskip

\noindent\textbf{Acknowledgements.} Piermarco Cannarsa was partly supported by the PRIN 2022 PNRR Project P20225SP98 ``Some Mathematical Approaches to Climate Change and Its Impacts'' (funded by the European Community-Next Generation EU), CUP E53D2301791 0001. He was also supported by the INdAM (Istituto Nazionale di Alta Matematica) Group for Mathematical Analysis, Probability and Applications, and by the MUR Excellence Department Project awarded to the Department of Mathematics, University of Rome Tor Vergata, CUP E83C23000330006.Wei Cheng was partly supported by the National Natural Science Foundation of China (Grant No. 12231010). Tianqi Shi was partly supported by the ``New Cornerstone Investigator Program'' led by Jinxin Xue. Wei Cheng also thanks Jiahui Hong for helpful discussions, Bangxian Han and Shibin Chen for comments on some literatures.

%%%%%%%%%%%%%%%%%%%%%%%%%%%%%%%%
%\newpage

\begin{center}
\begin{longtable}{p{2.9cm}p{11.4cm}p{8.6cm}}
\caption{Notations}\\
	\toprule
		$M$ & smooth Riemannian manifold without boundary (also compact in certain contexts)\\
		$\R^m$ and $\T^m$ & $m$-dimensional Euclidean space and torus\\
		$TM/T^*M$ & tangent/cotangent bundle of $M$\\
		$\pi$ & the canonical projection from $T^*M$ to $M$\\
		$\SCL(M)$ & the set of semiconcave functions with linear modulus on $M$\\
		$\SCL_{\mathrm{loc}}(\R^m)$ & the set of locally semiconcave functions with linear modulus on $\R^m$\\
		$\Lip(\R^m)$ & the set of Lipschitz functions on $\R^m$\\
		$\mathrm{AC}$ & absolutely continuous\\
		$\eta=\eta(s)$ & the absolutely continuous deterministic curve on $M$ or $\R^m$\\
		$\phi$ & a function on $M$ or $\R^m$ (usually denoted as a semiconcave function)\\
		$u$ & the viscosity solution of a stationary type HJ equation on $M$ or $\R^m$\\
		$D^{+}\phi(x)$ & the superdifferential of a function $\phi$\\
		$\mathbf{p}^{\#}_{\phi}(x)$ & the minimal-energy element of $D^{+}\phi(x)$ w.r.t. the Hamiltonian $H$\\
		$\text{Sing}\,(\phi)$ & the set of non-differentiability points of a function $\phi$\\
		$\Cut(u)$ & the cut locus of viscosity solution $u$\\
		$\tau_{u}(x)$ & the cut time function of  viscosity solution $u$\\
		$B(x,R)$ & open ball with radius $R$ \\
		$C^{1,1}(M)$ & the space of differentiable functions with Lipschitz differential on $M$ \\
		$\UC(M)$ & the space of uniformly continuous real-valued functions on $M$ \\
		$\Phi^t_H$ & the Hamiltonian flow associated to the Hamiltonian $H$\\
		$X,Y$ & Polish space, i.e., complete and separable metric space\\
		$\mathscr P(X)$ & the set of all Borel probability measures on $X$\\
		$\Gamma(\mu,\nu)$ & the set of transport plans between Borel probability measures $\mu$ and $\nu$\\
		$\mathscr L^m$ & Lebesgue measure on $m$-dimensional manifold\\
		$(\mathscr P_p(X),W_p)$ & $p$-Wasserstein space of $X$ ($p\geqslant 1$)\\	
		$\Gamma_p(\mu,\nu)$ & the set of minimizers of $W_p(\mu,\nu)$\\
		$\pi_i$ & the canonical projection from a product space to its $i$-th component\\
		$\pi_{i,j}$ & the canonical projection from product space to another product space of $i$-th and $j$-th component of the former\\
		$\mathscr P_c(\R^m)$ & the set of all Borel probability measures on $\R^m$ with compact supports\\
		$U$ & the functional on $\mathscr P_c(\R^m)$\\
		$\partial^{\pm}U(\mu)$ & the localized Fr\'echet super/sub-differential of $U$ at $\mu$\\
		$\mathscr S(U)$ & the set of singular points of $U$\\
		$\phi(\cdot)$ & the potential energy functional induced by the function $\phi$ on $M$ or $\R^m$\\
		$C^t(\mu,\nu)$ & the value of dynamical cost functional at $\mu,\nu$\\
		$\Gamma_o^t(\mu,\nu)$ & the set of minimizers of $C^t(\mu,\nu)$\\
		$(\Omega,\mathscr F,\mathbb P)$ & the probability space of sample space $\Omega$\\
		$\E$ & the expectation of random variables on $(\Omega,\mathscr F,\mathbb P)$\\
		$\xi=\xi(s)=\xi(s,\omega)$ & the absolutely continuous random curve of $(\Omega,\mathscr F,\mathbb P)$ on $M$ or $\R^m$\\
		$T_u(\mu)$ & the cut time function of measure associated to potential energy $u(\cdot)$\\ 
		$\mathscr C(u(\cdot))$ & the cut locus of potential energy functional $u(\cdot)$\\
		\bottomrule
\end{longtable}

\end{center}

\section{Preliminaries}

Let $(M,g)$ be a smooth Riemannian manifold. $d_g$ is the metric induced by $g$ on $M$. For any $x\in M$, $T_xM$ and $T_x^*M$ are denoted by the tangent and cotangent spaces of $M$ at $x$, respectively. $TM$ and $T^*M$ are the tangent and cotangent bundle of $M$, respectively. Let $\pi:T^*M\to M$ be the canonical projection from $T^*M$ to $M$, defined as $\pi(x,p):=x$, for all $(x,p)\in T^*M$.

\subsection{Some aspects on semiconcave functions}\label{sec:scl}Some basic relevant notions on semiconcavity are listed below.
\begin{enumerate}[--]
	\item Let $\Omega$ be an open convex subset of $\R^m$. A continuous function $\phi:\Omega\to\R$ is called a \emph{semiconcave function (of linear modulus) with constant $C\geqslant0$} if for any $x\in\Omega$ there exists $p\in\R^m$ such that
	\begin{equation}\label{eq:semiconcave}
	    \phi(y)\leqslant\phi(x)+\langle p,y-x\rangle+\frac C2|x-y|^2,\qquad \forall y\in\Omega,
	\end{equation}
	and we denote as $\phi\in\SCL(\Omega)$. The set of covectors $p$ satisfying \eqref{eq:semiconcave} is called the \emph{proximal superdifferential of $\phi$ at $x$} and we denote it by $D^+\phi(x)$.
	\item Similarly, $\phi$ is \textit{semiconvex} if $-\phi$ is semiconcave. The set of $D^-\phi(x)=-D^+(-\phi)(x)$ is called \emph{the proximal subdifferential of $\phi$ at $x$}. 
	\item For any $\phi\in\SCL(\Omega)$, the set $D^+\phi(x)$ is a singleton if and only if $\phi$ is differentiable at $x$, and $D^+\phi(x)=\{D\phi(x)\}$. A point $x$ is called a \textit{singular point} of a semiconcave function $\phi$ if $D^+\phi(x)$ is not a singleton. We denote by $\text{Sing}\,(\phi)$ the set of all singular points of $\phi$.
	\item Given any $\phi\in\SCL(\Omega)$. We call $p\in D^*\phi(x)$, the set of \emph{reachable gradients}, if there exists a sequence $x_k\to x$ as $k\to\infty$, $\phi$ is differentiable at each $x_k$ and $p=\lim_{k\to\infty}D\phi(x_k)$. We have $D^*\phi(x)\subset D^+\phi(x)$ and $D^+\phi(x)=\text{co}\,D^*\phi(x)$. 
\end{enumerate}
For more on the semiconcave functions in Euclid spaces, the readers can refer to \cite{Cannarsa_Sinestrari_book}.

\begin{Pro}\label{pro:inf}
A function $\phi:\Omega\to\R$ is a semiconcave function with constant $C\geqslant0$ if and only if there exists a family of $C^2$-functions $\{\phi_i\}$ with $D^2\phi_i\leqslant CI_m$ such that
\begin{align*}
	\phi=\inf_i\phi_i.
\end{align*}
\end{Pro}

Proposition \ref{pro:inf} is an important and useful characterization of semiconcavity. Let $S$ be a compact topological space and $F:S\times\R^m\to\R$ be a continuous function such that
\begin{enumerate}[\rm (a)]
	\item $F(s,\cdot)$ is of class $C^2$ for all $s\in S$ and $\|F(s,\cdot)\|_{C^2}$ is uniformly bounded by some constant $C$,
	\item $D_xF(s,x)$ is continuous on $S\times\R^m$,
	\item $\phi(x)=\inf\{F(s,x): s\in S\}$.
\end{enumerate}
We call such a function $\phi$ a \emph{marginal function} of a family of $C^2$-functions.

\begin{Pro}\label{pro:marginal}
Let $\phi(x)=\inf\{F(s,x): s\in S\}$ be a marginal function of the family of $C^2$-functions $\{F(s,\cdot)\}_{s\in S}$.
\begin{enumerate}[\rm (1)]
	\item $\phi$ is semiconcave with constant $C$.
	\item $D^*\phi(x)\subset\{D_xF(s,x): x\in M(x)\}$, where $M(x):=\arg\min\{F(s,x): s\in S\}$. 
	\item For each $x\in\R^m$, 
	\begin{align*}
	    D^+\phi(x)=
		\begin{cases}
			D_xF(s,x),& M(x)=\{s\}\ \text{is a singleton};\\
			\mbox{\rm co}\,\{D_xF(s,x): s\in M(x)\},&\text{otherwise.}
		\end{cases}
	\end{align*}
\end{enumerate} 	
\end{Pro}

\begin{Rem}
A function $\phi:M\to\R$, with $M$ a manifold, is called a semiconcave function, if there exists a family of $C^2$-functions $\{\phi_i\}$ such that $\phi=\inf_i\phi_i$, and the Hessians of $\phi_i$'s is uniformly bounded above. The readers can refer to \cite{Bangert1979,Fathi_Figalli2010,Ohta2009} for more  discussion on the semiconcavity of functions on manifolds, including the equivalence of various definitions of semiconcave functions. Thus, because of the local nature of our main discussion and for convenience, we will work on Euclidean space instead, then get similar conclusions on manifolds. 
\end{Rem}

\subsection{Weak KAM theory}

Let $M$ be a smooth manifold without boundary. A function $L=L(x,v):TM\to\R$ is a time-independent \emph{Tonelli Lagrangian}, if $L$ is of class $C^2$ and satisfies the following conditions:
 \begin{enumerate}[(L1)]
	\item The function $v\mapsto L(x,v)$ is strictly convex for all $x\in M$.
	\item There exist convex, nondecreasing and superlinear functions $\theta_0,\theta_1:[0,+\infty)\to[0,+\infty)$ and positive constants $c_0,c_1$, such that
	\begin{align*}
		\theta_0(|v|_x)-c_0\leqslant L(x,v)\leqslant \theta_1(|v|_x)+c_1,\qquad \forall (x,v)\in TM,
	\end{align*}where $|v|_x=\sqrt{g_x(v,v)}$.	
\end{enumerate}
We call $H$ associated with the Lagrangian $L$ a \emph{Tonelli Hamiltonian}, defined by
\begin{align*}
	H(x,p)=\sup_{v\in T_xM}\{p(v)-L(x,v)\},\qquad (x,p)\in T^*M.
\end{align*}

Given any $x,y\in M$, and $t>0$, we denote by $\Gamma_{x,y}^t$ the set of absolutely continuous curves $\eta\in \text{AC}([0,t],M)$ with $\eta(0)=x$ and $\eta(t)=y$. We call the function $A_t(x,y):\R^+\times M\times M\to\R$ the \emph{fundamental solution} of the associated Hamilton-Jacobi equation, defined by
\begin{align*}
	A_t(x,y):=\inf_{\eta\in\Gamma_{x,y}^t}\int_0^tL(\eta(s),\dot\eta(s))\,ds,\qquad t>0,\quad x,y\in M.
\end{align*}
For any $\phi:M\to\R$, $x\in M$, and $t>0$, we define
\begin{align}
	T_t^-\phi(x)&:=\inf_{y\in M}\{\phi(y)+A_t(y,x)\},\label{inf-conv}\\
	T_t^+\phi(x)&:=\sup_{y\in M}\{\phi(y)-A_t(x,y)\}\label{sup-conv}.	
\end{align}
	
The families $\{T_t^-\}_{t>0}$ and $\{T_t^+\}_{t>0}$ of operators are called \emph{negative and positive Lax-Oleinik evolution} respectively. 
\begin{The}[\cite{Fathi2012}]\label{thm:hje}
	Let $H: T^*M\to\R$ be a Tonelli Hamiltonian, then $u(t,x):=T_t^-\phi(x)$ and $\breve u(t,x)=T_t^+\phi(x)$ are the unique viscosity solution to Cauchy problem of evolutionary type Hamilton-Jacobi equation
\begin{align}\tag{HJ$_{e}-$}\label{hje-}
	\left\{
\begin{array}{l}
		\partial_tu(t,x)+H(x,\partial_xu(t,x))=0,\\
		u(0,x)=\phi(x).
\end{array}
\right.
	\end{align}
	and 
\begin{align}\tag{HJ$_{e}+$}\label{hje+}
\left\{
\begin{array}{l}
		\partial_t\breve u(t,x)-H(x,\partial_x\breve u(t,x))=0\\
		\breve u(0,x)=\phi(x),
\end{array}
\right.
\end{align}
respectively.
\end{The}

\begin{The}[\cite{Fathi_book,Fathi_Maderna2007,Cannarsa_Cheng_Hong2025}]\label{thm:hjs}
Assume that $M$ is compact and $H: T^*M\to\R$ is a Tonelli Hamiltonian, then there exists a unique $c=c[0]\in\R$ such that
\begin{equation}\tag{HJ$_{s}$}\label{hjs}
	H(x,Du(x))=c,\qquad x\in M
\end{equation}
admits a viscosity solution, where $c[0]$ is the Ma\~n\'e critical value. Moreover, the following statements are equivalent:
\begin{enumerate}[\rm (1)]
	\item $u$ is the viscosity solution to \eqref{hjs};
	\item $u=T_t^-u+c[0]t$ for any $t\geqslant 0$. We also call $u$ a weak KAM solution of \eqref{hjs};
	\item $T_t^-\circ T_t^+u=u$ for any $t\geqslant 0$.
\end{enumerate}

If $M=\R^m$, then there exists $c[0]\in\R$ such that \eqref{hjs} admits a viscosity solution $u\in\Lip(\R^m)\cap\SCL(\R^m)$ satisfying that $T_t^-u+c[0]t=u$ for all $t\geqslant 0$. If $c<c[0]$, there is no viscosity solution of \eqref{hjs}.
\end{The}
	
%Combining with the marginal function mentioned in Section \ref{sec:scl}, we know that the viscosity solutions of both \eqref{hje-} and \eqref{hjs} are semiconcave functions of $M$. In noncompact case, the existence of weak KAM solutions is due to \cite{Fathi_Maderna2007}.
%
%\begin{The}[]\label{thm:hjs-exist-Rm}
%Assume that $M=\R^m$. Then, there exists $c[0]\in\R$ such that \eqref{hjs} admits a viscosity solution $u\in\Lip(\R^m)\cap\SCL(\R^m)$ satisfying that $T_t^-u+c[0]t=u$ for all $t\geqslant 0$. If $c<c[0]$, there is no viscosity solution of \eqref{hjs}.
%\end{The}

For more conclusions on the weak KAM theory and viscosity solutions of Hamilton-Jacobi equations in Euclid space $\R^m$ or other non-compact manifold, readers can refer to \cite{Fathi_Maderna2007,Fathi2024}. We collect some regularity properties of Lax-Oleinik operators as follows.

\begin{Pro}[\protect{\cite[Proposition 4.6.6]{Fathi_book}}]\label{pro:lax-oleinik}Suppose $M$ is compact, $\phi\in C(M)$, $t_0>0$ and $c[0]=0$.
	%proposition 4.6.6
	\begin{enumerate}[\rm (1)]
		\item $\{T_t^-\}_{t>0}$ and $\{T_t^+\}_{t>0}$ are the semigroups with respect to $t\in(0,+\infty)$; 
		\item $T_t^-,T_t^+:C(M)\to C(M)$; 
		\item $\lim_{t\to 0^+}T_t^-\phi=\lim_{t\to 0^+}T_t^+\phi=\phi$. Then we define $T_0^-\phi=T_0^+\phi=\phi$;
		\item For any $t>0$, $T_t^-\phi\in\SCL(M)$ (resp. $-T_t^+\phi\in\SCL(M)$) and $\{T_t^-\phi\}_{t\geqslant t_0}$ (resp. $\{-T_t^+\phi\}_{t\geqslant t_0}$) are semiconcave uniformly to $t$; 
		\item $t\mapsto T_t^-\phi$ ($t\mapsto T_t^+\phi$) is uniformly continuous on $[0,+\infty)$;
		\item $(t,x)\mapsto T_t^-\phi(x)$ ($(t,x)\mapsto T_t^+\phi(x)$) is continuous on $[0,+\infty)\times M$, locally Lipschitz on $(0,+\infty)\times M$ and equi-Lipschitz on $[t_0,+\infty)\times M$ with respect to $\phi$.
	\end{enumerate}
\end{Pro}

\subsection{Cut locus}\label{sub:cut locus}

This subsection is inspired by \cite{Cannarsa_Cheng_Hong2025}. For any weak KAM solution \( u \) of equation \eqref{hjs}, we define the function \( B_u: [0, +\infty) \times M \to \mathbb{R} \) as follows:
\begin{align*}
    B_u(t, x) := u(x) - T_t^+u(x)+c[0]t, \quad (t, x) \in [0, +\infty) \times M.
\end{align*}
According to the definition of weak KAM solution, it is straightforward to verify that
\begin{align*}
    B_u(t, x) = T_t^- \circ T_t^+ u(x) - T_t^+ \circ T_t^- u(x).
\end{align*}
We denote by \( \Cut(u) \) the \textit{cut locus of \( u \)}, which consists of points \( x \in M \) where any calibrated curve ending at \( x \) cannot be extended further as a calibrated curve. More precisely, let \( \tau: M \to \mathbb{R} \) be the \emph{cut time function of \( u \)}, defined for any \( x \in M \) as
\begin{align*}
    \tau_u(x) := \sup\{ t \geqslant 0 : \exists \, \eta \in C^1([0, t]; M), u(\eta(t)) - u(x) = A_t(x, \eta(t)) + c[0]t \}.
\end{align*}
Thus \(\Cut(u)=\{x\in M: \tau_u(x)=0\}\).

\begin{Pro}\label{pro: repre-cali}
Suppose $\tau_u(x)>0$ and $t\in(0,\tau_u(x))$. Then there exists a unique $(u,L,c[0])$-calibrated curve $\eta_x: [0,t]\to M$ with $\eta_x(0)=x$, and we have
\begin{align}\label{eq:u-cut-cali}
	\eta_x(s)=\pi\circ\Phi_H^s(x,Du(x))=\pi\circ\Phi_H^{s-t}(\eta_x(t),Du(\eta_x(t))),\,\,s\in[0,t].
\end{align}  
\end{Pro}

Recall that the \textit{(projected) Aubry set} of a weak KAM solution $u$ of \eqref{hjs}
\begin{align*}
	\mathcal I(u):=\{x\in M:\exists\,u\text{-calibrated curve } \eta:[0,+\infty)\to M,\,\eta(0)=x\}.
\end{align*}

\begin{Pro}[\cite{Cannarsa_Cheng_Fathi2017,Cannarsa_Cheng_Fathi2021}]\label{pro:cut time}
Suppose $u$ is the weak KAM solution of \eqref{hjs}.
\begin{enumerate}[\rm (1)]
	\item Given $t>0$ and $x\in M$, then $T_t^+u(x)-c[0]t=u(x)$ if and only if there exists a $u$-calibrated curve $\eta:[0,t]\to M$ such that $\eta(0)=x$;
		\item $\tau_u(x)=\sup\{t\geqslant 0:B_u(t,x)=0\}$. If $t\in[0,\tau_u(x)]$, then $B_u(t,x)=0$;
		\item $\tau_u$ is upper semicontinuous and $\Cut(u)$ is a $G_\delta$-set;
		\item $\mathcal I(u)=\{x\in M:\tau_u(x)=+\infty\}$.
\end{enumerate}
\end{Pro}

\subsection{A priori estimates on fundamental solutions}\label{sub:funda-esti}

The following estimates about $A_t(x,y)$ are well-known.

\begin{Pro}[\protect{\cite[Section 2.4]{CCJWY2020}}]\label{pro:lip-esti}
	Any least action curve of $A_t(x,y)$ is class of $C^2$. Moreover, there exists a monotone nondecreasing, superlinear and convex function of $F:[0,+\infty)\to[0,+\infty)$, such that for any $t,R>0$, when $d(x,y)\leqslant R$, the minimizer $\eta$ of $A_t(x,y)$ satisfies 
	\begin{align*}
		|\dot\eta(s)|\leqslant F\left(\frac{R}{t}\right),\,\,\,\,\forall s\in[0,t].
	\end{align*}
\end{Pro}

\begin{defn}[Lipschitz in the large]
	Let $\kappa_1\geqslant 0,\kappa_2>0$. We say $\phi$ is a \textit{$(\kappa_1,\kappa_2)$-Lipschitz in the large} function in metric space $(X,d)$, if for any $x,y\in X$,
	\begin{align*}
		|\phi(x)-\phi(y)|\leqslant \kappa_1d(x,y)+\kappa_2.
	\end{align*}
\end{defn}

Recall that if $(X,d)=(M,g)$, $\phi\in\UC(M)$ if and only if for arbitrary $\varepsilon>0$, there exists $K_\varepsilon>0$, such that $\phi$ is $(K_\varepsilon,\varepsilon)$-Lipschitz in the large.

%\begin{Rem}[\protect{\cite[Example 3.1]{CCJWY2020}}]\
%	\begin{itemize}[--]
%		\item When $\kappa_2=0$, $\phi\in\Lip(X,\R)$, the Lipschitz real-valued function in $X$;
%		\item Suppose $(X,d)=(M,g)$, $\phi\in\UC(M,\R)$, uniformly continuous real-valued function in $M$, if and only if for arbitrary $\varepsilon>0$, there exists $K_\varepsilon>0$, such that $\phi$ is $(K_\varepsilon,\varepsilon)$- in the large;
%		\item If $M$ is compact, then $\phi$ is bounded if and only if $\phi$ is Lipschitz in the large.
%	\end{itemize}
%\end{Rem}

\begin{Pro}[\protect{\cite[Lemma 3.1]{CCJWY2020}}]\label{pro:funda-minimizer}
	Assume that $\phi$ is lower semicontinuous and $(\kappa_1,\kappa_2)$-Lipschitz in the large. Then for any $x\in M$, there exists $y\in M$, $T_t^-\phi(x)=\phi(y)+A_t(y,x)$. Furthermore, there exists $\lambda_\phi>0$, which only depends on $L$ and $\kappa_1>0$ satisfies that for such minimizer $y_{t,x}$, 
	\begin{align*}
		d(y_{t,x},x)\leqslant \lambda_\phi t+\kappa_2.
	\end{align*}
	More specifically, $\lambda_\phi:=c_0+\theta_1(0)+\theta_0^*(\kappa_1+1)+c_1
$, where $\theta_0^*$ is the Fenchel-Legendre transformation of $\theta_0$.  

Similarly, If $\phi$ is upper semicontinuous and $(\kappa_1,\kappa_2)$-Lipschitz in the large, then for any $x\in M$, there exists $z\in M$, $T_t^+\phi(x)=\phi(z)-A_t(x,z)$ and such maximizer $z_{t,x}$ safisfies
\begin{align*}
	d(z_{t,x},x)\leqslant \lambda_\phi t+\kappa_2.
\end{align*}
\end{Pro}
\begin{Pro}[\cite{Chen_Cheng_Zhang2018,Cannarsa_Cheng3}]\label{pro:cvpde}
Assume that $K$ is a compact and convex subset of $\R^m$ and $\lambda>0$. Then there exist $\tau_1\geqslant \tau_2>0$ and $C_{\lambda},C_{\lambda}^{'},C_{\lambda}^{''}>0$, such that $x\in K$\footnote{In general, $\tau_1$ and $\tau_2$ depend on $x$ in non-compact case.},
\begin{enumerate}[\rm (1)]
	\item  Functions $(t,y)\mapsto A_t(x,y)$ and $(t,y)\mapsto A_t(y,x)$ are both semiconcave and semiconvex (thus $C_{\mathrm{loc}}^{1,1}$) on 
	  \begin{align*}
	  	S(x,\lambda,\tau_1):=\{(t,y)\in\R^+\times K: t\in(0,\tau_1],d(y,x)\leqslant\lambda t\}
	  \end{align*}with constant $C_{\lambda}/t$ and $-C_{\lambda}^\prime/t$ respectively; In this case, for any $(t,y)\in S(x,\lambda_\phi,\tau_1)$,
	  \begin{align*}
	  	D_yA_t(x,y)=p(t),\,\,D_tA_t(x,y)=-H(\eta(t),p(t)),\,\,D_xA_t(x,y)=-p(0),
	  \end{align*}where $\eta\in\Gamma_{x,y}^t$ is the unique minimizer of $A_t(x,y)$, and $p(s):=L_v(\eta(s),\dot\eta(s))$, $s\in[0,t]$ is the dual arc of $\eta$;
	\item  When $t\in(0,\tau_2]$, $A_t(x,\cdot)$ is strictly convex on $B(x,\lambda t)$ with constant $C_{\lambda}^{''}/t$.
\end{enumerate}
	
In particlular, if $\phi\in\SCL(K)$ and $\lambda=\lambda_\phi$, then we can find $\tau_1(\phi)=\tau_1$, $\tau_2(\phi)=\tau_2$ and $C_{\lambda_\phi},C_{\lambda_\phi}^{'},C_{\lambda_\phi}^{''}>0$ depend on $\lambda_\phi$, which make (1) and (2) hold true, as well as $\phi(\cdot)-A_t(x,\cdot)$ is strictly convex on $B(x,\lambda_\phi t)$, thus has a unique maximizer.
\end{Pro}

\subsection{Lasry-Lions regularization, Arnaud's theorem and intrinsic singular characteristic}
In the case of compact manifold $M$, the following Lasry-Lions regularization type result in the context of weak KAM theory is due to Patrick Bernard.

\begin{Pro}[\cite{Bernard2007}]\label{pro:bernard}
Suppose $M$ is compact and $\phi\in\SCL(M)$. Then there exists $0<\tau_3(\phi)\leqslant \tau_2(\phi)$, such that when $t\in(0,\tau_3(\phi)]$, $T_t^+\phi\in C^{1,1}(M)$.
\end{Pro}

Since $\phi$ is semiconcave, $D^+\phi(x)\neq\varnothing$ holds for all $x\in M$. Now we recall the evolution of the pseudo-graphs 
\begin{align*}
	\text{graph}\, (D^+\phi)&:=\{(x,p)\in T^*M:p\in D^+\phi(x),x\in M\},\\
	\text{graph}\, (DT_t^+\phi)&:=\{(x,DT_t^+\phi(x)):x\in M\},\,\,\,\,t\in(0,\tau_3(\phi)].
\end{align*}
under the Hamiltonian flow $\{\Phi_H^t\}_{t\in(0,\tau_3(\phi)]}$.

\begin{Pro}[\cite{Arnaud2011}]\label{pro:arnaud}
For any $t\in(0,\tau_3(\phi)]$, we have 
\begin{align*}
	\Phi_H^t(\emph{graph}\, (DT_t^+\phi))=\emph{graph}\, (D^+\phi).
\end{align*}
\end{Pro}

Suppose $M$ is compact, for any fixed $x\in M$, we define curve $\mathbf y_x(t):[0,\tau_3(\phi)]\to M$,
\begin{align*}
	\mathbf y_x(t):=\left\{
	\begin{array}{ll}
	x,&t=0,\\
	\arg\max\{\phi(y)-A_t(x,y):y\in M\}, &t\in(0,\tau_3(\phi)].
	\end{array}
	\right.
\end{align*}
We call $\mathbf y_x$ an \textit{intrinsic singular characteristic} from $x$.

\begin{Pro}[\cite{CCHW2024}]\label{pro:intrinsic}
Suppose $\phi\in\SCL(M)$. Then 
\begin{align*}
	\mathbf y_x(t)=\pi\circ \Phi_H^t(x,DT_t^+\phi(x)),\,\,x\in M, t\in(0,\tau_3(\phi)].
\end{align*}
Let $\eta_{t,x}\in\Gamma_{x,\mathbf y_x(t)}^t$ be the unique minimizer of $A_t(x,\mathbf y_x(t))$, then
\begin{align*}
	\eta_{t,x}(s)=\pi\circ \Phi_H^s(x,DT_t^+\phi(x))=\pi(\eta_{t,x}(s),DT_{t-s}^+\phi(\eta_{t,x}(s))),\,\,\forall s\in[0,t).
\end{align*}
In local coordinates, $\eta_{t,x}$ satisfies differential equation
\begin{align*}
	\left\{
	\begin{array}{ll}
	  \dot\eta_{t,x}(s)=H_p(\eta_{t,x}(s),DT_{t-s}^+\phi(\eta_{t,x}(s))),&s\in[0,t],\\
	  \eta_{t,x}(0)=x.
	\end{array}
	\right.
\end{align*}Moreover, there exists a constant $C>0$, which is independent to $x$, such that 
\begin{align*}
	\|\mathbf y_x-\eta_{t,x}\|_\infty\leqslant Ct.
\end{align*}
\end{Pro}

In \cite{Cannarsa_Cheng3}, authors show that if $\phi$ is the viscosity solution to \eqref{hjs}, then $x\in\Sing(\phi)$ implies $\mathbf y_x(t)\in\Sing(\phi)$ for $t\in [0,\tau_3(\phi)]$. This is also true for $\phi\in\SCL(M)$ only with
\begin{align*}
	\mathop{\arg\min}_{p\in D^+\phi(x)}\{H(x,p)\}\notin D^*\phi(x),
\end{align*}
where $D^*\phi(x)$ is the set of reachable gradients of $\phi$ at $x$.

\subsection{Optimal transport}
Suppose $(X,d)$ is a metric space. We denote $\mathscr P(X)$ the set of all probability measures on $X$.
 
\subsubsection{Monge-Kantorovich problem and Kantorovich duality}

For two fixed Polish spaces $X$ and $Y$, the modern formulation of Monge's problem has probability measures $\mu\in \mathscr P(X)$, $\nu\in\mathscr P(Y)$ and a Borel cost function $c(x,y): X\times Y\to \R$, representing the cost of shipping a unit mass from $x$ to $y$. Given these data, the problem is
\begin{equation}\label{eq:Monge}\tag{M}
	\inf\left\{\int_Xc(x,T(x))\,d\mu:T:X\to Y\text{ is Borel}, T_\#\mu=\nu\right\},
\end{equation}
where $T$ is called the transport map from $\mu$ to $\nu$. Under certain conditions (for example, \cite[Theorem 10.28]{Villani_book2009}), Monge's problem \eqref{eq:Monge} can have a solution. Brenier (\cite{Brenier1987,Brenier1991}) and Knott-Smith (\cite{Knott_Smith1984}) proved that \eqref{eq:Monge} has a solution when $X=Y=\R^m$, $c(x,y)=\frac{1}{2}|x-y|^2$ and $\mu\ll\mathscr L^m$. 

\begin{Rem}\label{rem:find the law}
We remark that the conclusion due to Brenier and Knott-Smith provides us with the existence of a probability space $(\Omega,\mathscr F,\mathbb P)$ such that for arbitrary $\mu\in\mathscr P_2(\R^m)$ (see Definition \ref{defn:wasserstein}), we can always find a random variable $T:\Omega\to\R^m$, which satisfies $\law(T)=\mu$. Indeed, we choose $\Omega=\R^m$, $\mathscr F$ is the set of Lebesgue measurable subsets of $\R^m$ and $\mathbb P=\mathscr L^m|_{B(0,R)}$, where $R>0$ makes $\mathscr L^m(B(0,R))=1$. If $c(x,y)=\frac{1}{2}|x-y|^2$, for the corresponding Monge's transport problem, there exists an optimal map $T:\R^m\to\R^m$, i.e., $\law(T)=\mu$. By Example 10.36 of \cite{Villani_book2009}, the argument for the case $X=Y=M$ with $M$ compact manifold and $\mu\in\mathscr P(M)$ is similar.
\end{Rem}

For Polish spaces $X$ and $Y$, we call $\gamma\in\mathscr P(X\times Y)$ the \textit{transport plan} between $\mu\in\mathscr P(X)$ and $\nu\in\mathscr P(Y)$, if 
\begin{align*}
	(\pi_X)_{\#}\gamma=\mu,\qquad (\pi_Y)_{\#}\gamma=\nu,
\end{align*}
where $\pi_X$ and $\pi_Y$ are the canonical projection maps from $X\times Y$ to $X$ and $Y$ respectively. We denote $\Gamma(\mu,\nu)$ the set of all transport plans between $\mu$ and $\nu$. Obviously, $\mu\times\nu\in\Gamma(\mu,\nu)$.

Kantorovich's problem asks to find 
\begin{equation}\label{eq:Kantorovich}\tag{K}
	C(\mu,\nu):=\inf\left\{\int_{X\times Y}c(x,y)\ d\gamma: \gamma\in\Gamma(\mu,\nu)\right\}.
	%=\inf_{\substack{\law(x)=\mu\\ \law(y)=\nu}}\mathbb{E}(x,y),
\end{equation}
%where $x(\omega)$ and $y(\omega)$ in the last term of \eqref{eq:Kantorovich} are the $X$-value and $Y$-value random variables of $(\Omega,\mathscr F,\mathbb P)$ respectively. 

The following is one of the most basic results of the theory of optimal transport, namely the Kantorovich-Rubinstein duality theorem.

\begin{The}[\cite{Villani_book2009,Ambrosio_Brue_Semola_book2021}]
If $c: X\times Y\to\R$ is lower-semicontinuous, then
\begin{equation}\label{eq:KB}\tag{KR}
	C(\mu,\nu)=\sup_{(\phi,\psi)\in I_c}\left\{\int_Y\psi\ d\nu-\int_X\phi\ d\mu\right\}=\sup_{(\phi,\psi)\in K_c}\left\{\int_Y\psi\ d\nu-\int_X\phi\ d\mu\right\},
\end{equation}
where 
\begin{align*}
	I_c:=&\,\{(\phi,\psi): \psi(y)-\phi(x)\leqslant c(x,y)\ \forall\ x\in X, y\in Y\},\\
	K_c:=&\,\left\{(\phi,\psi)\in I_c: \psi=\inf_{x\in X}\{\phi(x)+c(x,\cdot)\}, \phi=\sup_{y\in Y}\{\psi(y)-c(\cdot,y)\}\right\}.
\end{align*}
\end{The}

In many literatures, elements in $K_c$ are referred to as \textit{Kantorovich admissible pairs}. 

\begin{defn}[\cite{Cheng_Hong_Shi2024}]\label{def:abstract lax-oleinik}
For any $\phi:X\to\R$, $\psi:Y\to\R$, the \textit{positive and negative type of abstract Lax-Oleinik operators} are defined as
\begin{align*}
	T^+\psi(x):=\sup_{y\in Y}\{\psi(y)-c(x,y)\},\,\,\,\,T^-\phi(y):=\inf_{x\in X}\{\phi(x)+c(x,y)\}.
\end{align*}
\end{defn}

Let $\phi:X\to\R$ and $\psi:Y\to\R$. For any $\mu\in\mathscr{P}(X)$ and $\nu\in\mathscr{P}(Y)$, we try to find a function $\phi:X\to\R$ and a function $\psi:Y\to\R$ such that
\begin{align}
	\int_YT^-\phi\ d\nu=&\,\inf_{\gamma\in\Gamma(\mu,\nu)}\int_{X\times Y}\phi(x)+c(x,y)\ d\gamma,\label{eq:K-}\tag{K$^-$}\\
	\int_XT^+\psi\ d\mu=&\,\sup_{\gamma\in\Gamma(\mu,\nu)}\int_{X\times Y}\psi(y)-c(x,y)\ d\gamma.\label{eq:K+}\tag{K$^+$}
\end{align}
In terms of abstract Lax-Oleinik operators, the equivalence of the problems \eqref{eq:K-} and \eqref{eq:K+} and well known Rubinstein-Kantorovich duality can be showed as follows.

\begin{The}[\cite{Cheng_Hong_Shi2024}]\label{thm:K+-}
\hfill
\begin{enumerate}[\rm (1)]
	\item If $\phi:X\to\R$ is a solution of \eqref{eq:K-}, then $(\phi,T^-\phi)\in I_c$ is a solution of \eqref{eq:KB}. Conversely, if $(\phi,\psi)\in I_c$ is a solution of \eqref{eq:KB}, then $\phi$ solves \eqref{eq:K-}.
	\item If $\psi:Y\to\R$ is a solution of \eqref{eq:K+}, then $(T^+\psi,\psi)\in I_c$ is a solution of \eqref{eq:KB}. Conversely, if $(\phi,\psi)\in I_c$ is a solution of \eqref{eq:KB}, then $\psi$ solves \eqref{eq:K+}.
\end{enumerate}	

\end{The}

Invoking the equivalence above, it is useful to refine the analysis of the dynamical nature of the associated optimal transport problem, as initiated by Bernard and Buffoni (\cite{Bernard_Buffoni2006, Bernard_Buffoni2007a, Bernard_Buffoni2007b}), which also forms the basis of the main content in this paper.

\subsubsection{Wasserstein space}
\begin{defn}[Wasserstein space]\label{defn:wasserstein}
Let $X$ be a Polish space. For $p\in[1,+\infty)$, set
\begin{align*}
	\mathscr P_p(X):=\left\{\mu\in\mathscr P(X): \int_Xd^p(x,x_0)\,d\mu<+\infty\right\},\,\,\,\,\exists x_0\in X,
\end{align*}
We call
\begin{align*}
	W_p(\mu,\nu):=\left\{\inf_{\gamma\in\Gamma(\mu,\nu)}\int_{X^2}d^p(x,y)\,d\gamma\right\}^{\frac{1}{p}}.
\end{align*}
For all $p\geqslant 1$, we denote $\Gamma_p(\mu,\nu)$ is the set of all optimal plans of $W_p(\mu,\nu)$. The set $\mathscr P_p(X)$ endowed the \emph{Wasserstein distance} $W_p$ is called \emph{$p$-Wasserstein space}.
\end{defn}

\begin{Rem}
If $X$ is a Polish space, then so is $\mathscr P_p(X)$. If $X$ is compact, $\mathscr P_p(X)$ is tight in the sense of Prokhorov. However, when $X$ is only locally compact, $(\mathscr P_p(X),W_p)$ may not be. See, for instance, \cite{Villani_book2009,Ambrosio_GigliNicola_Savare_book2008}.
\end{Rem}

Given $\boldsymbol{\mu}$-measurable maps $\boldsymbol{r},\boldsymbol{s} :(\boldsymbol{X},\boldsymbol{\mu})\to X$, a very useful inequality giving an estimate from above of the Wasserstein distance is
\begin{equation}\label{eq:Wass_dist_push_forw}
	W_p(\boldsymbol{r}_\#\boldsymbol{\mu},\boldsymbol{s}_\#\boldsymbol{\mu})\leqslant\boldsymbol{d}(\boldsymbol{r},\boldsymbol{s})_{L^p(\boldsymbol{\mu},X)}=\bigg\{\int d^p(x_1,x_2)d\gamma\bigg\}^{1/p}
\end{equation}
where $\gamma=(\boldsymbol{r},\boldsymbol{s})_\#\boldsymbol{\mu}\in\Gamma(\boldsymbol{r}_\#\boldsymbol{\mu},\boldsymbol{s}_\#\boldsymbol{\mu})$.

A \textit{constant speed geodesic} of a metric space $(X,d)$ is a curve $\eta:[0,1]\to X$, which satisfies that for any $0\leqslant s\leqslant t\leqslant 1$,
\begin{equation}\label{eq:constant_speed}
	d(\eta(s),\eta(t))=(t-s)d(\eta(0),\eta(1)).
\end{equation}
By triangle inequality, if the equality is replaced by ``$\leqslant$'' in \eqref{eq:constant_speed} for any $0\leqslant s\leqslant t\leqslant 1$, then $\eta$ is a constant speed geodesic equivalently. 

For $\mu_i\in\mathscr P(X_i)$, $i=1,\cdots,N$, we call $\boldsymbol{\mu}\in\Gamma(\mu_1,\mu_2,\cdots,\mu_N)$, if $(\pi_i)_\#\boldsymbol{\mu}=\mu_i$, where $\pi_i$ is the canonical projection from $X_1\times \cdots\times X_N$ to $X_i$. In the case $X_i=X$, a Hilbert space, $i=1,\cdots,N.$ and $\boldsymbol{\mu}\in\mathscr P(X^N)$, $N\geqslant 2$, $1\leqslant i,j,k\leqslant N$ and $\lambda\in [0,1]$, set
\begin{align}
	\pi_{i\to j}^\lambda &:=(1-\lambda)\pi_i+\lambda\pi_j: X^N\to X,\label{eq:pi-ij}\\
	\pi_{i\to j,k}^\lambda &:=(1-\lambda)\pi_{i,k}+\lambda\pi_{j,k}: X^N\to X^2,\label{eq:pi-ijk}\\
	\pi_{k,i\to j}^\lambda &:=(1-\lambda)\pi_{k,i}+\lambda\pi_{k,j}: X^N\to X^2,\label{eq:pi-kij}\\
	\boldsymbol{\mu}_{i\to j}^\lambda &:=(\pi_{i\to j}^\lambda)_{\#}\boldsymbol{\mu}\in\mathscr P(X),\label{eq:mu-ij}\\
	\boldsymbol{\mu}_{i\to j,k}^\lambda &:=(\pi_{i\to j,k}^\lambda)_{\#}\boldsymbol{\mu}\in\mathscr P(X^2),\label{eq:mu-ijk}\\
	\boldsymbol{\mu}_{k,i\to j}^\lambda &:=(\pi_{k,i\to j}^\lambda)_{\#}\boldsymbol{\mu}\in\mathscr P(X^2),\label{eq:mu-kij}
\end{align}
where $\pi_i$ is the canonical projection from $X^N$ to the $i$-th $X$, and $\pi_{i,j}$ is the canonical projection from $X^N$ to the product space of $i$-th and $j$-th $X$.

The following theorem implies that, when $X$ is a Hilbert space, $(\mathscr P_p(X),W_p)$ is a geodesic space, and geodesics are determined by the optimal plans between the two prescribed measures.

\begin{The}[\protect{\cite[Lemma 7.2.1 \& Theorem 7.2.2]{Ambrosio_GigliNicola_Savare_book2008}}]\label{thm:wasser-geo}
	Assume that $(X,\langle\cdot,\cdot\rangle)$ is a Hilbert space, and $p>1$. If $\mu_1,\mu_2\in\mathscr P_p(X)$ and $\boldsymbol{\mu}\in\Gamma_p(\mu_1,\mu_2)$, then the curve of $\mathscr P_p(X)$, $\eta(t):=\boldsymbol{\mu}_{1\to 2}^t$ is a constant speed geodesic connecting $\mu_1$ and $\mu_2$. 
	
	Conversely, for each constant speed geodesic $\gamma:[0,1]\to \mathscr P_p(X)$ connecting $\mu_1$ and $\mu_2$, there exists $\boldsymbol{\mu}\in\Gamma_p(\mu_1,\mu_2)$ such that $\gamma(t)=\boldsymbol{\mu}^t_{1\to 2}$, $\Gamma_p(\mu_1,\gamma(t))=\{\boldsymbol{\mu}_{1,1\to 2}^t\}$ and $\Gamma_p(\gamma(t),\mu_2)=\{\boldsymbol{\mu}_{1\to 2,2}^t\}$, $t\in(0,1)$.
\end{The}

\begin{The}[\protect{\cite[Theorem 7.3.2]{Ambrosio_GigliNicola_Savare_book2008}}]\label{thm:wass2-scl}
For arbitrary $\mu_1,\mu_2,\nu\in\mathscr P_2(X)$ and $\boldsymbol{\mu}_{1\to 2}\in\Gamma_2(\mu_1,\mu_2)$, 
\begin{align*}
	(1-\lambda)W_2(\mu_1,\nu)+\lambda W_2(\mu_2,\nu)-W_2(\boldsymbol{\mu}_{1\to 2}^\lambda,\nu)\leqslant \lambda(1-\lambda)W_2^2(\mu_1,\mu_2),
\end{align*}i.e., $(\mathscr P_2(X),W_2)$ is a positively curved space (see \cite[Definition 3.19]{Ambrosio_Gigli2013}). 
\end{The}
It can be shown that when $(X,d)$ is a positively curved space, so is $(\mathscr P_2(X),W_2)$; However, it is worth noting that if $(X,d)$ is a non-positively curved space, $(\mathscr P_2(X),W_2)$ is not necessarily non-positively curved (\cite[Example 3.21]{Ambrosio_Gigli2013}). 

\section{Generalized differential and semiconcavity of functionals of $\mathscr P_c(\R^m)$}

In the section, we assume $M=\R^m$ and we work on the space 
\begin{align*}
	\mathscr P_c(\R^m):=\left\{\mu\in\mathscr P(\R^m): \supp(\mu)\ \text{is a compact subset of }\,\R^m\right\}.
\end{align*}
It is obvious that $\mathscr P_c(\R^m)\subset\mathscr{P}_p(\R^m)$ for all $p\geqslant 1$. If $M$ is a compact manifold, $\mathscr P_c(M)=\mathscr P(M)$. It is worth noting that for any $p\geqslant 1$, $(\mathscr P_c(\R^m),W_p)$ as a metric space, is not locally compact. See, for instance, the example in \cite[Remark 7.1.9]{Ambrosio_GigliNicola_Savare_book2008}.

%In classical calculus of variation, 
%
%This also implies the case of Dirac measures in calculus of variations and the probability measures $\mathscr P(M)$ on the compact manifold $M$ (e.g. $m$-dimensional torus: $\T^m=\R^m/\Z^m$).
\subsection{Generalized differential and semiconcavity}

To deal with singularities of the functionals on Wasserstein spaces, we need some notions of $\mathscr{P}_2$-subdifferential calculus in the setting of compact supported measures, developed in \cite[Section 9.1, 9.2 \& 10.3]{Ambrosio_GigliNicola_Savare_book2008} (see also the seminal work \cite{McCann1997}). We follow the setting in \cite{Bonnet_Frankowska2022}. We write a functional $U:\mathscr{P}_c(\R^m)\to\R$ to mean the restriction of an extended real-valued functional $U:\mathscr{P}_2(\R^m)\to\R\cup\{\pm\infty\}$ such that $U$ is finite valued on $\mathscr{P}_c(\R^m)$. Given $\mu\in\mathscr P_c(\R^m)$ and $R>0$, we call $B_R(\mu):=\cup_{x\in\supp(\mu)}B(x,R)$ the $R$-fattening of $\supp(\mu)$.

\begin{defn}[Localized Fr\'{e}chet generalized differential]\label{def:local frechet}
Assume that $U: \mathscr P_c(\R^m)\to\R$, $\mu\in\mathscr P_c(\R^m)$, we say $\alpha\in L^2(\R^m;\mu)$ belongs to $\partial^+U(\mu)=\partial^+_{\rm loc}U(\mu)$, the set of \textit{localized Fr\'{e}chet superdifferential} of $U$ at $\mu$, if for any $R>0$ and any $\nu\in\mathscr{P}(B_R(\mu))$, 
\begin{equation}\label{eq:wass-supdiff}
	U(\nu)-U(\mu) \leqslant \sup_{\gamma \in \Gamma_2(\mu, \nu)}\left\{\int_{\mathbb{R}^{2m}}\langle\alpha(x), y-x\rangle\,d\gamma\right\}+o_R(W_2(\mu,\nu)),
\end{equation}
where $o_R(W_2(\mu,\nu))$, depending on $R$, represents the high order infinitesimal of $W_2(\mu,\nu)$.% in relation to $R$. 

Similarly, $\beta\in L^2(\R^m;\mu)$ belongs to $\partial^-U(\mu)=\partial^-_{\rm loc}U(\mu)$, the set of \textit{localized Fr\'{e}chet subdifferential} of $U$ at $\mu$, if for any $R>0$ and any $\nu\in\mathscr{P}(B_R(\mu))$,
\begin{equation}\label{eq:wass-subdiff}
	U(\nu)-U(\mu) \geqslant \inf_{\gamma \in \Gamma_2(\mu, \nu)}\left\{\int_{\mathbb{R}^{2m}}\langle\beta(x), y-x\rangle\,d\gamma\right\}+o_R(W_2(\mu,\nu)).
\end{equation}

For a time dependent functional $U:\R\times\mathscr P_c(\R^m)\to\R$, we say  $(q,\alpha)\in\R\times L^2(\R^m;\mu)\in \partial^+U(t,\mu)$, if for any $R>0$ and $a,b\in\R$ satisfying $t\in(a,b)$ and $(s,\nu)\in[a,b]\times\mathscr P(B_R(\mu))$,
\begin{align*}
	U(s,\nu)-U(t,\mu)\leqslant\sup_{\gamma \in \Gamma_2(\mu, \nu)}\left\{\int_{\mathbb{R}^{2m}}\langle\alpha(x), y-x\rangle\,d\gamma\right\}+q(s-t)+o_R(W_2(\mu,\nu)+|t-s|).
\end{align*}The definition of $\partial^-U(t,\mu)$ is similar.
\end{defn}

\begin{Rem}
We regard two elements in $\partial^\pm U(\mu)$ to be the same if they coincide up to a $\mu$-negligible set.
%
%We use the word ``localized'' in Definition \ref{def:local frechet} to emphasize the approximation approach to fixed $\mu$ is always in $\mathscr{P}_c(\R^m)$, which indicates that the approximation is local, rather than in $\mathscr P_2(\R^m)$, taking into account the case of non-compact supported measures with more approximate ways.
\end{Rem}

\begin{Pro}\label{pro:wass-supdiff-comparision}
Assume functionals $U,V:\mathscr P_c(\R^m)\to\R$ and $U\leqslant V$. If there exists $\mu\in\mathscr P_c(\R^m)$ such that $U(\mu)=V(\mu)$, we have $\partial^+V(\mu)\subset\partial^+U(\mu)$ and $\partial^-U(\mu)\subset\partial^-V(\mu)$.
\end{Pro}

\begin{proof}
For any $\alpha\in\partial^+V(\mu)$, due to \eqref{eq:wass-supdiff} and the relation $U\leqslant V$, we get 
\begin{align*}
	U(\nu)-U(\mu)&\leqslant V(\nu)-V(\mu)\\ 
	&\leqslant \sup_{\gamma \in \Gamma_2(\mu, \nu)}\left\{\int_{\mathbb{R}^{2m}}\langle\alpha(x), y-x\rangle\,d\gamma\right\}+o_R(W_2(\mu,\nu)).
\end{align*}
This implies that $\alpha\in\partial^+U(\mu)$. One can show that $\partial^-U(\mu)\subset\partial^-V(\mu)$ in a similar way.
\end{proof}

\begin{defn}[Singularity of functional in $\mathscr P_c(\R^m)$]\label{def:Pc-diff-sing}
Let $U:\mathscr P_c(\R^m)\to\R$ and $\mu\in\mathscr P_c(\R^m)$. 
\begin{enumerate}[\rm (1)]
	\item $\mu$ is called a \emph{regular point} of $U$ if both $\partial^+ U(\mu)$ and $\partial^-U(\mu)$ are non-empty. We also say $U$ is \emph{locally differentiable} at $\mu$.
	\item $\mu$ is called a \emph{singular point} of $U$ if $U$ is not locally differentiable at $\mu$. We denote by $\mathscr{S}(U)$ the set of all singular points of $U$.
\end{enumerate}
\end{defn}

\begin{Lem}\label{lem:smooth-approx}
Suppose $\mu\in\mathscr P_c(\R^m)$, $f\in C_c^\infty(\R^m)$ and $R>0$. There exists $\delta_{f,R}>0$ such that $\nu_t:=(\mathrm{id}+tDf)_\#\mu\in\mathscr P(B_{R}(\mu))$ and
\begin{align*}
	\Gamma_2(\mu,\nu_t)=\left\{\left[\mathrm{id}\times \left(\mathrm{id}+tDf\right)\right]_\#\mu\right\},\qquad \forall t\in[0,\delta_{f,R}].
\end{align*}
\end{Lem}

\begin{Pro}\label{pro:diff-singleton}
Suppose $U:\mathscr P_c(\R^m)\to\R$ is locally differentiable at $\mu$, then $\partial^+U(\mu)$ and $\partial^-U(\mu)$ are singleton. In this case, we denote
\begin{align*}
	\partial^+U(\mu)=\partial^-U(\mu)=\{\partial U(\mu)\}.
\end{align*}
\end{Pro}

We provide proofs of Lemma \ref{lem:smooth-approx} and Proposition \ref{pro:diff-singleton} in the appendix.

\begin{defn}[\cite{Bonnet_Frankowska2022}, Semiconcavity of functionals in $\mathscr P_c(\R^m)$]\label{def:local-scl}
$U:\mathscr P_c(\R^m)\to\R$ is \textit{locally (geodesically) semiconcave} with constant $C_R$, if for any fixed $R>0$, there exists $C_R>0$, such that for all $\mu_1,\mu_2\in\mathscr P(B(0,R))$ and $\gamma\in\Gamma_2(\mu_1,\mu_2)$ and $t\in[0,1]$,
\begin{align*}
	(1-t)U(\mu_1)+tU(\mu_2)-U(\gamma_{1\to 2}^t)\leqslant C_Rt(1-t)W_2^2(\mu_1,\mu_2),
\end{align*}where $\gamma_{1\to 2}^t$ is given by \eqref{eq:mu-ij}. We say $U:\mathscr P_c(\R^m)\to\R$ is \textit{locally strongly semiconcave}, if for any $\mu_1,\mu_2\in\mathscr P(B(0,R))$ and $\gamma\in\Gamma(\mu_1,\mu_2)$,
\begin{align*}
	(1-t)U(\mu_1)+tU(\mu_2)-U(\gamma_{1\to 2}^t)\leqslant C_Rt(1-t)W_{2,\gamma}^2(\mu_1,\mu_2)
\end{align*}
for each $t\in[0,1]$, where $\gamma_{1\to 2}^t$ is defined in \eqref{eq:mu-ij}, $W_{2,\gamma}^2(\mu_1,\mu_2):=\int_{\R^{2m}}d^2(x,y)\,d\gamma$.
\end{defn}

Theorem \ref{thm:wass2-scl} implies 2-Wasserstein metric $W_2$ is globally semiconcave in $\mathscr P_2(\R^m)$.

\begin{Rem}
For time-dependent functional $U:\R\times\mathscr P_c(\R^m)\to\R$, we only need to restrict to $[a,b]\times B(0,R)$, where $a,b\in\R$ and $R>0$. We can similarly define the semiconcavity of functionals of $\R\times\mathscr P_c(\R^m)$.
\end{Rem}

\subsection{Potential energy functional}\label{subsec:potential-energy-diff-scl}

\begin{defn}[\cite{Ambrosio_GigliNicola_Savare_book2008}, potential energy]\label{def:potential energy}
Given a Borel measurable function $\phi:\R^m\to[-\infty,+\infty]$. A functional of $\mathscr P(\R^m)$ induced by $\phi$
\begin{align*}
	\Phi(\mu):=\int_{\R^m}\phi(x)\,d\mu\qquad (\mu\in\mathscr P(\R^m))
\end{align*}
is called \textit{potential energy functional of} $\phi$. We denote $\Phi(\cdot)$ by $\phi(\cdot)$ for brevity if there is no confusion. 
\end{defn}

The following properties on generalized differential and semiconcavity for potential energy functional are useful.

\begin{Pro}\label{pro:potential-scl}
If $\phi\in\SCL_{\mathrm{loc}}(\R^m)$, then the potential energy functional $\phi(\cdot):\mathscr P_c(\R^m)\to\R$ is locally strongly semiconcave.
\end{Pro}

\begin{proof}
Let $R>0$, $\mu_1,\mu_2\in\mathscr P(B(0,R))$ and $\gamma\in\Gamma(\mu_1,\mu_2)$. Then
\begin{align*}
	&(1-t)\phi(\mu_1)+t\phi(\mu_2)-\phi(\gamma_{1\to 2}^t)\\
		&=\int_{\R^{2m}}(1-t)\phi(x)+t\phi(y)-\phi((1-t)x+ty)d\gamma\\
		&\leqslant\int_{\R^{2m}}C_Rt(1-t)|x-y|^2\,d\mu=C_Rt(1-t)W_{2,\gamma}^2(\mu_1,\mu_2).
\end{align*}
Thus, $\phi(\cdot)$ is locally strongly semiconcave by Definition \ref{def:local-scl}.
\end{proof}

The following proposition is firstly claimed in \cite{Ambrosio_Gangbo2008}, which discussed analogous conclusions on convex analysis in $\mathscr P_2(\R^m)$. We afford a proof in the appendix.

\begin{Pro}\label{pro:potential-scl-equiv}
Suppose $\mu\in\mathscr P_c(\R^m)$, $\phi\in\SCL_{\mathrm{loc}}(\R^m)$ and $\alpha\in L^2(\R^m;\mu)$. Then for the potential energy functional $\phi(\cdot)$, the following statements are equivalent:
\begin{enumerate}[\rm (1)]
	\item $\alpha\in\partial^+\phi(\mu)$;
	\item For arbitrary $R>0$ and $\nu\in\mathscr P(B_R(\mu))$,
	\begin{equation}\label{eq:potential-scl-equiv2}
		\phi(\nu)-\phi(\mu)\leqslant \inf_{\gamma \in \Gamma_2(\mu, \nu)}\left\{\int_{\mathbb{R}^{2m}}\langle\alpha(x), y-x\rangle\, d\gamma\right\}+o_R(W_2(\mu,\nu));
	\end{equation}
	\item For any $R>0$ and $\nu\in\mathscr P(B(0,R))$, there exists $C_R>0$, such that
	\begin{equation}\label{eq:potential-scl-equiv3}
		\phi(\nu)-\phi(\mu)\leqslant \inf_{\gamma\in\Gamma_2(\mu,\nu)}\left\{\int_{\mathbb{R}^{2m}}\langle\alpha(x), y-x\rangle\,d\gamma\right\}+C_RW_2^2(\mu,\nu).
	\end{equation}
\end{enumerate}
\end{Pro}

Generalized differential of functionals in $\mathscr P_p(X)$ was studied in literature \cite{Ambrosio_GigliNicola_Savare_book2008,Ambrosio_Gigli2013}. We restate the relevant conclusions in $\mathscr P_c(\R^m)$ and give a proof for the sake of rigorousness and completeness.

\begin{The}\label{thm:partial+D+}
Suppose $\phi\in\SCL_{\mathrm{loc}}(\R^m)$ and $\mu\in\mathscr P_c(\R^m)$. Then $\alpha\in\partial^+\phi(\mu)$ if and only if $\alpha(x)\in D^+\phi(x)$ holds true for $\mu$-a.e. $x\in\R^m$.
\end{The}

See the appendix for the proof of Theorem \ref{thm:partial+D+}. For the potential energy functional $\phi(\cdot)$ (i.e., $\Phi(\cdot)$) with $\phi$ a semiconcave function, Theorem \ref{thm:partial+D+} clarifies the relation between $\partial^+\phi$ and $D^+\phi$. This is a key point allowing us to study the singularities of potential energy functional $\phi(\cdot)$ by means of the singularities of the semiconcave function $\phi$.

\begin{Cor}\label{cor:partial+singleton}
Let $\phi\in\SCL_{\mathrm{loc}}(\R^m)$ and $\mu\in\mathscr P_c(\R^m)$. Then $\partial^+\phi(\mu)\neq\varnothing$ and $\phi(\cdot)$ is locally differentiable at $\mu$ if and only if $\partial^+\phi(\mu)$ is a singleton. Furthermore, $\alpha=\partial\phi(\mu)$ if and only if $\alpha(x)=D\phi(x)$ holds for $\mu$-a.e. $x\in\R^m$.
\end{Cor}

\begin{proof}
For the first assertion, we note that $D^+\phi(x)\neq\varnothing$ for all $x\in\R^m$. If $\phi(\cdot)$ is locally differentiable at $\mu$, by Proposition \ref{pro:diff-singleton}, $\partial^+\phi(\mu)$ is a singleton. Conversely, if $\partial^+\phi(\mu)$ is a singleton, we can find a $\mu$-full-measure set $E$ such that for $x\in E$, $D^+\phi(x)=\{\alpha(x)\}$. This implies $\phi$ is differentiable at $x$, $D^-\phi(x)=\{\alpha(x)\}$ and $\alpha\in\partial^-\phi(\mu)$. Then Proposition \ref{pro:diff-singleton} shows that $\phi(\cdot)$ is localized differeniable at $\mu$. The rest of proof is obvious. 	
\end{proof}

\begin{The}\label{thm:singular measure}
Let \(\phi \in \SCL_{\mathrm{loc}}(\R^m)\). Then \(\mu \in \mathscr{S}(\phi(\cdot))\), i.e., \(\mu\) is a singular point of \(\phi(\cdot): \mathscr{P}_c(\R^m) \to \mathbb{R}\), if and only if \(\mu(\Sing(\phi)) > 0\), where \(\Sing(\phi)\) denotes the set of singularities of \(\phi\). Consequently, \(\mu \notin \mathscr{S}(\phi(\cdot))\) if \(\mu \ll \mathscr{L}^m\), meaning that \(\mu\) is absolutely continuous with respect to the Lebesgue measure \(\mathscr{L}^m\).
\end{The}

\begin{proof}
If $\mu(\Sing(\phi))=0$, there exists $E$ with $\mu(E)=1$ such that $\phi$ is differentiable at $x$ for all $x\in E$. In other words, $D\phi\in\partial^+\phi(\mu)\cap\partial^-\phi(\mu)$ and $\phi(\cdot)$ is locallly differentiable at $\mu$. If $\mu(\Sing(\phi))>0$, then one can find $\alpha_1,\alpha_2\in \partial^+\phi(\mu)$ which are distinct by removing a $\mu$-null set. It yields $\mu$ is the singular point of $\phi(\cdot)$. 

Finally, recall that $\Sing(\phi)$ is an $(m-1)$-countably rectifiable set (see \cite[Corollary 4.1.13]{Cannarsa_Sinestrari_book}) and $\mathscr L^m(\Sing(\phi))=0$. If $\mu\ll\mathscr L^m$, then $\mu(\Sing(\phi))=0$ and this implies $\mu$ is a regular point of $\phi(\cdot)$.
\end{proof}

We need also consider the time-dependent potential energy functionals. Suppose $\phi:\R\times\R^m\to[-\infty,+\infty]$ is Borel measurable. The potential energy functional of $\mathscr P_c(\R^m)$ induced by $\phi$ is
\begin{align*}
	\Phi(t,\mu):=\int_{\R^m}\phi(t,x)\,d\mu=\int_{\R^{m+1}}\phi(r,x)\,d(\delta_t\times\mu),\,\,\,\,\forall(t,\mu)\in\R\times\mathscr P_c(\R^m).
\end{align*}
Again, we write $\Phi(t,\mu)=\phi(t,\mu)$ for the sake of brevity if without confusion.

\begin{Rem}\label{rem:time-dependent-D+}	
Assume that $\phi\in\SCL_{\mathrm{loc}}(\R\times\R^m)$ and $\mu\in\mathscr P_c(\R^m)$. 
\begin{enumerate}[\rm (1)]
	\item According to Definition \ref{def:local frechet} and Theorem \ref{thm:partial+D+}, $(q,\alpha)\in\partial^+\phi(t,\mu)$ if and only if $(q,\alpha(x))\in D^+\phi(t,x)$ for $\mu$-a.e. $x\in\R^m$;
    \item Similar to Corollary \ref{cor:partial+singleton}, $\partial^+\phi(t,\mu)\neq\varnothing$, and $\phi(\cdot,\cdot)$ is locally differentiable at $(t,\mu)$ if and only if $\partial^+\phi(t,\mu)$ is a singleton, and $(q,\alpha)=\partial\phi(t,\mu)$ if and only if $(q,\alpha(x))=(D_t\phi(t,x),D_x\phi(t,x))$ for $\mu$-a.e. $x\in\R^m$.
\end{enumerate}
\end{Rem}

\subsection{Dynamical cost functionals}

%The following definitions and conclusions are derived from Chapter 7 of \cite{Villani_book2009}, which is a good reference for more information on dynamical cost functions and functionals in general theory of optimal transport. Assume that $(M,g)$ is a Riemannian manifold, $\mu,\nu\in\mathscr P(M)$ and $t>0$.

\begin{defn}[Dynamical cost functional]The \textit{dynamical cost functional} of $c^{t}(\cdot,\cdot)$ is defined as
\begin{align*}
	C^{t}(\mu,\nu):=\inf_{\gamma\in\Gamma(\mu,\nu)}\int_{M\times M}c^{t}(x,y)\,d\gamma
	%=\inf_{\substack{\mathrm{law}(X)=\mu\\\mathrm{law}(Y)=\nu}}\mathbb E(c^{t}(X,Y)),
\end{align*}
We denote by $\Gamma_o^{t}(\mu,\nu)$ the set of minimizers of $C^{t}(\mu,\nu)$. 
%where $X(\omega)$ and $Y(\omega)$ in the rightmost item are $M$-valued random variables from a probability measure space $(\Omega,\mathscr F,\mathbb P)$. 
%We denote by $\Gamma_o^{t}(\mu,\nu)$ the set of minimizers of $C^{t}(\mu,\nu)$. 
\end{defn}

Now we consider the case that $L$ is a Tonelli Lagrangian and the time-dependent cost function $c^t(x,y)=A_t(x,y)$ is the fundamental solution with respect to $L$. In this context, a connection from optimal transportation to Mather theory and weak KAM theory is firstly studied by Bernard and Buffoni in \cite{Bernard_Buffoni2007a}. 
\begin{defn}[Random curves]
For a given probability space $(\Omega,\mathscr F,\mathbb P)$ and $T>0$, An absolutely continuous \textit{random curve} on $[0,T]$ refers to a measurable map $\xi: [0,T]\times\Omega\to M$, which satisfies for any fixed $\omega\in\Omega$, sample path $\xi(\cdot,\omega)$ is absolutely continuous curve on $M$. Moreover, $\xi=\xi(s)$ can also be treated as a random process for $s\in[0,T]$.
\end{defn}

For the setting on work space $\mathscr P_c(\R^m)$, We select the probability measure space $(\Omega,\mathscr F,\mathbb P)$ for the random curves according to Remark \ref{rem:find the law}. In the subsequent discussion, we usually omit the component $\omega\in\Omega$ of the random curves and denote the expectation of a random variable (or vector) on 
$(\Omega,\mathscr F,\mathbb P)$ by $\mathbb E$ for simplicity when no confusion arises.

\begin{defn}[Dynamical coupling]
We say an absolutely continuous random curve $\xi$ is a \textit{dynamical coupling} of $\mu$ and $\nu$ for $t>0$, if $\law(\xi(0))=\mu$ and $\law(\xi(t))=\nu$. The set of all dynamical couplings of $\mu$ and $\nu$  for $t>0$ is denoted by $L_{\mu,\nu}^t$. 
\end{defn}
%
%\begin{defn}[Dynamical coupling and transference plan]
%We say a random curve $\xi$ is a \textit{dynamical coupling} of $\mu$ and $\nu$ for $t>0$, if $\xi:\Omega\to AC([0,t];M)$ is measurable, $\law(\xi(0,\cdot))=\mu$ and $\law(\xi(t,\cdot))=\nu$. The set of all dynamical couplings of $\mu$ and $\nu$  for $t>0$ is denoted by $L_{\mu,\nu}^t$. 
%\end{defn}

%Combining Theorem 7.21 and Remark 7.25 in \cite{Villani_book2009}, we have the following statements.

\begin{Pro}[\cite{Bernard_Buffoni2007a,Villani_book2009}]\label{pro:exist-displacement}
For any $\nu_1,\nu_2\in\mathscr P(\R^m)$ such that $C^{t}(\nu_1,\nu_2)$ is finite, 
\begin{align*}
	C^t(\nu_1,\nu_2)=\min_{\xi\in L_{\nu_1,\nu_2}^t}\E\left(\int_0^tL(\xi(s),\dot\xi(s))\,ds\right).
\end{align*}
%\begin{align*}
%	C^t(\mu,\nu)=\min_{\gamma\in\Gamma(\mu,\nu)}\int_{\R^{2m}}A_t(x,y)\,d\gamma=\min_{\xi\in L_{\mu,\nu}^t}\int_{\Omega}\int_0^tL(\xi(s,\omega),\dot\xi(s,\omega))\,dsd\mathbb P.
%\end{align*}
%where 
%\begin{align*}
%	L_{\mu,\nu}^t:=\{\xi\in L^0(\Omega;AC([0,t];\R^m)):\law(\xi(0,\cdot))=\mu,\law(\xi(t,\cdot))=\nu\}.
%\end{align*}
The minimizer of the last term is called dynamical optimal coupling of $\nu_1$ and $\nu_2$ for $t>0$. $\mu_s:=\law(\xi_*(s))$ $(s\in[0,t])$ as the law of a dynamical optimal coupling $\xi_*$ is called a displacement interpolation of $C^t(\nu_1,\nu_2)$. In this case, we have
\begin{itemize}
	\item $\xi_*$ is a random solution of the Euler-Lagrange equation;
	\item the path $\{\mu_s\}_{s\in[0,t]}$ is a minimizing curve for the action functional defined on $\mathscr P(\R^m)$ by
	\begin{align*}
		\mathbb A^t(\{\mu_s\}_{s\in[0,t]}):=\min_{\xi}\E\left(\int_0^tL(\xi(s),\dot\xi(s))\,ds\right),
		%=\sup_{N\in\N}\sup_{s=t_0<t_1<\cdots<t_N=t}\sum_{i=0}^{N-1}C^{t_{i+1}-t_i}(\mu_{t_i},\mu_{t_{i+1}})=\inf
	\end{align*}where the last minimum is over all random curves $\xi$ such that $\law(\xi(s))=\mu_s$ for $s\in[0,t]$.
\end{itemize}
\end{Pro}
For more discussions on dynamical optimal coupling and displacement interpolation, readers can refer to Chapter 7 in \cite{Villani_book2009}.

%\textcolor{blue}{
%In what follows, we let $(\Omega,\mathscr F,\mathbb P)$ be the probability space mentioned in 
%Remark \ref{rem:find the law} for our settings on $\mathscr P_c(\R^m)$. 
%}

\begin{Cor}\label{cor:dis-inter}
Suppose $\mu,\nu\in\mathscr P_c(\R^m)$ and $t>0$. 
	\begin{enumerate}[\rm (1)]
		\item There exists a compact set $K_{\mu,\nu,t}\subset\R^m$ depending on $\mu$, $\nu$ and $t$, such that for any displacement interpolation $\{\mu_s\}_{s\in[0,t]}$ of $C^t(\mu,\nu)$, $\supp(\mu_s)\subset K_{\mu,\nu,t}$ for all $s\in[0,t]$. Thus, $\{\mu_s\}_{s\in[0,t]}\subset\mathscr P_c(\R^m)$;
		\item $\{\mu_s\}_{s\in[0,t]}\in\Lip([0,t];\mathscr P_c(\R^m))$.
	\end{enumerate}	
\end{Cor}

\begin{proof}
Every displacement interpolation $\{\mu_s\}_{s\in[0,t]}$ of $C^t(\mu,\nu)$ is the law of some random action-minimizing curve $\xi$ at times $s$. Invoking Proposition \ref{pro:lip-esti}, we have that $\{\xi(\cdot,\omega)\}_{\omega\in\Omega}$ is equi-Lipschitz, and for each $\omega\in\Omega$, $\Lip(\xi(\cdot,\omega))\leqslant F(|\xi(0,\omega)-\xi(t,\omega)|/t)$, where $F$ is a non-decreasing function given in Proposition \ref{pro:lip-esti}. Since $\mu=\law(\xi(0))$ and $\nu=\law(\xi(t))$, the equi-Lipschitz constant $\Lip(\xi)$ of $\{\xi(\cdot,\omega)\}_{\omega\in\Omega}$ is determined by $\mu,\nu$ and $t>0$. Due to the compactness of $\supp(\mu)$ and $\supp(\nu)$, $\supp(\mu_s)\subset B_{\Lip(\xi)t}(\mu)=:K_{\mu,\nu,t}$ for all $s\in[0,t]$.
	
Now we turn to show the Lipschitz property of $\{\mu_s\}_{s\in[0,t]}$. Since $\mu_s=\law(\xi(s))$ for $s\in[0,t]$, we have
\begin{align*}
	W_2(\mu_s,\mu_{s'})\leqslant \E(d^2(\xi(s),\xi(s')))^{1/2}\leqslant\Lip(\xi)|s-s'|,\qquad \forall s,s'\in[0,t].
\end{align*}
%\begin{align*}
%	W_2(\mu_s,\mu_{s'})\leqslant \left\{\int_{\Omega}d^2(\xi(s,\omega),\xi(s',\omega))\,d\mathbb P\right\}^{1/2}\leqslant\Lip(\xi)|s-s'|,\qquad \forall s,s'\in[0,t].
%\end{align*}
It means $\{\mu_s\}_{s\in[0,t]}$ is a Lipschitz curve in $\mathscr P_c(\R^m)\subset\mathscr P_2(\R^m)$.
\end{proof}

The following essentially known results are the technical support of the study of dynamical cost functionals.

\begin{Pro}[\protect{\cite[Theorem 2.4]{Ambrosio_Brue_Semola_book2021}}, disintegration of measures] Let $X,Y$ be polish spaces, $\mu\in\mathscr P(X)$ and $f:X\to Y$ be a Borel function. Then there exists a family $\{\mu_y\}_{y\in Y}\subset\mathscr P(X)$ such that
\begin{enumerate}[\rm (1)]
	\item $y\mapsto \mu_y$ is Borel, i.e. $y\mapsto \mu_y(A)$ is Borel for any Borel set $A\subset X$;
	\item $\mu=\int_Y\mu_yd(f_\#\mu)$, i.e. for any Borel subset $A\subset X$, $\mu(A)=\int_Y\mu_y(A)d(f_\#\mu)$;
	\item $\mu_y$ is concentrated on $f^{-1}(y)$ for $(f_\#\mu)$-a.e. $y\in Y$.
\end{enumerate}
	Any other family $\{\mu_y'\}_{y\in Y}\subset\mathscr P(X)$ with these properties satisfies $\mu_y'=\mu_y$ for $(f_\#\mu)$-a.e. $y\in Y$. 
	
	Particularly, let $X=X_1\times X_2$, $Y=X_i$ 
 and $f=\pi_i$, $i=1,2$. For $\mu\in\mathscr P(X_1\times X_2)$, we say $\mu=\int_{X_1}\mu_{x_1}d(\pi_1)_\#\mu$ and $\mu=\int_{X_2}\mu_{x_2}d(\pi_2)_\#\mu$ are the disintegrations of $\mu$ with respect to $(\pi_1)_\#\mu$ and $(\pi_2)_\#\mu$ (or $X_1$ and $X_2$), respectively.
 \end{Pro}

%\begin{Lem}[\protect{\cite[Lemma 5.3.2]{Ambrosio_GigliNicola_Savare_book2008}}]\label{lem:gluing lemma}
%Let $X_i$ $(i=1,2,3)$ be Polish spaces, $\gamma^{1,2}\in\mathscr P(X_1\times X_2)$ and $\gamma^{1,3}\in\mathscr P(X_1\times X_3)$ with  $(\pi_1)_\#\gamma^{1,2}=(\pi_1)_\#\gamma^{1,3}=\mu^{1}$. Then there exists $\boldsymbol{\mu}\in\mathscr P(X_1\times X_2\times X_3)$, which satisfies
%\begin{equation}\label{eq:disintegration}
%	(\pi_{1,2})_\#\boldsymbol{\mu}=\gamma^{1,2},\qquad(\pi_{1,3})_\#\boldsymbol{\mu}=\gamma^{1,3}.
%\end{equation}
%		
%Moreover, if $\gamma^{1,2}=\int\gamma_{x_1}^{1,2}d\mu^1$, $\gamma^{1,3}=\int\gamma_{x_1}^{1,3}d\mu^1$ and $\boldsymbol{\mu}=\int\boldsymbol{\mu}_{x_1}d\mu^1$ are the disintegrations of $\gamma^{1,2}$, $\gamma^{1,3}$ and $\boldsymbol{\mu}$ with respect to $\mu^1$ respectively, then \eqref{eq:disintegration} is equivalent to 
%\begin{align*}
%	\boldsymbol{\mu}_{x_1}\in\Gamma(\gamma_{x_1}^{1,2},\gamma_{x_1}^{1,3})\subset\mathscr P(X_2\times X_3)\qquad \text{for}\ \mu^1-a.e.\ x_1\in X_1.
%\end{align*}
%\end{Lem}

\begin{defn}[Gluing measure]
Let $\mu_i\in\mathscr P(X_i)$ for $i=1,2,3$. Set $\gamma\in\Gamma(\mu_1,\mu_2)$ and $\gamma'\in\Gamma(\mu_1,\mu_3)$. Then we denote
\begin{align*}
	\gamma*\gamma':=\int_{X_1}\gamma_{x_1}\times\gamma'_{x_1}d\mu_1\in\Gamma(\mu_1,\mu_2,\mu_3)\subset\mathscr P(X_1\times X_2\times X_3)
\end{align*}
as \textit{the gluing measure of $\gamma$ and $\gamma'$ (with respect to $\mu_1$)}, where $\gamma=\int_{X_1}\gamma_{x_1}d\mu_1$ and $\gamma'=\int_{X_1}\gamma'_{x_1}d\mu_1$ are the disintegration of $\gamma$ and $\gamma'$ with respect to $\mu_1$ respectively.
\end{defn}

%Hereafter, we denote $\boldsymbol{\mu}$ mentioned above by $\boldsymbol{\mu}=\gamma^{1,2}*\gamma^{1,3}$, the gluing measure of $\gamma^{1,2}$ and $\gamma^{1,3}$, which have the same marginal $\mu^1$ on $X_1$.

%\begin{Lem}\label{lem:gluing_lambda}
%Assume that $X_i=\R^m$, $\mu_i\in\mathscr P_c(X_i)$ with $i=1,2,3$ and $\gamma\in\Gamma(\mu_2,\mu_3)$. For any fixed $\lambda\in[0,1]$ and $t>0$, let $\mu_{\lambda}\in\Gamma_c^{t}(\mu_1,\gamma_{2\to 3}^\lambda)$. Then there exists $\gamma_\lambda\in\Gamma(\mu_1,\gamma)$ such that $(\gamma_\lambda)_{1,2\to 3}^\lambda=\mu_\lambda$. 
%\end{Lem}
%
%\begin{proof}
%The idea of the proof can be derived from that of Proposition 7.3.1 in \cite{Ambrosio_GigliNicola_Savare_book2008} and the technique required is given by Lemma \ref{lem:gluing lemma}.	
%\end{proof}

\begin{The}\label{thm:dyn-cost}
Suppose $t>0$, $K\subset\R^m$ is a convex and compact subset.
\begin{enumerate}[\rm (1)]
	\item  $C^{t}(\mu,\cdot)$ and $C^{t}(\cdot,\mu)$ are superlinear on $\mathscr P_c(\R^m)$ endowed with $W_1$ metric, where $\mu\in\mathscr P_c(\R^m)$;
	\item  \( (t, \nu) \mapsto C^t(\mu, \nu) \in \mbox{\rm Lip}([a,b] \times \mathscr{P}(K)) \), where $\mu\in\mathscr P(K)$;
	\item  for $t\in(0,1)$, there exists $C_K>0$ depends on $K$, such that for any $\mu\in\mathscr P(K)$, $C^t(\mu,\cdot)$ is locally strongly semiconcave on $\mathscr P(K)$ with constant $\frac{C_K}{t}$;
\end{enumerate}	
\end{The}

See the appendix for the proof of Theorem \ref{thm:dyn-cost}.

\begin{The}[Generalized superdifferential of $C^{t}(\mu,\nu)$]\label{thm:cost-supdiff}
Define $\mu,\nu\in\mathscr P_c(\R^m)$, $t>0$ and $\gamma\in\Gamma_o^{t}(\mu,\nu)$, as well as $\gamma=\int_{\R^m}\gamma_x\,d\mu(x)=\int_{\R^m }\gamma_y\,d\nu(y)$ are the disintegrations of $\gamma$ with respect to $\mu$ and $\nu$ respectively. Moreover, we choose $\eta_{x,y}\in\Gamma_{x,y}^t$, the minimizer of $A_t(x,y)$, which is measurable with respect to $(x,y)\in\supp(\mu)\times\supp(\nu)$\footnote{This can be guaranteed by standard measurable selection theorem (see \cite{Clarke_book2013}).}. Then	\begin{align*}
		\mathbf p_\nu(y):=\int_{\R^m}p_{x,y}(t)\,d\gamma_{y}&\in\partial_\nu^+C^{t}(\mu,\nu),\\
		-\mathbf p_\mu(x):=\int_{\R^m}-p_{x,y}(0)\,d\gamma_{x}&\in\partial_\mu^+C^{t}(\mu,\nu),\\
		-\int_{\R^{2m}}H(y,p_{x,y}(t))\,d\gamma &\in D_t^+C^{t}(\mu,\nu),
	\end{align*}where $p_{x,y}(s):=L_v(\eta_{x,y}(s),\dot\eta_{x,y}(s))$, $s\in[0,t]$.
\end{The}

For the proof of Theorem \ref{thm:cost-supdiff}, see the appendix.

\begin{Cor}\label{cor:HJ_ineq}
If $(t,\nu)\mapsto C^t(\mu,\nu)$ is locally differentiable, we have
\begin{align*}
	D_tC^t(\mu,\nu)+\int_{\R^m}H(y,\partial_\nu C^t(\mu,\nu))\,d\nu\leqslant 0.
\end{align*}
The strict inequality holds in general.
\end{Cor}

\begin{proof}
If $(t,\nu)\mapsto C^t(\mu,\nu)$ is localized differentiable, according to Theorem \ref{thm:cost-supdiff},
\begin{align*}
	D_tC^{t}(\mu,\nu)=-\int_{\R^{2m}}H(y,p_{x,y}(t))\,d\gamma,\qquad
	\partial_\nu C^{t}(\mu,\nu)=\int_{\R^m}p_{x,y}(t)\,d\gamma_{y}.
\end{align*} 
Due to the fact that $H(x,\cdot)$ is strictly convex, Jensen's inequality implies
\begin{align*}
	&\,D_tC^t(\mu,\nu)+\int_{\R^m}H(y,\partial_\nu C^t(\mu,\nu))\,d\nu\\
		=&\,-\int_{\R^{2m}}H(y,p_{x,y}(t))\,d\gamma+\int_{\R^m}H\left(y,\int_{\R^m}p_{x,y}(t)\,d\gamma_{y}\right)\,d\nu\\
		\leqslant&\,-\int_{\R^{2m}}H(y,p_{x,y}(t))\,d\gamma+\int_{\R^m}\int_{\R^m}H\left(y,p_{x,y}(t)\right)\,d\gamma_{y}d\nu=0.
\end{align*}
In general, the inequality above is strict if $\mu$ and $\nu$ are not Dirac measures.
\end{proof}

\begin{The}\label{thm:semi-convex}
Given $\mu_1\in\mathscr P_c(\R^m)$. Suppose $\lambda,R>0$ with $\supp(\mu_1)\subset B(0,R)$ and let $t\in(0,t_{\lambda,R})$, where $t_{\lambda,R}$ is determined by Proposition \ref{pro:cvpde}. Then for all $\mu_2,\mu_3\in\mathscr P(B(0,R))$ those satisfy the following condition: one can find $\gamma_{1,i}\in \Gamma_o^t(\mu_1,\mu_i)$, $i=2,3$ such that $|x_1-x_i|\leqslant\lambda t$ holds true for $\gamma_{1,i}$-a.e. $(x_1,x_i)$, we have
\begin{align*}
	(1-s)C^t(\mu_1,\mu_2)+sC^t(\mu_1,\mu_3)-C^t(\mu_1,\gamma_{2\to 3}^s)\geqslant \frac{C_{\lambda,R}}{t}s(1-s)W^2_{2,\gamma}(\mu_2,\mu_3),
\end{align*}
where $\gamma:=\gamma_{1,2}*\gamma_{1,3}\in\Gamma(\mu_1,\mu_2,\mu_3)$ and $W^2_{2,\gamma}(\mu_2,\mu_3):=\int_{\R^{3m}}|x_2-x_3|^2\,d\gamma$. In this case,  we call $C^t(\mu_1,\cdot)$ locally generalized semiconvex with constant $\frac{C_{\lambda,R}}{t}$.
\end{The}

\begin{proof}
Taking $\mu_2,\mu_3\in\mathscr P(B(0,R))$ satisfying the condition in the theorem above, then for all $s\in[0,1]$, $\gamma_{1,2\to 3}^s\in\Gamma(\mu_1,\gamma_{2\to 3}^s)$ and
\begin{align*}
	&\,(1-s)C^t(\mu_1,\mu_2)+sC^t(\mu_1,\mu_3)-C^t(\mu_1,\gamma_{2\to 3}^s)\\
		\geqslant&\,\int_{\R^{2m}}(1-s)A_t(x,y)\,d\gamma_{1,2}+\int_{\R^{2m}}sA_t(x,y)\,d\gamma_{1,3}-\int_{\R^{2m}}A_t(x,y)\,d\gamma_{1,2\to 3}^s\\
		\geqslant&\,\int_{\R^{3m}}(1-s)A_t(x_1,x_2)+sA_t(x_1,x_3)-A_t(x_1,(1-s)x_2+sx_3)\,d\gamma\\
		\geqslant&\,\int_{\R^{3m}}\frac{C_{\lambda,R}}{t}s(1-s)|x_2-x_3|^2\,d\gamma\\
		=&\,\frac{C_{\lambda,R}}{t}s(1-s)W_{2,\gamma}^2(\mu_2,\mu_3).
\end{align*}
Here we used the semiconvexity result of fundamental solutions in Proposition \ref{pro:cvpde}.
\end{proof}

\section{Random Lax-Oleinik operators and viscosity solutions}

In this section, we define Lax-Oleinik operators in $\mathscr P(\R^m)$ and investigate its properties in $\mathscr P_c(\R^m)$. It can be seen that they inherit many characteristics of the classical Lax-Oleinik operators, so it can be proved that they are viscosity solutions of the corresponding Hamilton-Jacobi equation in $\mathscr P_c(\R^m)$. We highlight the approach of using the Hopf-Lax semigroup, as developed in previous works \cite{Kuwada2010, Gigli_Kuwada_Ohta2013}, within the context of optimal transport.

\subsection{Definitions and properties on Random Lax-Oleinik operators}

\begin{defn}[Random Lax-Oleinik operators]\label{def:Pt-,Pt+}
For any Borel measurable function $\phi:\R^m\to [-\infty,+\infty]$, $t>0$ and $\mu\in\mathscr P(\R^m)$, we define
\begin{align}
	P_t^+\phi(\mu)&:=\sup_{\nu\in\mathscr P(\R^m)}\{\phi(\nu)-C^t(\mu,\nu)\},\label{eq:def-pt+}\\
	P_t^-\phi(\mu)&:=\inf_{\nu\in\mathscr P(\R^m)}\{\phi(\nu)+C^t(\nu,\mu)\}.\label{eq:def-pt-}
\end{align}
$P_t^+$ and $P_t^-$ are called \emph{positive} and \emph{negative type Lax-Oleinik operators} respectively.
\end{defn}

The following theorem and corollary show that $\mathscr P_c(\R^m)$ is a suitable working space for the operator $P_t^{\pm}$. Before proceeding, we introduce a lemma necessary for later discussion.

\begin{Lem}[\protect{\cite[Lemma 4.3]{Cheng_Hong_Shi2024}}]\label{lem:proj}
Assume that $X$ and $Y$ are two Polish spaces and $\gamma\in\mathscr{P}(X\times Y)$.
\begin{enumerate}[\rm (1)]
	\item $\pi_X(\supp(\gamma))\subset\supp\,((\pi_X)_\#\gamma)$;
	\item $\pi_X(\supp(\gamma))=\supp((\pi_X)_\#\gamma)$ if $Y$ is compact.
\end{enumerate}
\end{Lem}

\begin{The}\label{thm:pt-}
Suppose that $t>0$, $\phi$ is lower semicontinuous on $\R^m$ and $(\kappa_1,\kappa_2)$-Lipschitz in the large. Then the followings hold.
\begin{enumerate}[\rm (1)]
	\item For any $\nu\in\mathscr P_c(\R^m)$, there exists a minimizer $\mu=\mu(t,\nu)\in\mathscr P_c(\R^m)$ such that
	\begin{align*}
		P_t^-\phi(\nu)=\phi(\mu)+C^t(\mu,\nu).
	\end{align*}
	\item $P_t^-\phi(\nu)=T_t^-\phi(\nu)$ is finite for all $\nu\in\mathscr P_c(\R^m)$. That is $P_t^-\phi$ is the potential energy functional induced by $T_t^-\phi$ on $\mathscr P_c(\R^m)$.
	\item For $\nu\in\mathscr P_c(\R^m)$, each minimizer $\mu(t,\nu)$ is contained in $\mathscr P_c(\R^m)$ and makes \eqref{eq:K-} hold true.
	\item For any minimizer $\mu(t,\nu)$, there exist constants $C_1,C_2>0$ such that
	\begin{align*}
		W_1(\mu(t,\nu),\nu)\leqslant C_1t+C_2,
	\end{align*}
	where $C_1=c_0+\theta_0^*(\kappa_1+1)+\theta_1(0)+c_1$ and $C_2=\kappa_2$, with $\theta_0^*$ the convex conjugate function of $\theta_0$.
\end{enumerate}
\end{The}

\begin{proof}
Firstly, to prove part (1), we begin by noting that there exists a probability space \((\Omega, \mathscr{F}, \mathbb{P})\) and a random variable \(y_\nu: \Omega \to \mathbb{R}^m\) such that the law of \(y_\nu\) is \(\nu\), by Remark \ref{rem:find the law}. Leveraging the a priori estimates on the minimizer of \(T_t^-\phi\) (refer to Lemma \ref{pro:funda-minimizer}), we can express the marginal function as:
\begin{align*}
	T_t^-\phi(y_\nu(\omega))=\min_{x\in B(y_\nu(\omega),\lambda_\phi t+\kappa_2)}\{\phi(x)+A_t(x,y_\nu(\omega))\},\,\,\,\forall\omega\in\Omega,
\end{align*}
where \(\lambda_\phi > 0\) depends only on \(L\) and \(\kappa_1\). It's straightforward to verify that the map \(\omega \mapsto B(y_\nu(\omega), \lambda_\phi t + \kappa_2)\) is a measurable and compact set-valued map. By invoking Theorem 18.19 from \cite{Aliprantis_Border_book2006}, we can admit a measurable selection \(x_{t,\nu}\)
\begin{align*}
	\omega\mapsto x_{t,\nu}(\omega)\in \mathop{\arg\min}_{x\in B(y_\nu(\omega),\lambda_\phi t+\kappa_2)}\{\phi(x)+A_t(x,y_{\nu}(\omega))\}\subset B(y_\nu(\omega),\lambda_\phi t+\kappa_2).
\end{align*}
Consequently, defining \(\mu(t,\nu) :=\law(x_{t,\nu})\), we derive the inequality:
%
%Consequently, defining \(\mu(t, \nu) := (x_{t,\nu})_\#\mathbb{P}\), we derive the inequality:
\begin{equation}\label{eq:geq_p}
	\begin{split}
		\int_{\R^m}T_t^-\phi(y)\,d\nu &=\E(T_t^-\phi(y_{\nu}))\\
		&=\E(\phi(x_{t,\nu})+A_t(x_{t,\nu},y_\nu))\\
		&\geqslant\phi(\mu(t,\nu))+C^{t}(\mu(t,\nu),\nu)\geqslant P_t^-\phi(\nu).
	\end{split}
\end{equation}
%\begin{equation}\label{eq:geq_p}
%	\begin{split}
%		\int_{\R^m}T_t^-\phi(y)\,d\nu &=\int_\Omega T_t^-\phi(y_{\nu}(\omega))\,d\mathbb P\\
%		&=\int_\Omega\phi(x(\omega))+A_t(x_{t,\nu}(\omega),y_\nu(\omega))\,d\mathbb P\\
%		&=\int_{\R^{2m}}\phi(x)+A_t(x,y)\,d(x_{t,\nu},y_\nu)_\#\mathbb P\\
%		&\geqslant\phi(\mu(t,\nu))+C^{t}(\mu(t,\nu),\nu)\geqslant P_t^-\phi(\nu).
%	\end{split}
%\end{equation}
On the other hand, since $T_t^-\phi(y)\leqslant \phi(x)+A_t(x,y)$ for all $x,y\in\R^m$, then given $\nu\in\mathscr P_c(\R^m)$, for arbitrary $\mu\in\mathscr P_c(\R^m)$ and $\gamma\in\Gamma(\mu,\nu)$, we have
\begin{align*}
	\int_{\R^m}T_t^-\phi(y)\,d\nu=\int_{\R^{2m}}T_t^-\phi(y)\,d\gamma\leqslant\int_{\R^{2m}}\phi(x)+A_t(x,y)\,d\gamma.
\end{align*} 
This implies that
\begin{equation}\label{eq:leq_p}
	T_t^-\phi(\nu)=\int_{\R^m}T_t^-\phi(y)\,d\nu\leqslant\inf_{\mu\in\mathscr P(\R^m)}\inf_{\gamma\in\Gamma(\mu,\nu)}\int_{\R^{2m}}\phi(x)+A_t(x,y)\,d\gamma=P_t^-\phi(\nu).
\end{equation}
Combining \eqref{eq:geq_p} and \eqref{eq:leq_p} we find that
\begin{align*}
	P_t^-\phi(\nu)=\phi(\mu(t,\nu))+C^{t}(\mu(t,\nu),\nu).
\end{align*}
Additionally, we claim that $\mu(t,\nu)\in\mathscr P_c(\R^m)$. Indeed, by Corollary 5.22 from \cite{Villani_book2009} and the existence of random curves, we can find a dynamical optimal coupling \(\xi \in L_{\mu, \nu}^t\) such that

%Additionally, we claim that \(\mu \in \mathscr{P}_c(\mathbb{R}^m)\). Indeed, the map \(\omega \mapsto (x_{t,\nu}(\omega), y_\nu(\omega))\) is measurable, and by Corollary 5.22 from \cite{Villani_book2009} and the existence of random curves, we can find \(\xi(\omega) \in \Gamma_{x_{t,\nu}(\omega), y_\nu(\omega)}^t\) such that
\begin{align*}
	T_t^-\phi(y_\nu)=\phi(x_{t,\nu})+\int_0^tL(\xi(s),\dot\xi(s))\,ds
\end{align*}

%\begin{align*}
%	T_t^-\phi(y_\nu(\omega))=\phi(x_{t,\nu}(\omega))+\int_0^tL(\xi(s,\omega),\dot\xi(s,\omega))\,ds
%\end{align*}
Furthermore, by Proposition \ref{pro:lip-esti}, for any \(\omega \in \Omega\), \(\xi(\cdot,\omega) \in \Lip([0, t]; \mathbb{R}^m)\) with a Lipschitz constant dependent only on \(t\), \(L\), and \((\kappa_1, \kappa_2)\). Therefore, setting the uniform Lipschitz constant as \(\Lip(\xi)\), we have that
\begin{align*}
	\supp(\mu(t,\nu))&=\supp(\law(\xi(0)))\\
		&\subset \bigcup_{\omega\in\Omega}B(\xi(t,\omega),\Lip(\xi)t)=\bigcup_{y\in\supp(\nu)}B(y,\Lip(\xi)t),
\end{align*} 
%
%\begin{align*}
%	\supp(\mu(t,\nu))&=\supp(\law(\xi(0)))\\
%		&\subset \bigcup_{\omega\in\Omega}B(\xi(t),\Lip(\xi)t)=\bigcup_{y\in\supp(\nu)}B(y,\Lip(\xi)t),
%\end{align*} 

Conclusion (2) follows directly from \eqref{eq:geq_p} and \eqref{eq:leq_p}.

Now, we turn to prove part (3). Let \(\mu\in\arg\min_{\mathscr{P}(\mathbb{R}^m)} \{\phi(\cdot) + C^t(\cdot, \nu)\}\) and let \(\gamma\in\Gamma_o^t(\mu, \nu)\). From part (2), it follows that
\begin{align*}
	\int_{\R^m}T_t^-\phi(y)\,d\nu=P_t^-\phi(\nu)=\phi(\mu)+C^t(\mu,\nu)=\int_{\R^{2m}}\phi(x)+A_t(x,y)\,d\gamma
\end{align*}
We also have that
\begin{equation}\label{eq:gamma-a.e.}
	T_t^-\phi(y)=\phi(x)+A_t(x,y),\,\,\,\,\gamma-\mathrm{a.e.}\,(x,y)\in\R^{2m}.
\end{equation}
This implies that the minimizer \(\mu\) satisfies equation \eqref{eq:K-}. Furthermore, we can deduce that
\begin{equation}\label{eq:gamma-supp}
	T_t^-\phi(y)=\phi(x)+A_t(x,y),\,\,\,\,\,\forall(x,y)\in\supp(\gamma).
\end{equation}
In fact, by combining the lower semicontinuity of \(\phi\) with the properties of \(T_t^-\phi\), we can define the function
\begin{align*}
	F_\phi(x,y):=\phi(x)+A_t(x,y)-T_t^-\phi(y)
\end{align*}
which is a non-negative and lower semicontinuous function on \(\mathbb{R}^{2m}\). If there were some \((x_0, y_0) \in \supp(\gamma)\) such that \(F_\phi(x_0, y_0) > 0\), then we could find some \(R_0 > 0\) such that \(F_\phi(x, y)\) is positive on the ball \(B((x_0, y_0), R_0)\) with \(\gamma(B((x_0, y_0), R_0)) > 0\), which would contradict equation \eqref{eq:gamma-a.e.}. According to part (1) of Lemma \ref{lem:proj}, we know that
\begin{equation}\label{eq:proj-compact}
	\pi_2(\supp(\gamma))\subset\supp((\pi_2)_\#\gamma)=\supp(\nu),
\end{equation}
For any fixed \(y \in \mathbb{R}^m\) satisfying equation \eqref{eq:gamma-a.e.}, by Lemma \ref{pro:funda-minimizer}, when \(x\) satisfies the same equation, \(x\) must be in \(B(y, \lambda_\phi t + \kappa_2)\). Due to the compact projection, we have:
\begin{align*}
	\supp(\gamma)\subset \bigcup_{y\in\supp(\nu)}\{(z,y):z\in B(y,\lambda_\phi t+\kappa_2)\}.
\end{align*}
Additionally, from the compact projection, there exists a compact set \(K \subset \mathbb{R}^m\) such that \(\supp(\gamma) \subset \mathbb{R}^m \times K\). Following from part (2) of Lemma \ref{lem:proj}, we get:
\begin{equation}\label{eq:mu-compact support}
	\begin{split}
		\supp(\mu)&=\supp((\pi_1)_{\#}\gamma)\\
		&=\pi_1(\supp(\gamma))\\
		&\subset\pi_1\left(\bigcup_{y\in\supp(\nu)}\{(z,y):z\in B(y,\lambda_\phi t+\kappa_2)\}\right)=\bigcup_{y\in\supp(\nu)}B(y,\lambda_\phi t+\kappa_2),
	\end{split}
\end{equation}
which means that \(\supp(\mu)\) is compact and every minimizer of \(P_t^-\phi(\nu)\) is contained in \(\mathscr{P}_c(\mathbb{R}^m)\).

Finally, we aim to prove part (4). Since \(\mu\) is the minimizer, we have \(\phi(\mu) + C^t(\mu, \nu) \leqslant \phi(\nu) + C^t(\nu, \nu)\). Moreover, by Theorem \ref{thm:dyn-cost} and the superlinearity of \(C^t\), for any \(k \in \mathbb{R}\):
\begin{align*}
	C^t(\mu,\nu)\geqslant kW_1(\mu,\nu)-(c_0+\theta_0^*(k))t.
\end{align*}
We also have:
\begin{equation}\label{eq:w1-lip in the large}
	kW_1(\mu,\nu)-(c_0+\theta_0^*(k))t-C^t(\nu,\nu)\leqslant C^t(\mu,\nu)-C^t(\nu,\nu)\leqslant \phi(\nu)-\phi(\mu).
\end{equation}
Suppose \(\gamma \in \Gamma_1(\mu, \nu)\). Because \(\phi\) is \((\kappa_1, \kappa_2)\)-Lipschitz in the large, \eqref{eq:w1-lip in the large} yields:
\begin{align*}
	kW_1(\mu,\nu)-(c_0+\theta_0^*(k))t-C^t(\nu,\nu)&\leqslant\int_{\R^{2m}}\phi(y)-\phi(x)\,d\gamma\\
		&\leqslant\int_{\R^{2m}}\kappa_1|x-y|+\kappa_2\,d\gamma\\
		&=\kappa_1W_1(\mu,\nu)+\kappa_2.
\end{align*}
Taking \(k = \kappa_1 + 1\), we obtain from property (L2) of \(L\):
\begin{align*}
	W_1(\mu,\nu)&\leqslant (c_0+\theta_0^*(\kappa_1+1))t+C^t(\nu,\nu)+\kappa_2\\
		&\leqslant(c_0+\theta_0^*(\kappa_1+1))t+\int_{\R^{m}}A_t(y,y)\,d\nu+\kappa_2\\
		&\leqslant(c_0+\theta_0^*(\kappa_1+1))t+\int_{\R^{m}}\int_0^tL(y,0)\,ds\,d\nu+\kappa_2\\
		&\leqslant (c_0+\theta_0^*(\kappa_1+1)+\theta_1(0)+c_1)t+\kappa_2.
\end{align*}
This completes the proof.
\end{proof}

By similar reasoning, we can also derive the corresponding conclusions for the positive random Lax-Oleinik operator \(P_t^+\).

\begin{Cor}\label{cor:pt+}
Suppose $t>0$, $\phi$ is upper semicontinuous on $\R^m$ and $(\kappa_1,\kappa_2)$-Lipschitz in the large. Then the followings hold.
\begin{enumerate}[\rm (1)]
	\item For any $\mu\in\mathscr P_c(\R^m)$, there exists a maximizer $\nu=\nu(t,\mu)\in\mathscr P_c(\R^m)$ such that
	\begin{align*}
		P_t^+\phi(\mu)=\phi(\nu)-C^t(\mu,\nu).
	\end{align*}
	\item $P_t^+\phi(\mu)=T_t^+\phi(\mu)$ is finite for all $\mu\in\mathscr P_c(\R^m)$. That is $P_t^+\phi$ is the potential energy functional induced by $T_t^+\phi$ on $\mathscr P_c(\R^m)$.
	\item For $\mu\in\mathscr P_c(\R^m)$, each maximizer $\nu(t,\mu)$ is contained in $\mathscr P_c(\R^m)$ and makes \eqref{eq:K+} hold true.
	\item For any maximizer $\nu(t,\mu)$, there exist constants $C_1,C_2>0$ such that
	\begin{align}\label{eq:argmax-wasser}
		W_1(\nu(t,\mu),\mu)\leqslant C_1t+C_2,
	\end{align}where $C_1=c_0+\theta_0^*(\kappa_1+1)+\theta_1(0)+c_1$ and $C_2=\kappa_2$.
\end{enumerate}
\end{Cor}

\begin{Pro}\label{pro:probability lax-oleinik}
Let $\phi\in \UC(\R^m)$, $t_0>0$ and $c[0]=0$. Random Lax-Oleinik operators have the following properties:
\begin{enumerate}[\rm (1)]
		\item $\{P_t^-\}_{t>0}$ and $\{P_t^+\}_{t>0}$ are semigroups with respect to $t\in(0,+\infty)$;
%		\item $P_t^-,P_t^+:\UC(\mathscr P_c(\R^m))\to \UC(\mathscr P_c(\R^m))$;
		\item $\lim_{t\to 0^+}P_t^-\phi=\lim_{t\to 0^+}P_t^+\phi=\phi$, thus we can define $P_0^-\phi=P_0^+\phi=\phi$;
		\item for arbitrary $t,s\geqslant 0$, $P_t^-\circ P_s^+\phi=T_t^-\circ T_s^+\phi(\cdot)$; $P_t^+\circ P_s^-\phi=T_t^-\circ T_s^-\phi(\cdot)$;
		\item for any $t>0$, $P_t^-\phi$ (resp. $-P_t^+\phi$) is locally (resp. strongly) semiconcave and $\{P_t^-\phi\}_{t\geqslant t_0}$ (resp. $\{-P_t^+\phi\}_{t\geqslant t_0}$) uniformly localized (resp. strongly) semiconcave;
		\item $t\mapsto P_t^-\phi$ ($t\mapsto P_t^+\phi$) uniformly continuous on $[0,+\infty)$; 
		\item suppose $K\subset\R^m$ is compact set. $(t,\mu)\mapsto P_t^-\phi(\mu)$ ($(t,x)\mapsto P_t^+\phi(x)$) is continuous on $[0,+\infty)\times \mathscr P(K)$, locally Lipschitz on $(0,+\infty)\times \mathscr P(K)$ and equi-Lipschitz on $[t_0,+\infty)\times \mathscr P(K)$ with respect to $\phi$.
	\end{enumerate}
\end{Pro}

\begin{proof}
It should be noted that the operators \( P_t^\pm \) represent the potential energies of the functionals induced by \( T_t^\pm \) on the space of compactly supported probability measures \( \mathscr{P}_c(\mathbb{R}^m) \), as established in Theorem \ref{thm:pt-} and Corollary \ref{cor:pt+}. Consequently, by drawing an analogy with the results presented in Proposition \ref{pro:lax-oleinik}, we can readily deduce the aforementioned conclusions.
\end{proof}

\begin{Pro}\label{pro:sub-hjs}
The following statements are equivalent.
\begin{enumerate}[\rm (1)]
	\item $\phi$ is the subsolution to \eqref{hjs}; 
	\item $\phi(\nu)-\phi(\mu)\leqslant C^t(\mu,\nu)+c[0]t$ for all $t>0$ and $\mu,\nu\in\mathscr P_c(\R^m)$;
	\item For any $\mu\in\mathscr P_c(\R^m)$, $[0,+\infty)\ni t\mapsto P_t^-\phi(\mu)+c[0]t$ is non-decreasing;
	\item For any $\mu\in\mathscr P_c(\R^m)$, $[0,+\infty)\ni t\mapsto P_t^+\phi(\mu)-c[0]t$ is non-increasing.
\end{enumerate}
\end{Pro}

\begin{proof}
According to \cite{Fathi_book}, $\phi$ is a subsolution to equation \eqref{hjs} if and only if for any $x, y \in \mathbb{R}^m$ and $t>0$, the inequality $\phi(y) - \phi(x) \leqslant A_t(x, y)+c[0]t$ holds. Then, for any fixed $\mu, \nu \in \mathscr{P}_c(\mathbb{R}^m)$ and $\gamma \in \Gamma_o(\mu, \nu)$, integrating both sides of this inequality over $(x, y)$ with respect to $\gamma$ yields equation (2). Conversely, if equation (2) holds, setting $\mu = \delta_x$ and $\nu = \delta_y$ implies that inequality (1) is satisfied. The equivalence among equations (2), (3), and (4) can be derived from Definition \ref{def:Pt-,Pt+}.
\end{proof}

Building on Theorem~\ref{thm:pt-} and Corollary~\ref{cor:pt+}, the following statement is a direct corollary of \cite[Lemma 3.1]{Cannarsa_Cheng_Hong2025}.

\begin{Pro}\label{pro:probability commutator}
Suppose $\phi\in \UC(\R^m)$ and $\mu\in\mathscr P_c(\R^m)$.
\begin{enumerate}[\rm (1)]
	\item For any $t\geqslant 0$, $P_t^-\circ P_t^+\phi\geqslant\phi(\cdot)$ and $P_t^-\circ P_t^+\phi(\mu)=\phi(\mu)$ if and only if there exists $\nu\in\mathscr P_c(\R^m)$ such that $P_t^+\phi(\nu)=\phi(\mu)-C^t(\nu,\mu)$.
	\item For any $t\geqslant 0$, $P_t^+\circ P_t^-\phi\leqslant\phi(\cdot)$ and $P_t^+\circ P_t^-\phi(\mu)=\phi(\mu)$ if and only if there exists $\nu\in\mathscr P_c(\R^m)$ such that $P_t^-\phi(\nu)=\phi(\mu)+C^t(\mu,\nu)$;
	\item For any $t\geqslant 0$, $P_t^+\circ P_t^-\circ P_t^+\phi=P_t^+\phi$ and $P_t^-\circ P_t^+\circ P_t^-\phi=P_t^-\phi$.
\end{enumerate} 
\end{Pro}

\subsection{Hamilton-Jacobi equations in $\mathscr P_c(\R^m)$}

Let \(\phi \in C(\mathbb{R}^m)\) be a function \((\kappa_1, \kappa_2)\)-Lipschitz in the large. According to Theorem~\ref{thm:hje}, as the value function of the Bolza problem, \(u(t, x) := T_t^- \phi(x)\) is the viscosity solution to equation \eqref{hje-} with the initial condition \(\phi\). A natural question arises: whether \(U(t, \mu) := P_t^- \phi(\mu)\) will be the viscosity solution to some corresponding form of the Hamilton-Jacobi equations.

Assume that $H\in C(\R^{m}\times\R^{m})$ is a Tonelli Hamiltonian. We consider the following Hamilton-Jacobi equation: 
\begin{align}
	\partial_tU(t,\mu)+\int_{\R^m}H(x,\partial_\mu U(t,\mu)(x))\,d\mu =0,\tag{PHJ$_{e}$}\label{phje}
\end{align}
where $(t,\mu)\in\R^+\times\mathscr P_c(\R^m)$.

\begin{defn}[Viscosity solution of \eqref{phje}]\label{defn:viscosity}
We call the functional $U:\R^+\times\mathscr P_c(\R^m)$ a \textit{viscosity subsolution} of \eqref{phje}, if for every $(t,\mu)\in\R^+\times\mathscr P_c(\R^m)$, 
\begin{align}
	q+\int_{\R^m}H(x,\alpha(x))\,d\mu\leqslant0,\,\,\,\,\,\,\,\,\forall (q,\alpha)\in\partial^+U(t,\mu);\label{eq:phje-sub}
\end{align}
Similarly, $U:\R^+\times\mathscr P_c(\R^m)$ is called a \textit{viscosity supersolution} of \eqref{phje}, if for every $(t,\mu)\in\R^+\times\mathscr P_c(\R^m)$,
\begin{align}
	p+\int_{\R^m}H(x,\beta(x))\,d\mu\geqslant0,\,\,\,\,\,\,\,\,\forall (p,\beta)\in\partial^-U(t,\mu).\label{eq:phje-sup}
\end{align}
If $U$ is both a subsolution and supersolution to \eqref{phje}, we say $U$ is a \emph{viscosity solution} of \eqref{phje}. 
\end{defn}

In a more general and abstract frame, Gangbo, Nguyen and Tudorascu introduced concept of Hamilton-Jacobi equations in Wasserstein spaces in \cite{Gangbo_Nguyen_Tudorascu2008}, which is further studied in \cite{Gangbo_Tudorascu2019} later by Gangbo and Tudorascu. Although, \eqref{phje} can be regarded as a special case of their work, it is closely related to main contents of this paper.

\begin{The}\label{thm:phje-cauchy}
Let $\phi$ be lower semicontinuous and $(\kappa_1,\kappa_2)$-Lipschitz in the large. Then $U(t,\mu):=P_t^-\phi(\mu)$ is a viscosity solution of the following Cauchy problem of \eqref{phje}
\begin{align*}
	\left\{
		\begin{array}{ll}
			\partial_tU(t,\mu)+\int_{\R^m}H\left(x,\partial_\mu U(t,\mu)\right)\,d\mu=0,&(t,\mu)\in\R^+\times\mathscr P_c(\R^m),\\
			U(0,\mu)=\phi(\mu),& \mu\in \mathscr P_c(\R^m).
		\end{array}
	\right.
\end{align*}
\end{The}

Corollary \ref{cor:HJ_ineq} indicates that, unlike the fundamental solution \( A_t(x,y) \), the dynamical cost functional \( C^t(\mu, \nu) \) can only possibly be a (strict) viscosity subsolution in general. This is because it can be regarded as a convex combination of \( A_t(x,y) \) based on measures \( \mu \) and \( \nu \). However, under the action of the random Lax-Oleinik operators, which corresponds to the convolution with the potential energy functional, it becomes a viscosity solution.

\begin{proof}
Following Theorem~\ref{thm:pt-}, we can express \( U(t, \mu) \) as follows:
\begin{align*}
	U(t,\mu):=P_t^-\phi(\mu)=\int_{\R^m}T_t^-\phi(x)\,d\mu=\int_{\R^+\times\R^m}u(r,x)\,d(\delta_t\times\mu).
\end{align*}
It should be noted that for any fixed \( R > 0 \), \( u(t, x) := T_t^- \phi(x) \) is a semiconcave function on the domain \( (0, +\infty) \times B(0, R) \). When \( \mu \in \mathscr{P}(B(0, R)) \), according to Remark~\ref{rem:time-dependent-D+}, \( (q, \alpha) \in \partial^+ U(t, \mu) \) if and only if \( (q, \alpha(x)) \in D^+ u(t, x) \) for almost every \( x \in \mathbb{R}^m \) with respect to the measure \( \mu \). Furthermore, since \( u(t, x) \) is the viscosity solution to equation \eqref{hje-}, we have
\begin{align*}
	q+H(x,\alpha(x))\leqslant 0,\,\,\,\,\mu-\mathrm{a.e.}\,x\in\R^m.
\end{align*}
By integrating over \( \mu \), we derive equation \eqref{eq:phje-sub}. Given the arbitrariness of \( (q, \alpha) \in \partial^+ U(t, \mu) \), it follows that \( U(t, \mu) \) is a viscosity subsolution to equation \eqref{phje}.

If \((p, \beta) \in \partial^- U(t, \mu)\), then according to Remark~\ref{rem:time-dependent-D+}, since \(\partial^+ U(t, \mu) \cap \partial^- U(t, \mu) \neq \varnothing\), \(u\) is differentiable on a set of full measure with respect to \(\mu \times \delta_t\), and we have
\begin{align*}
	D^+u(t,x)\cap D^-u(t,x)=\{(p,\beta(x))\},\,\,\,\,\mu-\mathrm{a.e.}\,x\in\R^m.
\end{align*}
Given that \(u(t, x)\) is the viscosity supersolution, it follows that
\begin{align*}
    p+H(x,\beta(x))\geqslant 0,\,\,\,\,\mu-\mathrm{a.e.}\,x\in\R^m.
\end{align*}
By integrating over \(\mu\), we obtain equation \eqref{eq:phje-sup}, which implies that \(U(t, \mu)\) is the viscosity supersolution to \eqref{phje}. Consequently, \(U(t, \mu)\) is the viscosity solution to \eqref{phje}. Furthermore, by Proposition~\ref{pro:probability lax-oleinik}, \(U(0, \cdot) = \phi(\cdot)\), thus \(U\) is the solution to the Cauchy problem for \eqref{phje}.
\end{proof}

\begin{Rem}
If $\phi$ is upper semicontinuous and $(\kappa_1,\kappa_2)$-Lipschitz in the large, $V(t,\mu):=P_t^+\phi(\mu)$ is the viscosity solution to the following Cauchy problem
\begin{align*}
	\left\{
		\begin{array}{ll}
			\partial_tV(t,\mu)-\int_{\R^m}H\left(x,\partial_\mu V(t,\mu)\right)\,d\mu=0,&(t,\mu)\in\R^+\times\mathscr P_c(\R^m),\\
			V(0,\mu)=\phi(\mu),& \mu\in \mathscr P_c(\R^m).
		\end{array}
	\right.
\end{align*}	
\end{Rem}

\begin{Pro}\label{pro:phjs}
The following statements are equivalent.
\begin{enumerate}[\rm (1)]
	\item $u$ is a weak KAM solution (viscosity solution) of \eqref{hjs} on $\T^m$;
	\item For any $t\geqslant 0$, $T_t^-u+c[0]t=u$;
	\item For any $t\geqslant 0$, $P_t^-u+c[0]t=u(\cdot)$;
	\item For any $t\geqslant 0$, $P_t^-\circ P_t^+u=u(\cdot)$.
\end{enumerate}	In this case, $u(\cdot)$ is the viscosity solution of the corresponding stationary type of Hamilton-Jacobi equation 
\begin{align}\tag{PHJ$_{s}$}\label{phjs}
	\int_{\T^m}H(x,\partial_\mu u(\mu)(x))\,d\mu=c[0],\,\,\,\,\mu\in\mathscr P(\T^m).
\end{align}
\end{Pro}

\begin{proof}
The equivalence of conditions (1) and (2) is discussed in \cite{Cannarsa_Cheng_Hong2025}. If condition (2) holds, for any \(\mu \in \mathscr{P}(\mathbb{T}^m)\), we have
\begin{align*}
	P_t^-u(\mu)+c[0]t=\int_{\T^m}T_t^-u+c[0]t\,d\mu=\int_{\T^m}u\,d\mu=u(\mu),
\end{align*}
which implies condition (3). On the other hand, let \(t \geqslant 0\) and \(\mu \in \mathscr{P}(\mathbb{T}^m)\), then \(P_t^- u(\mu) + c[0] t = u(\mu)\). By choosing \(\mu = \delta_x\), we get \(T_t^- u(x) + c[0] t = u(x)\). The arbitrariness of \(x \in \mathbb{T}^m\) confirms that condition (2) is satisfied. Due to Theorem \ref{thm:hjs} and condition (4) of Proposition \ref{pro:probability lax-oleinik}, condition (4) is equivalent to conditions (1)-(3).

For the final part of this proposition, we observe that \(P_t^- u(\mu) + c[0] t = u(\mu)\) and since \(P_t^- u(\mu)\) is the viscosity solution to equation \eqref{phje}, it follows that \(u(\cdot)\) is the viscosity solution to equation \eqref{phjs}.
\end{proof}

A simple fact is: if $u$ satisfies $u=T_t^-u+c[0]t$ for some $t\geqslant 0$, then $u=T_t^-\circ T_t^+u$ for such $t$. As a direct consequence of this fact, we have the following corollary.

\begin{Cor}\label{cor:wkam-fixed}
Assume that $u\in\Lip(\R^m)\cap\SCL(\R^m)$ is a viscosity solution of \eqref{hjs}. Then for such $u$, we can also get $T_t^-\circ T_t^+u=u$ for arbitrary $t\geqslant 0$. As a direct consequence,
\begin{align*}
	P_t^-u+c[0]t=u(\cdot),\,\,P_t^-\circ P_t^+u=u(\cdot),\,\,\forall t\geqslant 0.
\end{align*}
\end{Cor}

\section{Cut locus and singularity propagation}

In Section~\ref{subsec:potential-energy-diff-scl}, we identified an equivalent statement regarding the singular points of potential energy functionals induced by semiconcave functions. Subsequently, we introduce the concept of the cut locus of potential energy functionals induced by viscosity solutions to Hamilton-Jacobi equations in \(\mathscr{P}_c(\mathbb{R}^m)\), drawing inspiration from the definition of the cut locus of viscosity solutions. We then investigate the propagation of its singularities in the context of mass transport.

\subsection{Cut locus and cut time function}

In this subsection, we let $H:\R^{2m}\to\R$ be a Tonelli Hamiltonian and $u\in\Lip(\R^m)\cap\SCL(\R^m)$ be a viscosity solution of \eqref{hjs}
\begin{align*}
	H(x,Du(x))=c[0].
\end{align*}
Theorem \ref{thm:hjs} provides the existence of such solutions. 

\begin{defn}[Calibrated curve]\label{def:u-measure-cali}
An absolutely continuous curve\footnote{For the definition on absolutely continuous curves of metric space, readers can see \cite[Definition 1.1.1]{Ambrosio_GigliNicola_Savare_book2008}.} $\{\mu_s\}_{s\in I}\subset  \mathscr P_c(\R^m)$ defined on some interval $I\subset\R$ is called a \textit{$(u(\cdot),C,c[0])$-calibrated curve}, if for any $a,b\in I$, $a\leqslant b$, we have 
\begin{align*}
	u(\mu_b)-u(\mu_a)=C^{b-a}(\mu_a,\mu_b)+c[0](b-a).
\end{align*}
\end{defn}

\begin{Rem}\label{rem:wass-cali-exist}
Suppose $u$ is a viscosity solution.  For $\mu\in\mathscr P_c(\R^m)$, there is a $(u(\cdot),C,c[0])$-calibrated curve ending at $\mu$. This is a direct consequence of the following fact mentioned in \cite{Fathi_book}: for each $x\in\R^m$, there exists a $(u,L,c[0])$-calibrated curve $\eta\in C^2((-\infty,0],\R^m)$ with $\eta(0)=x$ in $\R^m$. In other words, for any $a\leqslant b\leqslant 0$, 
\begin{align*}
u(\eta(b))-u(\eta(a))=\int_a^bL(\eta(s),\dot\eta(s))\,ds.	
\end{align*}
Define $\mu=\law(x(\cdot))$. Then there exists a random curve $\xi=\xi(t)=\xi(t,\omega)$ such that $\xi(0,\omega)=x(\omega)$ and $\xi(\cdot,\omega)\in C^2((-\infty,0],\R^m)$ is a calibrated curve for each $\omega\in\Omega$. Suppose $\mu_s:=\law(\xi(s))$, then $\{\mu_s\}_{s\in(-\infty,0]}$ is we need.
\end{Rem}

\begin{Pro}
Given a $(u(\cdot),C,c[0])$-calibrated curve $\{\mu_s\}_{s\in [a,b]}\subset\mathscr P_c(\R^m)$. 
\begin{enumerate}[\rm (1)]
	\item $\{\mu_s\}_{s\in [a,b]}$ is a displacement interpolation of $C^{b-a}(\mu_a,\mu_b)$;
	\item for any subinterval $I'\subset [a,b]$, $\{\mu_s\}_{s\in I'}$ is also a $(u(\cdot),C,c[0])$-calibrated curve.
\end{enumerate}
\end{Pro}

\begin{proof}
According to Definition \ref{def:u-measure-cali}, for any $a\leqslant t_1<t_2<t_3\leqslant b$, we have that
\begin{align*}
	u(\mu_{t_2})-u(\mu_{t_1})&=C^{t_2-t_1}(\mu_{t_1},\mu_{t_2})+c[0](t_2-t_1),\\
		u(\mu_{t_3})-u(\mu_{t_2})&=C^{t_3-t_2}(\mu_{t_2},\mu_{t_3})+c[0](t_3-t_2),\\
        u(\mu_{t_3})-u(\mu_{t_1})&=C^{t_3-t_1}(\mu_{t_1},\mu_{t_3})+c[0](t_3-t_1),
\end{align*}
then $C^{t_3-t_1}(\mu_{t_1},\mu_{t_3})=C^{t_2-t_1}(\mu_{t_1},\mu_{t_2})+C^{t_3-t_2}(\mu_{t_2},\mu_{t_3})$. Due to the arbitrariness of $t_1,t_2$ and $t_3$, $\{\mu_s\}_{s\in[a,b]}$ is a displacement interpolation of $C^{b-a}(\mu_a,\mu_b)$. The rest of proposition is an immediate consequence of Definition \ref{def:u-measure-cali}.
\end{proof}

\begin{defn}[Cut point]\label{def:cut measure}
We call $\mu\in\mathscr P_c(\R^m)$ a \textit{cut point} of $u(\cdot)$, the potential energy functional induced by viscosity solution $u$, if all $(u(\cdot),C,c[0])$-calibrated curves ending at $\mu$ are unable to extend forward still as $(u(\cdot),C,c[0])$-calibrated curves. The set of all cut points of $u(\cdot)$ is denoted by $\mathscr C(u(\cdot))$, which is called the \textit{cut locus} of functional $u(\cdot)$.
\end{defn}

\begin{Lem}\label{lem:int-calibrated-diff}
$u(\cdot)$ is locally differentiable on the interior of $(u(\cdot),C,c[0])$-calibrated curve. Furthermore, $\mathscr S(u(\cdot))\subset\mathscr C(u(\cdot))$.
\end{Lem}

\begin{proof}
Let $I=[a,b]$ and $\{\mu_s\}_{s\in I}\subset  \mathscr P_c(\R^m)$ be a $(u(\cdot),C,c[0])$-calibrated curve, where $-\infty\leqslant a<b<+\infty$. Then for any $c,d\in[a,b]$ with $c\leqslant d$, we have
\begin{equation}\label{eq:phi-cali}
	u(\mu_d)-u(\mu_c)=C^{d-c}(\mu_c,\mu_d)+c[0](d-c).
\end{equation}
According to Corollary \ref{cor:partial+singleton}, $\partial^+u\neq\varnothing$ on $\mathscr P_c(\R^m)$. Thus, we only need to prove $\partial^-u(\mu_s)\neq\varnothing$ for all $s\in(a,b)$.
	
Note that $u$ is obviously a subsolution, which means
\begin{align*}
	u(y)\leqslant u(x)+A_t(x,y)+c[0]t,\,\,\,\,\forall x,y\in\R^m,\,\,t>0.
\end{align*}
Integrating the inequality above by $\gamma\in\Gamma_o^{b-s}(\mu,\mu_b)$ for any fixed $\mu\in\mathscr P_c(\R^m)$, we obtain
\begin{equation}\label{eq:phi-sub}
	u(\mu_b)-u(\mu)\leqslant C^{b-s}(\mu,\mu_b)+c[0](b-s).
\end{equation}
We set $d=b$ and $c=s$ in \eqref{eq:phi-cali}, then
\begin{equation}\label{eq:phi-cali-2}
	u(\mu_b)-u(\mu_s)=C^{b-s}(\mu_s,\mu_b)+c[0](b-s).
\end{equation}
By combining equations \eqref{eq:phi-sub} and \eqref{eq:phi-cali-2}, we derive
\begin{align*}
    U(\mu) := u(\mu_s) - u(\mu) &\leqslant C^{b-s}(\mu, \mu_b) - C^{b-s}(\mu_s, \mu_b) =: V(\mu),
\end{align*}
with equality when \(\mu = \mu_s\). Using Proposition \ref{pro:wass-supdiff-comparision}, we find that \(\partial^+ V(\mu_s) \subset \partial^+ U(\mu_s)\), which implies
\begin{align*}
    \partial^+ C^{b-s}(\mu, \mu_b) \bigg|_{\mu = \mu_s} &\subset \partial^+ (-u(\mu_s)) = -\partial^- u(\mu_s).
\end{align*}
According to Theorem \ref{thm:cost-supdiff}, \(-\mathbf{p}_{\mu_s} \in \partial^+ C^{b-s}(\mu_s, \mu_b) \subset -\partial^- u(\mu_s)\), hence \(\partial^- u(\mu_s) \neq \varnothing\). Since \(s \in (a, b)\) is arbitrary, \(u(\cdot)\) is locally differentiable on the interior of the \((u(\cdot), C, c[0])\)-calibrated curve. Consequently, if \(u(\cdot)\) is singular at \(\mu\), \(\mu\) is not in the interior of the \((u(\cdot), C, c[0])\)-calibrated curve, indicating that the \((u(\cdot), C, c[0])\)-calibrated curve ending at \(\mu\) cannot be extended further forward; that is, \(\mu\) is the cut point of \(u(\cdot)\).
\end{proof}

Inspired by the concept of cut time function, we introduce the \textit{cut time function of measure} with respect to the potential energy induced by viscosity solution $u$:
\begin{equation}\label{eq: cut time measure}
    \begin{split}
    	T_u(\mu):=&\,\sup\{t\geqslant 0:\exists\nu(\cdot)\in\mathrm{AC}([0,t];\mathscr P_c(\R^m)),\\
    	&\,\qquad\qquad\qquad\qquad\qquad \text{such that}\,u(\nu(t))=u(\mu)+C^t(\mu,\nu(t))+c[0]t\}.
    \end{split}
\end{equation}
Then $\mu\in\mathscr C(u(\cdot))$ if and only if $T_u(\mu)=0$. 

\begin{The}\label{thm: measure cut time}
Assume that $u$ is the viscosity solution to \eqref{hjs} on $\R^m$ and $\mu\in\mathscr P_c(\R^m)$. Then
\begin{align}
	T_u(\mu)&=\sup\{t\geqslant 0: P_t^-\circ P_t^+u(\mu)=P_t^+\circ P_t^-u(\mu)\}\label{eq: cut time measure 1}\\
		&=\sup\{t\geqslant 0:\mu(\{x\in\R^m:B_u(t,x)=0\})=1\}\label{eq: cut time measure 2}\\
		&=\inf\{\tau_u(x):x\in\supp(\mu)\}.\label{eq: cut time measure 3}
\end{align}
\end{The}

\begin{proof}
We first define the set on the right-hand side of equation \eqref{eq: cut time measure} as \( S(u, \mu) \). Note that \( u \) is a viscosity solution. Corollary \ref{cor:wkam-fixed} states that \( P_t^- \circ P_t^+ u = u \) and \( u = P_t^- u + c[0]t \) for any \( t \geqslant 0 \). Thus, \( t \) satisfies \( P_t^- \circ P_t^+ u(\mu) = P_t^+ \circ P_t^- u(\mu) \) if and only if
\begin{align}\label{eq:vis-commu}
    0 = P_t^- \circ P_t^+ u(\mu) - P_t^+ \circ P_t^- u(\mu) = u(\mu) - P_t^+ u(\mu) + c[0]t.
\end{align}
According to the definition of \( P_t^+ \), the Lipschitz property of displacement interpolation (Corollary \ref{cor:dis-inter}) and Corollary \ref{cor:pt+}, there exists a function \( \nu(\cdot) \in \mathrm{AC}([0, t]; \mathscr{P}_c(\mathbb{R}^m)) \) such that
\begin{align*}
    0 = u(\mu) - u(\nu(t)) + C^t(\mu, \nu(t)) + c[0]t,
\end{align*}
which implies that equation \eqref{eq:vis-commu} is equivalent to \( t \in S(u, \mu) \), i.e., equality \eqref{eq: cut time measure 1} holds.

Recall that the equality \( P_t^- \circ P_t^+ u(\mu) - P_t^+ \circ P_t^- u(\mu) = 0 \) can be expressed as
\begin{align}\label{eq:B(t,x)=0}
    \int_{\mathbb{R}^m} T_t^- \circ T_t^+ u(x) - T_t^+ \circ T_t^- u(x) \, d\mu = 0.
\end{align}
According to \cite[Lemma 3.1]{Cannarsa_Cheng_Hong2025}, for any \( t \geqslant 0 \), we have \( T_t^- \circ T_t^+ u \geqslant u \geqslant T_t^+ \circ T_t^- u \). This implies that equation \eqref{eq:B(t,x)=0} holds if and only if for almost every \( x \in \mathbb{R}^m \) with respect to \(\mu\), \( B_u(t, x) = 0 \). In other words, equality \eqref{eq: cut time measure 2} is satisfied.

We now proceed to the proof of equation \eqref{eq: cut time measure 3}. When \( t_0 < T_u(\mu) \), we have
\begin{align*}
    \mu(\{x \in \mathbb{R}^m : B_u(t_0, x) = 0\}) = 1.
\end{align*}
We assert that for any \( x_0 \in \supp(\mu) \), \( B_u(t_0, x_0) = 0 \). Indeed, if \( B_u(t_0, x_0) > 0 \), the lower semicontinuity of \( B_u(t_0, \cdot) \) implies that there exists \( R_0 > 0 \) such that \( B_u(t_0, x) > 0 \) for all \( x \in B(x_0, R_0) \). Since \( x_0 \in \supp(\mu) \), we have \( \mu(B(x_0, R_0)) > 0 \), which contradicts \( \mu(\{x \in \mathbb{R}^m : B_u(t_0, x) = 0\}) = 1 \). This assertion implies that \( t_0 \leqslant \inf\{\tau_u(x) : x \in \supp(\mu)\} \). By the arbitrariness of \( t_0 \), we conclude that \( T_u(\mu) \leqslant \inf\{\tau_u(x) : x \in \supp(\mu)\} \).

On the other hand, if \( T_u(\mu) < \inf\{\tau_u(x) : x \in \supp(\mu)\} \), there exists some \( t \in \mathbb{R} \) such that
\begin{align}\label{eq:cut time inf}
    T_u(\mu) < t < \inf\{\tau_u(x) : x \in \supp(\mu)\}.
\end{align}
For all \( x \in \supp(\mu) \), we have \( t < \tau_u(x) \), which implies \( B_u(t, x) = 0 \). Therefore,
\begin{align*}
    \supp(\mu) \subset \{x \in \mathbb{R}^m : B_u(t, x) = 0\},
\end{align*}
and it follows that
\begin{align*}
    1 = \mu(\supp(\mu)) \leqslant \mu(\{x \in \mathbb{R}^m : B_u(t, x) = 0\}),
\end{align*}
i.e., \( T_u(\mu) \geqslant t \), but this contradicts equation \eqref{eq:cut time inf}. In conclusion, we have \( T_u(\mu) = \inf\{\tau_u(x) : x \in \supp(\mu)\} \).
\end{proof}

Following the concept of the Aubry set as discussed in Section \ref{sub:cut locus}, the \textit{Aubry set} of the potential energy function \( u(\cdot) \) is defined as the collection of all measures \(\mu\) for which there exists a \((u(\cdot), C, c[0])\)-calibrated curve \(\nu(\cdot): [0, +\infty) \to \mathscr{P}(\mathbb{T}^m)\) such that \(\nu(0) = \mu\). Consequently, a measure \(\mu\) is in the Aubry set of \( u(\cdot) \) if and only if \( T_u(\mu) = +\infty \).

\begin{Cor}\label{cor:cut-prop}
Under the assumptions in Theorem \ref{thm: measure cut time},
\begin{enumerate}[\rm (1)]
	\item if $\supp(\mu)\cap\Cut(u)\neq\varnothing$, then $\mu\in\mathscr C(u(\cdot))$;
		\item if $\tau_u$ is continuous, then $\supp(\mu)\cap\Cut(u)\neq\varnothing$ if and only if $\mu\in\mathscr C(u(\cdot))$;
		\item $\mu$ is in the Aubry set of $u(\cdot)$ if and only if $\supp(\mu)\subset\mathcal I(u)$;
		\item if \(\mu = \law(Z)\) is the distribution of an \(\mathbb{R}^m\)-valued random variable \(Z : \Omega \to \mathbb{R}^m\), then \(\mu \in \mathscr{C}(u(\cdot))\) if and only if, for any \(\varepsilon > 0\),
		\begin{align*}
			\mathbb{P}\left( \{\omega \in \Omega : \tau_u(Z(\omega)) < \varepsilon \} \right) > 0.
		\end{align*}
\end{enumerate}
\end{Cor}

\begin{proof}
If the support of \(\mu\) intersects the cut locus of \(u\), i.e., \(\supp(\mu) \cap \Cut(u) \neq \varnothing\), then there exists an \(x \in \supp(\mu)\) such that \(\tau_u(x) = 0\). This implies that \(T_u(\mu) = 0\) by equation \eqref{eq: cut time measure 3}, and thus \(\mu \in \mathscr{C}(u(\cdot))\). For the proof of part (2), we only need to consider the ``if'' direction. Since \(\mu\) is a cut point of \(u(\cdot)\), equation \eqref{eq: cut time measure 3} gives us
\begin{align*}
    T_u(\mu) = \inf\{\tau_u(x) : x \in \supp(\mu)\} = 0.
\end{align*}
Together with the continuity of \(\tau_u\), there exists an \(x_0 \in \supp(\mu)\) such that \(\tau_u(x_0) = 0\), which means \(x_0 \in \supp(\mu) \cap \Cut(u)\).

By using $T_u$, $\mu$ is in the Aubry set of $u(\cdot)$ if and only if $T_u(\mu)=+\infty$. Then (3) is a direct consequence.

We now turn to prove (4). Given any fixed \(\varepsilon > 0\), we first note that \(\tau_u(Z(\omega)) < \varepsilon\) if and only if \(B_u(\varepsilon, Z(\omega)) > 0\) for \(\omega \in \Omega\). Due to this fact, we have
\[
\mathbb{P}(\{\omega : B_u(\varepsilon, Z(\omega)) > 0\}) = \mathbb{P}(\{\omega : \tau_u(Z(\omega)) < \varepsilon\}) > 0.
\]
Therefore, we can deduce the following:
\begin{align*}
	\mu(\{x\in\R^m:B_u(\varepsilon,x)=0\})&=1-\mu(\{x\in\R^m:B_u(\varepsilon,x)>0\})\\
	&\leqslant 1-\mu(\{x\in\R^m:\exists\omega\in\Omega,Z(\omega)=x,B_u(\varepsilon,x)>0\})\\
	&=1-\mathbb P(\{\omega\in\Omega:\tau_u(Z(\omega))<\varepsilon\})<1.
\end{align*}
This implies \(T_u(\mu) \leqslant \varepsilon\) by \eqref{eq: cut time measure 2}. Since \(\varepsilon\) is arbitrary, we conclude that \(T_u(\mu) = 0\). On the other hand, if there exists some \(\varepsilon_0 > 0\) such that \(\mathbb{P}(\{\omega \in \Omega : \tau_u(Z(\omega)) < \varepsilon_0\}) = 0\), i.e., \(B_u(\varepsilon_0, Z(\omega)) = 0\) for \(\mathbb{P}\)-a.e. \(\omega \in \Omega\), then we have
\[
\mu(\{x \in \mathbb{R}^m : B_u(\varepsilon_0, x) = 0\}) \geqslant \mathbb{P}(\{\omega \in \Omega : B_u(\varepsilon_0, Z(\omega)) = 0\}) = 1,
\]
which means \(T_u(\mu) \geqslant \varepsilon_0 > 0\). This contradicts the assumption that \(\mu \in \mathscr{C}(u(\cdot))\). This completes the proof.
\end{proof}

\begin{Lem}\label{lem:homeo-supp}
	Suppose $X,Y$ are Polish spaces, $\mu\in\mathscr P(X)$ and $f: X\to Y$ be a homeomorphism. Then $\supp(f_\#\mu)=f(\supp(\mu))$.
\end{Lem}

\begin{proof}
We denote the ball with radius \( R \geqslant 0 \) and center \( x \in X \) (resp. \( y \in Y \)) in \( X \) (resp. \( Y \)) as \( B_X(x, R) \) (resp. \( B_Y(y, R) \)). Given \( y \in f(\supp(\mu)) \), there exists a unique \( x \in \supp(\mu) \) such that \( y = f(x) \). Moreover, for any \( R > 0 \), we can choose some \( r > 0 \) such that \( f^{-1}(B_Y(y, R)) \), as a neighborhood of \( x \), contains \( B_X(x, r) \). In this case,
\begin{align*}
    f_\#\mu(B_Y(y, R)) = \mu(f^{-1}(B_Y(y, R))) \geqslant \mu(B_X(x, r)) > 0,
\end{align*}
which implies that \( y \in \supp(f_\#\mu) \).

Conversely, let \( y \in \supp(f_\#\mu) \). There exists a unique \( x \in X \) such that \( y = f(x) \). We need to verify that \( x \in \supp(\mu) \). For any \( R > 0 \), there exists some \( r > 0 \) such that \( f(B_X(x, R)) \supset B_Y(y, r) \), which means \( B_X(f^{-1}(y),R) \supset f^{-1}(B_Y(y, r)) \). This implies
\begin{align*}
    \mu(B_X(x, R)) = \mu(B_X(f^{-1}(y), R)) \geqslant \mu(f^{-1}(B_Y(y, r))) = f_\#\mu(B_Y(y, r)) > 0.
\end{align*}
By the arbitrariness of \( R > 0 \), it follows that \( x \in \supp(\mu) \) and \( y = f(x) \in f(\supp(\mu)) \). 
\end{proof}

\begin{The}
Given a viscosity solution $u$ of \eqref{hjs} on $\R^m$. Suppose $\mu\in\mathscr P_c(\R^m)$ and $T_u(\mu)>0$. Let $\mu_t:=[\pi\circ\Phi_H^t(\cdot,Du(\cdot))]_\#\mu$ for $t\in[0,T_u(\mu)]$.  
\begin{enumerate}[\rm (1)]
	\item $\{\mu_t\}_{t\in[0,T_u(\mu)]}$ is a $(u(\cdot),C,c[0])$-calibrated curve;
	\item $\mu=[\pi\circ\Phi_H^{-t}(\cdot,Du(\cdot))]_\#\mu_t$ for $t\in[0,T_u(\mu))$;
	\item If $t\in[0,T_u(\mu))$ and $\mu\ll\mathscr L^m$, then $\mu_t\ll\mathscr L^m$.
\end{enumerate}
\end{The}

In this case, the expression $\pi\circ\Phi^t_H(\cdot,Du(\cdot)):\R^m\rightarrow \R$ should be interpreted as follows: there exists a measurable set $A$ with $\mu(A)=1$ such that $u$ is differentiable at each $x\in A$ and $\pi\circ\Phi^t_H(\cdot,Du(\cdot))$ is defined on $A$. One can extend the map $\pi\circ\Phi^t_H(\cdot,Du(\cdot))$ to $\R^m$ as a map $\varphi_t(\cdot)$ such that $\varphi_t\vert_{\R^m\setminus A}=0$.
%
%
% \begin{align*}
%  \pi\circ\Phi^t_H(x,Du(x)):=\left\{
%  \begin{array}{ll}
%    \pi\circ\Phi^t_H(x,Du(x)),&x\in A,\\
%   0, &x\in \R^m\setminus A,
%  \end{array}
%  \right. 
% \end{align*}
% where $A$ is a measurable set with $\mu(A)=1$ such that $u$ is differentiable at each $x\in A$. 
 Since the set-valued map $x\rightsquigarrow D^+u$ is upper semicontinuous, $\pi\circ\Phi^t_H(\cdot,Du(\cdot))\in C(A)$ and $\varphi_t(\cdot)$ is measurable. Therefore, $[\varphi_t]_\#\mu $ is well-defined. Similarly, $[\pi\circ\Phi^{-t}_H(\cdot,Du(\cdot))]_\#\mu_t$ is also well-defined.

\begin{proof}
Assume \( T_u(\mu) > 0 \). According to Theorem \ref{thm: measure cut time}, for any \( x \in \supp(\mu) \), we have \( \tau_u(x) \geqslant T_u(\mu) > 0 \). By Proposition \ref{pro: repre-cali}, for each \( x \in \supp(\mu) \), there exists a unique calibrated curve starting from \( x \):
\begin{align*}
    \varphi_t(x) = \pi \circ \Phi_H^t(x, Du(x)), \quad t \in [0, T_u(\mu)], \quad x \in \supp(\mu).
\end{align*}
For any fixed \( t \in [0, T_u(\mu)] \), the map
\begin{align*}
    \varphi_t(\cdot): \supp(\mu) &\to \mathbb{R}^m \\
    x &\mapsto \pi \circ \Phi_H^t(x, Du(x))
\end{align*}
is well-defined on \( \supp(\mu) \). Moreover, since \( u \) is of class \( C^{1,1} \) on the interior of the \((u, L, c[0])\)-calibrated curve, \(\varphi_t(\cdot)\in \text{Lip}(\supp(\mu); \mathbb{R}^m) \). %From Appendix \ref{app:relatively continuous}, the measure \( \mu_t := [\pi \circ \Phi_H^t(\cdot, Du(\cdot))]_\#\mu \) is well-defined for \( t \in [0, T_u(\mu)] \). 
Then, for any \( 0 \leqslant s < t \leqslant T_u(\mu) \), combining with Proposition \ref{pro:sub-hjs}, we obtain
\begin{align*}
    C^{t-s}(\mu_s, \mu_t) + c[0](t-s) &\geqslant u(\mu_t) - u(\mu_s) \\
    &= \int_{\mathbb{T}^m} u(\varphi_t(x)) - u(\varphi_s(x)) \, d\mu \\
    &= \int_{\mathbb{T}^m} A_{t-s}(\varphi_s(x),\varphi_t(x)) \, d\mu + c[0](t-s) \\
    &\geqslant C^{t-s}(\mu_s, \mu_t) + c[0](t-s),
\end{align*}
which implies that \( \{\mu_t\}_{t \in [0, T_u(\mu)]} \) is a \((u(\cdot), C, c[0])\)-calibrated curve.

Furthermore, by invoking Proposition \ref{pro: repre-cali}, for any fixed \( t \in [0, T_u(\mu)) \), we have
\begin{align*}
    \pi \circ \Phi_H^{-t}(\pi \circ \Phi_H^t(x, Du(x)), Du(\pi \circ \Phi_H^t(x, Du(x)))) = \pi \circ \Phi_H^{-t}(\varphi_t(x), Du(\varphi_t(x))) = x.
\end{align*}
The map \( \varphi_t(\cdot) \) is continuous and injective, and thus a bi-Lipschitz homeomorphism from \( \supp(\mu) \) onto its image for each \( t \in [0, T_u(\mu)) \), by the compactness of \( \supp(\mu) \). That is, \( \pi \circ \Phi_H^{-t}(\cdot, Du(\cdot)) \) is the Lipschitz inverse of \( \varphi_t(\cdot) \) on \( \supp(\mu_t) \) due to Lemma \ref{lem:homeo-supp}. It follows that
\begin{align*}
    [\pi \circ \Phi_H^{-t}(\cdot, Du(\cdot))]_\# \mu_t = [\pi \circ \Phi_H^{-t}(\cdot, Du(\cdot))]_\# [\pi \circ \Phi_H^t(\cdot, Du(\cdot))]_\# \mu = \mu, \quad t \in [0, T_u(\mu)).
\end{align*}
This proves part (2). Note that \( \pi \circ \Phi_H^{-t}(\cdot, Du(\cdot)) \) is Lipschitz and maps Lebesgue zero measure sets into Lebesgue zero measure sets. Consequently, if \( \mu = [\pi \circ \Phi_H^{-t}(\cdot, Du(\cdot))]_\# \mu_t \ll \mathscr{L}^m \), then \( \mu_t \ll \mathscr{L}^m \) for \( t \in [0, T_u(\mu)) \). This completes the proof of part (3).
\end{proof}

\subsection{Propagation of singularities}

\begin{The}\label{thm:wass-argmax}
Assume that $\phi\in\Lip(\R^m)\cap\SCL_{\mathrm{loc}}(\R^m)$. Then for any $\mu\in\mathscr P_c(\R^m)$, there exists $t_{\phi,\mu}>0$, depending on $\phi$ and $\mu$, such that 
\begin{align*}
	\mathop{\arg\max}_{\nu\in\mathscr P(\R^m)}\{\phi(\nu)-C^t(\mu,\nu)\},
\end{align*}
is a singleton for all $t\in(0,t_{\phi,\mu}]$. Define
\begin{align*}
	\nu_{\phi,\mu}(t):=\left\{
	\begin{array}{ll}
	    \mu,&t=0,\\
		\mathop{\arg\max}\limits_{\nu\in\mathscr P(\R^m)}\{\phi(\cdot)-C^t(\mu,\cdot)\},&t\in(0,t_{\phi,\mu}],
	\end{array}
	\right.
\end{align*}
then $\lim_{t\to 0^+}\nu_{\phi,\mu}(t)=\mu$ in the sense of metric $W_p$ for any fixed $p\geqslant 1$.
\end{The}

\begin{Rem}\label{rem:independent-mu}
When considering the case of a compact manifold, for example, \( M = \mathbb{T}^m \), the time \( t_\phi = t_{\phi, \mu} \) mentioned in Theorem \ref{thm:wass-argmax} can be independent of \(\mu\), being determined only by \(\phi\).
\end{Rem}

\begin{proof}%[Proof of Theorem \ref{thm:wass-argmax}]
First, Theorem \ref{thm:pt-} asserts that if \(\mu \in \mathscr{P}_c(\mathbb{R}^m)\),
\begin{align*}
    \varnothing \neq \mathop{\arg\max}_{\nu \in \mathscr{P}(\mathbb{R}^m)} \{\phi(\nu) - C^t(\mu, \nu)\} \subset \mathscr{P}_c(\mathbb{R}^m).
\end{align*}
Then, for any \(\nu \in \mathop{\arg\max}_{\nu \in \mathscr{P}(\mathbb{R}^m)} \{\phi(\nu) - C^t(\mu, \nu)\}\) and \(\gamma \in \Gamma_o^t(\mu, \nu)\),
\begin{align*}
    P_t^+ \phi(\mu) = \int_{\mathbb{R}^m} T_t^+ \phi(x) \, d\mu = \int_{\mathbb{R}^{2m}} T_t^+ \phi(x) \, d\gamma = \int_{\mathbb{R}^{2m}} \phi(y) - A_t(x, y) \, d\gamma.
\end{align*}
Moreover, for arbitrary \(x, y \in \mathbb{R}^m\), \(T_t^+ \phi(x) \geqslant \phi(y) - A_t(x, y)\). Thus, \(T_t^+ \phi(x) = \phi(y) - A_t(x, y)\) holds for \(\gamma\)-a.e. \((x, y) \in \mathbb{R}^{2m}\). According to Lemma \ref{pro:funda-minimizer}, \(|x - y| \leq \lambda_\phi t\) for \(\gamma\)-a.e. \((x, y) \in \mathbb{R}^{2m}\). In addition, there exists \(R > 0\) such that \(\supp(\mu) \subset B(0, R)\) (determined by \(\mu\)), and we can find some \(t_R > 0\) such that when \(t \in (0, t_R)\),
\begin{align*}
    \supp(\nu) \subset B_{\lambda_\phi t}(\mu) \subset B(0, R).
\end{align*}

We denote \( C_R > 0 \) as the semiconcave constant of \( \phi \) on \( B(0, R) \). Let \( t \in (0, t_R) \) and assume that \( \nu_1, \nu_2 \in \mathscr{P}_c(\mathbb{R}^m) \) satisfy
\begin{align*}
    \nu_1, \nu_2 \in \mathop{\arg\max}_{\nu \in \mathscr{P}(\mathbb{R}^m)} \{\phi(\nu) - C^t(\mu, \nu)\}.
\end{align*}
Then, \( \nu_1, \nu_2 \in \mathscr{P}(B(0, R)) \). We have two optimal plans \( \gamma_i \in \Gamma_o^t(\mu, \nu_i) \), such that \( |x - y_i| \leqslant \lambda_\phi t \) for \( \gamma_i \)-a.e. \( (x, y_i) \in \mathbb{R}^{2m} \), \( i = 1, 2 \). For \( s \in [0, 1] \), we define
\begin{align*}
    \gamma_{x,1\to 2}^s &:= ((1-s)\pi_{x,y_1} + s\pi_{x,y_2})_\#\gamma, \\
    \gamma_{1\to 2}^s &:= ((1-s)\pi_{y_1} + s\pi_{y_2})_\#\gamma,
\end{align*}
where we reset \( \gamma \in \Gamma(\mu, \nu_1, \nu_2) \) and \( (\pi_{x,y_1})_\#\gamma = \gamma_1 \), \( (\pi_{x,y_2})_\#\gamma = \gamma_2 \). In this case, \( \gamma_{1\to 2}^s \) is the generalized geodesic from \( \nu_1 \) to \( \nu_2 \) with \( \mu \) as the base point, and \( \gamma_{x,1\to 2}^s \in \Gamma(\mu, \gamma_{1\to 2}^s) \). Furthermore, Theorem \ref{thm:semi-convex} implies the existence of \( t_{\lambda_\phi, R} > 0 \) and \( C_{\phi, \mu} > 0 \) such that when \( t \in (0, t_{\lambda_\phi, R}] \),
\begin{equation}\label{eq:cost-geo-striscl}
	(1-s)C^t(\mu, \nu_1) + sC^t(\mu, \nu_2) - C^t(\mu, \gamma_{1\to 2}^s) \geqslant \frac{C_{\phi, \mu}}{t}s(1-s)W_{2, \gamma}^2(\nu_1, \nu_2).
\end{equation}

On the other hand, \(\phi\) is semiconcave on \(B(0, R)\) with constant \(C_R\). Invoking Proposition \ref{pro:potential-scl}, we obtain
\begin{align}
    (1-s)\phi(\nu_1) + s\phi(\nu_2) - \phi(\gamma_{1\to2}^s) \leqslant C_R s(1-s) W_{2,\gamma}^2(\nu_1, \nu_2).\label{eq:phi-geo-scl}
\end{align}
We set \( t_{\phi,\mu} := \min\left\{ t_{\lambda_\phi,R}, \frac{C_{\phi,\mu}}{C_R} \right\} \). Combining this with equations \eqref{eq:cost-geo-striscl} and \eqref{eq:phi-geo-scl}, when \( t \in (0, t_{\phi,\mu}] \),
\begin{align*}
    &(1-s)[\phi(\nu_1) - C^t(\mu, \nu_1)] + s[\phi(\nu_2) - C^t(\mu, \nu_2)] - [\phi(\gamma_{1\to 2}^s) - C^t(\mu, \gamma_{1\to 2}^s)] \\
    &\leqslant \left( C_R - \frac{C_{\phi,\mu}}{t} \right) s(1-s) W_{2,\gamma}^2(\nu_1, \nu_2) < 0.
\end{align*}
This implies that the function \( s \mapsto \phi(\gamma_{1\to 2}^s) - C^t(\mu, \gamma_{1\to 2}^s) \) is strictly concave on \([0, 1]\), and thus has a unique maximizer. Therefore, \(\nu_1 = \nu_2\), which means that \(\mathop{\arg\max}_{\nu \in \mathscr{P}(\mathbb{R}^m)} \{\phi(\nu) - C^t(\mu, \nu)\}\) is a singleton.

For \(\phi \in \text{Lip}(\mathbb{R}^m)\), we have \((\kappa_1, \kappa_2) = (\text{Lip}(\phi), 0)\). Equation \eqref{eq:argmax-wasser} states
\begin{align*}
    \lim_{t \to 0^+} W_1(\mu, \nu_{\phi, \mu}(t)) = \lim_{t \to 0^+} \lambda_\phi t = 0.
\end{align*}
Note that by Theorem \ref{thm:pt-}, the support of \(\nu_{\phi, \mu}(t)\) is contained in the union of balls \(\cup_{x \in \supp(\mu)} B(x, \lambda_\phi t)\). For all \(t \in [0, t_{\phi, \mu}]\), the supports \(\supp(\nu_{\phi, \mu}(t))\) are contained within the same compact subset \(K\), which implies that the \(W_p\)-distances for all \(p \geqslant 1\) are equivalent on \(K\). Consequently, for any fixed \(W_p\)-metric, we have that $\lim_{t \to 0^+} \nu_{\phi, \mu}(t) = \mu$. This completes the proof.
\end{proof}

\begin{Rem}\label{rem:sing-propgation-repre}
When considering a compact manifold, such as \(\mathbb{T}^m\), the singularity curve \(\nu_{\phi,\mu}(t)\) mentioned in Theorem \ref{thm:sing-propogation} below can be more precisely characterized. In fact, by the application of Propositions \ref{pro:arnaud} and \ref{pro:intrinsic}, we have
\begin{align*}
    \nu_{\phi,\mu}(t) = [\pi \circ \Phi_H^t(\cdot, DT_t^+ \phi(\cdot))]_\#\mu, \quad t \in (0, t_{\phi,\mu}].
\end{align*}
\end{Rem}

\begin{Lem}\label{lem:Dphi}
Suppose \(t\geqslant0 \). Given \(\mu \in \mathscr{P}_c(\mathbb{R}^m)\), a locally differentiable point of \(u(\cdot)\), we define 
\begin{align*}
  \psi_s(x):=\left\{
  \begin{array}{ll}
   \pi\circ\Phi^{s-t}_H(x,Du(x)),&x\in A,\\
   0, &x\in \R^m\setminus A,
  \end{array}
  \right.
\end{align*}
where $A$ is a measurable set with $\mu(A)=1$ such that $u$ is differentiable at each $x\in A$.
%we define \(\varphi_s(x) := \pi \circ \Phi_H^{s-b}(x, Du(x))\) for all \(s \in [a, b]\) and for \(\mu\)-almost every \(x \in \mathbb{R}^m\). 
%Then \(\varphi_s\) is measurable on a set of full measure with respect to \(\mu\). 
Let \(\{\mu_s\}_{s \in [0,t]}\) be a subset of \(\mathscr{P}_c(\mathbb{R}^m)\), defined as \(\mu_s = (\psi_s)_\#\mu\). Then for any \(s \in [0,t]\),
\begin{align}\label{eq:phi-wass-cali}
    u(\mu_t) = u(\mu_s) + C^{t-s}(\mu_s, \mu_t) + c[0](t-s).
\end{align}
Thus, \(\{\mu_s\}_{s \in [0,t]}\) is a \((u(\cdot), C, c[0])\)-calibrated curve.
\end{Lem}

\begin{proof}
%Note that the set-valued map \( x \rightsquigarrow D^+\phi(x) \) is upper semicontinuous. Consequently, the continuity of \( D\phi \) and the measurability of \( \varphi_t \) on a set of full measure with respect to \(\mu\) can be derived from Corollary \ref{cor:partial+singleton} and the results presented in Appendix \ref{app:relatively continuous}.
For any \(t \geqslant 0\), we first note that \(\mu_t = \mu\) and
\begin{align*}
    u(\mu_t) - u(\mu_0) \leqslant C^{t}(\mu_0, \mu_t) + c[0]t.
\end{align*}
Then, there exists a random curve \(\xi\) such that \(\law(\xi(s)) = \mu_s\) for \(s \in [0,t]\). This allows us to analyze the problem by considering individual points. For a differentiable point \(x\) of \(u\), there exists a unique \((u, L, c[0])\)-calibrated curve \(\eta_x \in C^2((-\infty, t]; \mathbb{R}^m)\) satisfying \(\eta_x(t) = x\). For all \(t>0\), \(\eta_x\) is the minimizer of the following Bolza problem:
\begin{align*}
    \min_{\eta \in AC([0,t]; \mathbb{R}^m)} \left\{ u(\eta(0)) + \int_0^t L(\eta(s), \dot{\eta}(s)) \, ds : \eta(t) = x \right\}.
\end{align*}
Simultaneously, for \(s \in [0,t]\),
\begin{align}\label{eq:phi-classic-cali}
    u(\eta_x(t)) - u(\eta_x(s)) = A_{t-s}(\eta_x(s), \eta_x(t)) + c[0](t-s).
\end{align}
Therefore, \(p_x(t) := L_v(\eta_x(t), \dot{\eta}_x(t)) = Du(x)\) and \(\eta_x\) solves the Euler-Lagrange equation. By invoking the well-posedness of the Cauchy problem for ODEs, we claim that \(\eta_s(x) = \psi_s(x)\) holds for any \(s \in [0,t]\) and \(\mu\)-almost every \(x \in \mathbb{R}^m\). Then, \((\psi_s)_\#\mu_t = \law(\xi(s))\) is the displacement interpolation of \(C^{t}(\mu_0, \mu_t)\) and \(\mu_s = (\psi_s)_\#\mu_t\). Finally, integrating \eqref{eq:phi-classic-cali} by \(\mu_t\), we obtain \eqref{eq:phi-wass-cali}.
\end{proof}

\begin{The}\label{thm:sing-propogation}
Assume that \(\mu \in \mathscr{C}(u(\cdot))\). Let \(\nu_{u,\mu}(t)\) and \(t_{u,\mu} > 0\) be as mentioned in Theorem \ref{thm:wass-argmax}, where \(u(\cdot)\) is the potential energy functional induced by the viscosity solution \(u\) of equation \eqref{hjs}. Then, \(\nu_{u,\mu}(t) \in \mathscr{S}(u(\cdot))\) for all \(t \in (0, t_{u,\mu}]\).
\end{The}

\begin{proof}
We prove this argument by contradiction: for any fixed \( t \in (0, t_{u,\mu}] \), we assume that \( u(\cdot) \) is locally differentiable at \( \nu_{u,\mu}(t) \). By Corollary \ref{cor:pt+}, there exists \( \gamma_t \in \Gamma_o^t(\mu, \nu_{u,\mu}(t)) \) such that
\begin{align*}
    P_t^+ u(\mu) = \int_{\mathbb{R}^m} T_t^+ u(x) \, d\mu = \int_{\mathbb{R}^{2m}} [u(y) - A_t(x, y)] \, d\gamma_t.
\end{align*}
Thus, we have for \( \gamma_t \)-almost every \( (x, y) \in \mathbb{R}^{2m} \),
\begin{align*}
    u(x) - A_t(x, y) = T_t^+ u(x).
\end{align*}
Moreover, because \( u(\cdot) \) is locally differentiable at \( \nu_{u,\mu}(t) \) by our assumption, \( u \) is differentiable at \( \nu_{u,\mu}(t) \)-almost every \( y \in \mathbb{R}^m \). Fermat's rule implies
\begin{align*}
    0 \in D_y^+ (u(y) - A_t(x, y)) = Du(y) - D_y^- A_t(x, y), \quad \gamma_t \text{-almost every } (x, y) \in \mathbb{R}^{2m}.
\end{align*}
Therefore, \( D_y A_t(x, y) = Du(y) \) for \( \gamma_t \)-almost every \( (x, y) \in \mathbb{R}^{2m} \). According to Lemma \ref{lem:Dphi}, \( x = \psi_t(y) \), \( (\psi_t)_\# (\nu_{u,\mu}(t)) = \mu \) and
\begin{align*}
    [0, t] \ni s \mapsto (\psi_s)_\# \nu_{u,\mu}(t)
\end{align*}
is a \((u(\cdot), C, c[0])\)-calibrated curve starting from \( \mu \). That is, a \((u(\cdot), C, c[0])\)-calibrated curve ending with \( \mu \) can be extended further (for the existence of such a curve, see Remark \ref{rem:wass-cali-exist}). This conclusion contradicts the assumption that \( \mu \in \mathscr{C}(u(\cdot)) \).
\end{proof}

\begin{Rem}
Although when \( t > t_{u,\mu} \), the set \(\mathop{\arg\max}_{\nu \in \mathscr{P}(\mathbb{R}^m)} \{\phi(\nu) - C^t(\mu, \nu)\}\) is not necessarily a singleton, we can still establish that
\begin{align*}
    \mathop{\arg\max}_{\nu \in \mathscr{P}(\mathbb{R}^m)} \{u(\nu) - C^t(\mu, \nu)\} \subset \mathscr{S}(u(\cdot)).
\end{align*}
\end{Rem}

\begin{Rem}
Suppose $u$ is a viscosity solution of \eqref{hjs} on $\R^m$ and $\mu\ll\mathscr L^m$.
\begin{enumerate}[\rm (1)]
	\item Let \(0 < T_u(\mu) \leqslant +\infty\). For \(t \in [0, T_u(\mu))\), we have
	\begin{align*}
		\mu_t := [\pi \circ \Phi_H^t(\cdot, Du(\cdot))]_\#\mu = \mathop{\arg\max}_{\nu \in \mathscr{P}(\mathbb{R}^m)} \{u(\nu) - C^t(\mu, \nu)\} \ll \mathscr{L}^m.
	\end{align*}
	In fact, this case is consistent with the DiPerna-Lions theory (\cite{DiPerna_Lions1989}). More precisely, \(\mu_t\) is driven by the flow \(\mathbf X_1(t, x) := \pi \circ \Phi_H^t(x, Du(x))\) with the vector field
	\begin{align*}
		\mathbf{b}_1(t, x) = H_p(x, Du(x)),
	\end{align*}
	where it is defined. Since \(u\) is of class \(C^{1,1}\) on the interior of the \((u, L, c[0])\)-calibrated curve, \(\mathbf{b}_1\) satisfies certain regularity conditions.
	\begin{enumerate}
	    \item $\Lip(\mathbf b_1(t,\cdot))\in L^1([0,t])$;
	    \item $\mathbf b_1\in L^1([0,t];\Lip(\R^m;\R^m))$;
	    \item $\mathrm{div}(\mathbf b_1(t,\cdot))\in L^1([0,t];L^\infty(\R^m))$;
	    \item $\frac{|\mathbf b_1|}{1+|x|}\in L^1([0,t];L^\infty(\R^m))$.
	\end{enumerate}
	According to the DiPerna-Lions theory, there exists a unique regular Lagrangian flow associated with the vector field \(\mathbf{b}_1\).
	\item Given \( T_u(\mu) = 0 \), it follows that \(\mu \in \mathscr{C}(u(\cdot))\) and there exists a positive time \(t_{u,\mu}\) determined by Theorem \ref{thm:wass-argmax}. For \(t \in (0, t_{u,\mu}]\),
	\begin{align*}
		\nu_{u,\mu}(t) = \mathop{\arg\max}_{\nu \in \mathscr{P}(\mathbb{R}^m)} \{ u(\nu) - C^t(\mu, \nu) \} \centernot{\ll} \mathscr{L}^m
	\end{align*}
	are the singular points of \(u(\cdot)\). Specifically, for the case of \(\mathscr{P}(\mathbb{T}^m)\), due to Arnaud's theorem, \(\nu_{u,\mu}\) is driven by the local Lipschitz flow on \(t \in (0,t_{u}] \times \mathbb{T}^m\)
	\begin{align*}
		\mathbf X_2(t, x) := \pi \circ \Phi_H^t(x, DT_t^+ u(x)),
	\end{align*}
	which can be extended to \([0,\infty) \times \mathbb{T}^m\), for instance, as a \(t_u\)-periodic vector field \(\mathbf{b}_2(t, x)\). Note that \(\Lip(DT_t^+ u) \leqslant \frac{C}{t}\) (see \cite{CCHW2024}). Therefore, \(\mathbf{b}_2\) does not generally satisfy conditions (a)-(d) when the time step \(t_u \to 0^+\), implying that the DiPerna-Lions theory is not applicable directly in this case.
\end{enumerate}
In other words, the transportation of an absolutely continuous measure (with respect to the Lebesgue measure) along the characteristics before the cut time can ensure absolute continuity. However, after the cut time, it loses this absolute continuity and becomes a singular measure.
\end{Rem}

\subsection{Variational construction of strict generalized measure singular characteristics}

In the last subsection, we constructed a Lipschitz curve \(\nu_{u,\mu}(\cdot)\) on \([0, t_u)\) (Remark \ref{rem:independent-mu}), which can be extended to \([0, +\infty)\) and propagates singularities globally, provided the initial data is singular. However, this construction depends critically on the time step \(t_u\). To obtain a semiflow evolution (at least for dynamical cost induced by mechanical Lagrangian) on the cut locus, we need to let \(t_u \to 0^+\). This idea is inspired by the concept of \emph{minimizing movements}, which is widely used in the theory of abstract gradient flows.

\begin{Lem}[\cite{CCHW2024}]\label{lem:direc-diff}
Suppose $M$ is a Riemannian manifold, $\phi\in\SCL_{\mathrm{loc}}(M)$ and $\eta:\R\to M$ is a locally absolutely continuous curve, then
\begin{align*}
	\frac{d}{dt}\phi(\eta(t))=p(\dot\eta(t)),\,\,\text{a.e. }t\in\R,\,\,\forall p\in D^+\phi(\eta(t)).
\end{align*}
\end{Lem}

\begin{The}\label{thm:conti-eq}
Suppose \(\phi \in \SCL(\mathbb{T}^m)\), \(T > 0\), and \(\mu_0 \in \mathscr{P}(\mathbb{T}^m)\). Then there exists a Lipschitz curve \(\mu(\cdot): [0, T] \to \mathscr{P}(\mathbb{T}^m)\) satisfying the following continuity equation:
\begin{equation}\label{eq:gene-sing-continuity}
\left\{
	\begin{array}{ll}
		\frac{d}{dt}\mu + \mathrm{div}(H_p(\cdot, \mathbf{p}_\phi^\#(\cdot)) \cdot \mu) = 0, \\
		\mu(0) = \mu_0,
	\end{array}
\right.
\end{equation}
where \(\mathbf{p}_\phi^\#(x) := \arg\min\{H(x, p): p \in D^+\phi(x)\}\) for any \(x \in \mathbb{T}^m\). 
\end{The}

\begin{proof}
Given \( T > 0 \) and \( \mu_0 \in \mathscr{P}(\mathbb{T}^m) \), let
\begin{align*}
    \Delta: 0 = t_0 < t_1 < \cdots < t_{N-1} < t_N = T
\end{align*}
be a partition of the interval \([0, T]\) with \(\|\Delta\| < t_\phi\), where \( t_\phi \) is as determined in Remark \ref{rem:independent-mu}. For \( N \in \mathbb{N} \) and \( k = 0, 1, \cdots, N-1 \), we define \(\mu_{k+1}\) inductively as
\begin{align*}
    \mu_{k+1} := \mathop{\arg\max}\{\phi(\nu) - C^{t_{k+1} - t_k}(\mu_k, \nu) : \nu \in \mathscr{P}(\mathbb{T}^m)\}.
\end{align*}
More precisely,
\begin{align*}
    P_{t_{k+1} - t_k}^+ \phi(\mu_k) = \phi(\mu_{k+1}) - C^{t_{k+1} - t_k}(\mu_k, \mu_{k+1}), \quad k = 0, 1, \cdots, N-1.
\end{align*}
In this case, by Corollary \ref{cor:pt+}, we have the following sequence of equalities:
\begin{equation}\label{eq:inductive}
\begin{aligned}
	\phi(\mu_{k+1}) - \phi(\mu_k) &= P_{t_{k+1}-t_k}^+\phi(\mu_k) - \phi(\mu_k) + C^{t_{k+1}-t_k}(\mu_k,\mu_{k+1}) \\
	&= \int_{\mathbb{T}^m} [T_{t_{k+1}-t_k}^+\phi(x) - \phi(x)] \, d\mu_k + C^{t_{k+1}-t_k}(\mu_k,\mu_{k+1}) \\
	&= \int_{\mathbb{T}^m} \int_{t_k}^{t_{k+1}} \frac{d}{ds} T_{s-t_k}^+\phi(x) \, ds \, d\mu_k + C^{t_{k+1}-t_k}(\mu_k,\mu_{k+1}) \\
	&= \int_{\mathbb{T}^m} \int_{t_k}^{t_{k+1}} H(x, DT_{s-t_k}^+\phi(x)) \, ds \, d\mu_k + C^{t_{k+1}-t_k}(\mu_k,\mu_{k+1}) \\
	&= \int_{t_k}^{t_{k+1}} \int_{\mathbb{T}^m} H(x, DT_{s-t_k}^+\phi(x)) \, d\mu_k \, ds + C^{t_{k+1}-t_k}(\mu_k,\mu_{k+1}).
\end{aligned}
\end{equation}
Summing this equation from \(k=0\) to \(k=N-1\), we obtain
\begin{align}\label{eq:inductive-2}
	\phi(\mu_N) - \phi(\mu_0) = \sum_{k=0}^{N-1} \left[ \int_{t_k}^{t_{k+1}} \int_{\mathbb{T}^m} H(x, DT_{s-t_k}^+\phi(x)) \, d\mu_k \, ds + C^{t_{k+1}-t_k}(\mu_k,\mu_{k+1}) \right].
\end{align}
It is noteworthy that for each \( k \), Remark \ref{rem:sing-propgation-repre} gives us
\begin{align*}
	\mu_{k+1} = [\pi \circ \Phi_H^{t_{k+1} - t_k}(\cdot, DT_{t_{k+1} - t_k}^+ \phi(\cdot))]_\# \mu_k.
\end{align*}
In this framework, we define
\begin{align}\label{eq:mu-delta}
    \mu_\Delta(s) &:= [\pi \circ \Phi_H^{s - t_k}(\cdot, DT_{t_{k+1} - t_k}^+ \phi(\cdot))]_\# \mu_k, \quad s \in [t_k, t_{k+1})
\end{align}
for \( k = 0, 1, \cdots, N-1 \) with \(\mu_\Delta(0):= \mu_0\). Then 
$\{\mu_{\Delta}(s)\}_{s\in[t_k,t_{k+1}]}$ is a displacement interpolation of $C^{t_{k},t_{k+1}}(\mu_k,\mu_{k+1})$ by Proposition \ref{pro:intrinsic}. Moreover, the family $\{\mu_\Delta:\|\Delta\|\ll 1\}$ is equi-Lipschitz (with a constant depending only on $L$ and $\phi$) in the sense of $W_1$ metric (and $W_p$ metric for all $p\geqslant 1$). If we set  $\tau_\Delta(\cdot): [0, T] \to [0, T]$  as  $\tau_\Delta(s) = t_k$  for  $s \in [t_k, t_{k+1})$ ,  $k = 0, 1, \ldots, N-1$ , and  $\tau_\Delta(T) = T$, then \eqref{eq:inductive-2} can be rewritten as 
\begin{equation}\label{eq:inductive-3}
\begin{split}
	\phi(\mu_N) - \phi(\mu_0) &= \int_{0}^{T} \int_{\mathbb{T}^m} H(x, DT_{s-\tau_\Delta(s)}^+\phi(x)) \, d\mu_{\Delta}(\tau_\Delta(s)) \, ds +\sum_{k=0}^{N-1} \A^{t_{k+1}-t_k}(\mu_\Delta)\\
	&= \int_{0}^{T} \int_{\mathbb{T}^m} H(x, DT_{s-\tau_\Delta(s)}^+\phi(x)) \, d\mu_{\Delta}(\tau_\Delta(s)) \, ds + \A^{T}(\mu_\Delta).\\
	%&=\int_{0}^{T} \int_{\mathbb{T}^m} H(x, DT_{s-\tau_\Delta(s)}^+\phi(x)) \, d\mu_{\tau_\Delta(s)} \, ds + \min_{\xi}\E\left(\int_0^TL(\xi(s),\dot\xi(s))\,ds\right),
\end{split}
\end{equation}
Now, consider a sequence of partitions  $\{\Delta_n\}_{n \geqslant 1}$  such that  $\lim_{n \to \infty} \|\Delta_n\| = 0$ . By invoking the Arzelà-Ascoli theorem (\cite[Theorem 10.3 \& 10.4]{Ambrosio_Brue_Semola_book2021}), we assume that  $\mu_{\Delta_n} \to \mu \in \text{Lip}([0, T]; \mathscr{P}(\mathbb{T}^m))$  uniformly. Also, since  $\tau_{\Delta_n}(t) \to t$  as  $\|\Delta_n\| \to 0$ ,  $\mu_{\Delta_n}(\tau_{\Delta_n}(t)) \mathrel{\ensurestackMath{\stackon[1pt]{\rightharpoonup}{\scriptstyle\ast}}} \mu(t)$  for any  $t \in [0, T]$ . Furthermore, from (\cite[Proposition 4.2]{CCHW2024}), we recall that for any  $x \in \mathbb{T}^m$,
\begin{align}\label{eq:limit-argmin}
    \lim_{t \to 0^+} DT_t^+ \phi(x) = \mathop{\arg\min}_{p \in D^+\phi(x)} \{H(x, p)\} =: \mathbf{p}_\phi^\#(x).
\end{align}
Equation \eqref{eq:inductive-3} and \eqref{eq:limit-argmin}, along with Fatou's lemma, imply that
\begin{equation}\label{eq:conti-geq}
\begin{split}
	\phi(\mu(T)) - \phi(\mu(0))&=\lim_{n\to\infty}  \int_{0}^{T} \int_{\mathbb{T}^m} H(x, DT_{s-\tau_{\Delta_n}(s)}^+\phi(x)) \, d\mu_{\Delta_n}(\tau_{\Delta_n}(s)) \, ds + \A^{T}(\mu_{\Delta_n})\\
	&\geqslant \int_{0}^{T} \int_{\mathbb{T}^m} H(x, \mathbf p_\phi^\#(x)) \, d\mu(s)\, ds + \A^{T}(\mu)\\
	&=\int_{0}^{T} \int_{\mathbb{T}^m} H(x, \mathbf p_\phi^\#(x)) \, d\mu(s)\, ds + \min_{\xi}\E\left(\int_0^TL(\xi(s),\dot\xi(s))\,ds\right)
\end{split}
\end{equation}
where the last minimum is over all absolutely continuous random curves $\xi:[0,T]\times\Omega\to M$ such that $\mu(s)=\law(\xi(s))$ for $s\in[0,T]$. We choose one of such minimizers and denote it by $\xi_*$. Now we can obtain from \eqref{eq:conti-geq} that
\begin{align}\label{eq:conti-eq}
	\E(\phi(\xi(T))-\phi(\xi(0)))\geqslant\E\left(\int_0^T[H(\xi_*(s), \mathbf p_{\phi}^\#(\xi_*(s)))+L(\xi_*(s),\dot\xi_*(s))]\,ds\right)
\end{align}
by Fubini's theorem. Due to Lemma \ref{lem:direc-diff} and Fenchel's inequality, \eqref{eq:conti-eq} is actually an equality. Furthermore, $\mu(t)$ is the law at time $t$ of a random solution $\xi_*$ of the following equation:
\begin{align*}
	\dot {\mathbf x}(t)=H_p(\mathbf x(t),\mathbf p_\phi^\#(\mathbf x(t))). 
\end{align*}In other words, $\xi_*$ satisfies
\begin{align*}
	\dot{\xi}_*(t, \omega) = H_p(\xi_*(t, \omega), \mathbf{p}_{\phi}^\#(\xi_*(t, \omega))), \quad \mathscr{L}^1\text{-a.e. } t \in [0, T], \quad \mathbb{P}\text{-a.e. } \omega \in \Omega.
\end{align*}
Consequently, \(\mu(\cdot): [0, T] \to \mathscr{P}(\mathbb{T}^m)\) satisfies \eqref{eq:gene-sing-continuity}. This completes the proof.
\end{proof}

\begin{Cor}\label{cor:conti-eq-1}
Assume that \(\phi=u\) in \eqref{eq:gene-sing-continuity} is a viscosity solution to \eqref{hjs}, , and \(\mu(\cdot): [0, T] \to \mathbb{T}^m\) is the solution of \eqref{eq:gene-sing-continuity}.
\begin{enumerate}[\rm (1)]
	\item If \(\mu_0(\overline{\Sing(u)}) > 0\), then \([\mu(t)](\overline{\Sing(u)}) > 0\) for all \(t \in [0, T]\).
	\item If \(\mu_0 \in \mathscr{C}(u(\cdot))\), then \(\mu(t) \in \mathscr{C}(u(\cdot))\) for all \(t \in [0, T]\).
\end{enumerate}
\end{Cor}

\begin{proof}
For the proof of (1), recall Theorem 1.2 in \cite{Albano2014_1}. If \(\xi(0, \omega) \in \overline{\Sing(u)}\) for a fixed \(\omega \in \Omega\), then \(\xi_*(t, \omega) \in \overline{\Sing(u)}\) for all \(t \in [0, T]\), where $\xi_*$ is the minimizer as in the proof of Theorem \ref{thm:conti-eq}. Thus,
\[
\{\omega \in \Omega : \xi_*(0, \omega) \in \overline{\Sing(u)}\} \subset \{\omega \in \Omega : \xi_*(t, \omega) \in \overline{\Sing(u)}\}, \quad \forall t \in [0, T].
\]
Since \(\mu_0(\overline{\Sing(u)}) > 0\), for any \(t \in [0, T]\), we have
\begin{align*}
		[\mu(t)](\overline{\Sing(u)})&=\mathbb P(\{\omega\in\Omega:\xi_*(t,\omega)\in\overline{\Sing(u)}\})\\
		&\geqslant\mathbb P(\{\omega\in\Omega:\xi_*(0,\omega)\in\overline{\Sing(u)}\})\\
		&		=\mu_0(\overline{\Sing(u)})>0.
\end{align*}

Now we turn to the proof of (2). Since \(\mu(\cdot): [0, T] \to \mathbb{T}^m\) is the solution of \eqref{eq:gene-sing-continuity}, we know that \(\mu(t)\) is a law at time $t$ of a random solution $\xi_*$ of the equation
\[
    \dot\xi_*(t, \omega) = H_p(\xi_*(t, \omega), \mathbf{p}_\phi^\#(\xi_*(t, \omega))) \in \text{co} \, H_p(\xi_*(t, \omega), D^+ \phi(\xi_*(t, \omega)))
\]
for a.e. \(t \in [0, T]\) and \(\mathbb{P}\)-a.e. \(\omega \in \Omega\). Since \(\mu_0 \in \mathscr{C}(\phi(\cdot))\), Corollary \ref{cor:cut-prop} implies that
\[
    \mathbb{P}(\{\omega \in \Omega : \tau_u(\xi(0, \omega)) < \varepsilon\}) > 0, \quad \forall \varepsilon > 0.
\]
Moreover, according to \cite[Theorem 5.7]{CCHW2024}, \(\tau_u(\xi_*(t, \omega)) \leqslant \tau_u(\xi_*(0, \omega))\) for all \(t \in [0, T]\) and \(\mathbb{P}\)-a.e. \(\omega \in \Omega\). Hence, for any \(\varepsilon > 0\),
\[
    \mathbb{P}(\{\omega \in \Omega : \tau_u(\xi_*(t, \omega)) < \varepsilon\}) \geqslant \mathbb{P}(\{\omega \in \Omega : \tau_u(\xi_*(0, \omega)) < \varepsilon\}) > 0.
\] 
This shows that \(\mu(t) \in \mathscr{C}(u(\cdot))\) by Corollary \ref{cor:cut-prop}.
\end{proof}

%\begin{Cor}\label{cor:conti-eq-2}
%Let \( H(x,p) = \frac{1}{2} \langle A(x) p, p \rangle + V(x) \) be a time-independent mechanical Hamiltonian. Suppose \(\mu(\cdot): [0,T] \to \mathbb{T}^m\) is the solution of \eqref{eq:gene-sing-continuity}. If \(\phi\) is a viscosity solution to \eqref{hjs} and \(\mu_0 \in \mathscr{S}(\phi(\cdot))\) in \eqref{eq:gene-sing-continuity}, then \(\mu(t) \in \mathscr{S}(\phi(\cdot))\) for all \(t \in [0,T]\).
%\end{Cor}
%
%\begin{proof}
%In the proof of Theorem \ref{thm:conti-eq}, we have \(\mu(t) = \law(\xi(t,\cdot))\) and
%\begin{align*}
%	\dot{\xi}(t,\omega) = A(\xi(t,\omega)) \mathbf{p}_{\phi}^\#(\xi(t,\omega)), \quad \mathscr{L}^1-\text{a.e. } t \in [0,T], \quad \mathbb{P}-\text{a.e. } \omega \in \Omega.
%\end{align*} 
%By \cite[Theorem 1.1]{Albano2016_1}, if \(\xi(0,\omega) \in \Sing(\phi)\) for a fixed \(\omega \in \Omega\), then \(\xi(t,\omega) \in \Sing(\phi)\) for all \(t \in [0,T]\). 
%The rest of the proof follows similarly to the last part of the proof of Theorem \ref{thm:conti-eq}.
%\end{proof}

\begin{Rem}
Given \(\phi \in \SCL(\T^m)\), we define\footnote{This definition can be extended similarly to the case of Euclidean space \(\R^m\).}
\begin{align*}
	\overline{\mathscr S}(\phi(\cdot)) := \left\{ \mu \in \mathscr P(\T^m) : \mu(\overline{\Sing(\phi)}) > 0 \right\}.
\end{align*}

\begin{enumerate}[(1)]
	\item \(\overline{\mathscr S}(\phi(\cdot))\) is an \(F_\sigma\)-set in the sense of weak-$\ast$ topology. In fact, we can consider the sets
	\[
		A_k := \left\{ \mu \in \mathscr P(\T^m) : \mu(\overline{\Sing(\phi)})\geqslant\frac{1}{k} \right\}, \quad k \in \mathbb{N}.
	\]
	By the Portmanteau theorem, for any sequence \(\{\mu_n\}_{n \in \mathbb{N}}\subset A_k\) with \(\mu_n \mathrel{\ensurestackMath{\stackon[1pt]{\rightharpoonup}{\scriptstyle\ast}}} \mu\),
	\[
		\mu_n(\overline{\Sing(\phi)})\geqslant\frac 1k \Rightarrow\ \mu(\overline{\Sing(\phi)}) \geqslant \limsup_{n \to \infty} \mu_n(\overline{\Sing(\phi)})\geqslant\frac{1}{k},
	\]
	i.e., \(\mu \in A_k\). Therefore, each \(A_k\) is a closed subset, and \(\overline{\mathscr S}(\phi(\cdot)) = \bigcup_{k=1}^\infty A_k\) is an \(F_\sigma\)-set.
	
	\item Suppose \(u\) is a viscosity solution of \eqref{hjs}. Then $\Sing(u)\subset\Cut(u)\subset\overline{\Sing(u)}$ (see \cite{CCHW2024}). However, \(\mathscr C(u(\cdot))\) and \(\overline{\mathscr S}(u(\cdot))\) are not mutually contained in general, and we have \(\mathscr S(u(\cdot)) \subsetneq \mathscr C(u(\cdot)) \cap \overline{\mathscr S}(u(\cdot))\). 
%	Theorem \ref{thm:conti-eq}, Corollary \ref{cor:conti-eq-1}, and Corollary \ref{cor:conti-eq-2} imply that \(\mathscr C(\phi(\cdot)) \cap \overline{\mathscr S}(\phi(\cdot))\) is a positively invariant set for the semi-flow induced by \eqref{eq:gene-sing-continuity} in \(\mathscr P(\T^m)\).
\end{enumerate}
\end{Rem}

\begin{Lem}[\cite{CCHW2024}]\label{lem:right-derivative}
Given \( T > 0 \). Let \( M \) be a closed manifold, \(\phi \in \SCL(M)\), and \( H \in C^2(T^*M) \) be a Tonelli Hamiltonian. If \(\eta: [0, T] \to M\) is a strict singular characteristic, i.e., \(\eta\) satisfies the equality in the following variational inequality
\begin{align*}
    \phi(\eta(t_2)) - \phi(\eta(t_1))\leqslant \int_{t_1}^{t_2} L(\eta(s), \dot{\eta}(s)) + H(\eta(s), \mathbf{p}_\phi^\#(\eta(s))) \, ds, \quad \forall t_1, t_2 \in [0, T], \, t_1 < t_2,
\end{align*}
or equivalently,
\begin{align*}
    \phi(\eta(T)) - \phi(\eta(0)) &\leqslant \int_{0}^{T} L(\eta(s), \dot{\eta}(s)) + H(\eta(s), \mathbf{p}_\phi^\#(\eta(s))) \, ds,
\end{align*}
then
\begin{align*}
    \dot{\eta}^+(t) = H_p(\eta(t), \mathbf{p}_\phi^\#(\eta(t))), \quad \forall t \in [0, T).
\end{align*}
\end{Lem}

\begin{The}\label{thm:right_derivative}
Under the same assumptions as in Theorem \ref{thm:conti-eq}, the Lipschitz curve mentioned in Theorem \ref{thm:conti-eq} that satisfies equation \eqref{eq:gene-sing-continuity} also has the following property: for any \(f \in C^\infty(\mathbb{T}^m)\), the map \(t \mapsto f(\mu(t))\) is in \(\text{\rm Lip}\,([0, T];\R)\) and its right derivative
\begin{align*}
    \frac{d^+}{dt} f(\mu(t)) = \int_{\mathbb{T}^m} \langle Df(x), H_p(x, \mathbf{p}_\phi^\#(x)) \rangle \, d\mu(t)
\end{align*}
exists for all \(t \in [0, T)\).
\end{The}

\begin{proof}
%Combining equations \eqref{eq:strict-chara-geq} and \eqref{eq:strict-chara-leq}, we obtain
%\begin{align*}
%    \int_\Omega \phi(\xi(T, \omega)) \, d\mathbb{P} - \int_\Omega \phi(\xi(0, \omega)) \, d\mathbb{P} = \int_\Omega \int_0^T \left[ H(\xi(t, \omega), \mathbf{p}_\phi^\#(\xi(t, \omega))) + L(\xi(t, \omega), \dot{\xi}(t, \omega)) \right] dt \, d\mathbb{P}.
%\end{align*}
%Since for any \(\omega \in \Omega\),
%\begin{align*}
%    \phi(\xi(T, \omega)) - \phi(\xi(0, \omega)) \leqslant \int_0^T \left[ L(\xi(s, \omega), \dot{\xi}(s, \omega)) + H(\xi(s, \omega), \mathbf{p}_\phi^\#(\xi(s, \omega))) \right] ds,
%\end{align*}
%it follows that for \(\mathbb{P}\)-almost every \(\omega \in \Omega\), the inequality above is actually an equality. 
By Lemma \ref{lem:right-derivative}, we can deduce from the proof of Theorem \ref{thm:conti-eq} that
\begin{align*}
    \dot{\xi}_*^+(t, \omega) = H_p(\xi_*(t, \omega), \mathbf{p}_\phi^\#(\xi_*(t, \omega))), \quad \forall t \in [0, T), \, \mathbb{P}-\text{a.e.} \, \omega \in \Omega.
\end{align*}
Now, for an arbitrary fixed \(f \in C^\infty(\mathbb{T}^m)\), consider the map \(F(t) := f(\mu(t))\) with \(t \in [0, T]\). According to Proposition 16.3 in \cite{Ambrosio_Brue_Semola_book2021} and Theorem \ref{thm:conti-eq}, \(F \in\text{Lip}([0,T];\R)\). Moreover,
\begin{align*}
    \lim_{\tau \downarrow t} \frac{F(\tau) - F(t)}{\tau - t} &= \lim_{\tau \downarrow t} \frac{1}{\tau - t} \left\{ \int_{\mathbb{T}^m} f(x) \, d[\mu(\tau)] - \int_{\mathbb{T}^m} f(x) \, d[\mu(t)] \right\} \\
    &= \lim_{\tau \downarrow t} \frac{1}{\tau - t} \E\left(  f(\xi_*(\tau)) - f(\xi_*(t)) \, \right) \\
    &= \E\left(\lim_{\tau \downarrow t} \frac{f(\xi_*(\tau)) - f(\xi_*(t))}{\tau - t}\right)\\
    &= \E\left(\langle Df(\xi_*(t)), \dot{\xi}_*^+(t) \rangle \right) \\
    &= \E\left(\langle Df(\xi_*(t)), H_p(\xi_*(t), \mathbf{p}_\phi^\#(\xi_*(t))) \rangle\right) \\
    &= \int_{\mathbb{T}^m} \langle Df(x), H_p(x, \mathbf{p}_\phi^\#(x)) \rangle \, d[\mu(t)].
\end{align*}
That is, \(\dot{F}^+(t) = \int_{\mathbb{T}^m} \langle Df(x), H_p(x, \mathbf{p}_\phi^\#(x)) \rangle \, d[\mu(t)]\) for any \(t \in [0, T)\). Since the argument above holds for any \(f \in C^\infty(\mathbb{T}^m)\), we have completed the proof.
\end{proof}

\appendix

\section{Proofs of some statements}

\begin{proof}[Proof of Proposition \ref{pro: repre-cali}]
According to the definition of the cut time function, for \( x \in M \) with \( \tau_u(x) > 0 \), there exists a curve \( \eta_x: [0, \tau_u(x)] \to M \) such that
\begin{align*}
    u(\eta_x(\tau_u(x))) - u(x) = A_{\tau_u(x)}(x, \eta_x(\tau_u(x))) + c[0] \tau_u(x).
\end{align*}
This implies that
\begin{align*}
    T_{\tau_u(x)}^+ u(x) - c[0] \tau_u(x) &\geqslant u(\eta_x(\tau_u(x))) - A_{\tau_u(x)}(x, \eta_x(\tau_u(x))) -c[0]\tau_u(x)= u(x) \\
    &\geqslant T_{\tau_u(x)}^+ u(x) - c[0] \tau_u(x).
\end{align*}
Thus, we have
\begin{align*}
    T_{\tau_u(x)}^+ u(x) - c[0] \tau_u(x) = u(x)
\end{align*}
and
\begin{align*}
    D T_{\tau_u(x)}^+ u(x) = D u(x).
\end{align*}
The term \( A_{\tau_u(x)}(x, \eta_x(\tau_u(x))) \) has a unique minimizer, which we denote by \( \eta_x \). Therefore, \( \eta_x \) is a \((u, L, c[0])\)-calibrated curve starting from \( x \). Additionally, \( \eta_x \) is also a \((u, L, c[0])\)-calibrated curve on any subinterval \([0, t]\) of \([0, \tau_u(x)]\).

Furthermore, we have
\begin{align*}
    T_{\tau_u(x)}^- u(\eta_x(\tau_u(x))) + c[0] \tau_u(x) &= u(\eta_x(\tau_u(x))) \\
    &= u(x) + A_{\tau_u(x)}(x, \eta_x(\tau_u(x))) + c[0] \tau_u(x),
\end{align*}
which implies that \( x \in \arg\min\{ u(\cdot) + A_{\tau_u(x)}(\cdot, \eta(\tau_u(x))) \}. \) By Fermat's rule,
\begin{align*}
    Du(x)= D_x A_{\tau_u(x)}(x, \eta_x(\tau_u(x)))= L_v(\eta_x(0), \dot{\eta}_x(0)).
\end{align*}
Thus, the first equality in \eqref{eq:u-cut-cali} is satisfied. Note that \( u \) is differentiable at the interior points of calibrated curves. Consequently, for \( t \in (0, \tau_u(x)) \),
\begin{align*}
    T_t^- u(\eta_x(t)) + c[0] t = u(\eta_x(t))= u(x) + A_t(x, \eta_x(t)) + c[0] t,
\end{align*}
and
\begin{align*}
    D T_t^- u(\eta_x(t)) = D u(\eta_x(t))= D_y A_t(x, \eta_x(t))= L_v(\eta_x(t), \dot{\eta}_x(t)),
\end{align*}
which confirms that the second equality in \eqref{eq:u-cut-cali} holds true.
\end{proof}

\begin{proof}[Proof of Lemma \ref{lem:smooth-approx}]
%The proof is partly inspired by \cite[Theorem 5.2]{Ambrosio_Brue_Semola_book2021}.
Without loss of generality we suppose $f\not\equiv0$. Observe that there exists $\delta_f>0$ depending on $f$ such that $I_m+tD^2f$ is positively definite for $t\in[0,\delta_f]$, i.e., $\mathrm{id}+tDf$ is the gradient of smooth convex function $\psi_t:=\frac{1}{2}|\cdot|^2+tf(\cdot)$. It is easy to check that
\begin{align*}
	\Gamma:=\{(x,D\psi_t(x)):x\in\R^m\}
\end{align*} 
is a $c$-cyclically monotone closed graph with $c(x,y)=\frac{1}{2}|x-y|^2$. Then for $t\in[0,\delta_f]$,  the map $\mathrm{id}+tDf$ is optimal and 
\begin{align*}
	[\mathrm{id}\times(\mathrm{id}+tDf)]_\#\mu\in\Gamma_2(\mu,\nu_t).
\end{align*}
By inequality \eqref{eq:Wass_dist_push_forw} for any $s,t\in[0,\delta_f]$,
\begin{align*}
	W_2(\nu_t,\nu_s)=&\,W_2((\mathrm{id}+t Df)_\#\mu,(\mathrm{id}+s Df)_\#\mu)\\
	\leqslant&\,\left\{\int_{\R^{m}}(t-s)^2|Df|^2d\mu\right\}^{\frac{1}{2}}\leqslant |t-s|\delta_f^{-1}W_2(\mu,\nu_{\delta_f}),
\end{align*}
it follows $\{(\mathrm{id}+tDf)_\#\mu\}_{0\leqslant t\leqslant \delta_f}$ is a constant speed geodesic in $\mathscr P_2(\R^m)$. Invoking  Theorem \ref{thm:wasser-geo}, $\Gamma_2(\mu,\nu_t)$ contains a unique transport plan for $t\in[0,\delta_f]$. That is
\begin{align*}
	\Gamma_2(\mu,\nu_t)=\left\{\left[\mathrm{id}\times \left(\mathrm{id}+tDf\right)\right]_\#\mu\right\},\qquad\forall t\in[0,\delta_f].
\end{align*} 

Since $\mu$ is compactly supported and $f\in C_c^\infty(\R^m)$, we claim that for any $R>0$, there exists a $\delta_{f,R}>0$ such that
\begin{equation}\label{eq:support1}
	(\mathrm{id}+tDf)(\supp(\mu))\subset B_R(\mu),\qquad t\in[0,\delta_{f,R}].
\end{equation}
Indeed, we take $\delta_{f,R}=\min\{\delta_f,R/\|f\|_{C^1}\}$. Thus
\begin{align*}
	\left|\left(\mathrm{id}+tDf\right)(x)-x\right|=\left|tDf(x)\right|\leqslant t\|f\|_{C^1}\leqslant R,\qquad\forall x\in\R^m
\end{align*}
and \eqref{eq:support1} follows. By \eqref{eq:support1} we obtain
\begin{align*}
	\nu_t(\R^m\setminus B_R(\mu))=\mu\left(\left[\mathrm{id}+tDf\right]^{-1}(\R^m\setminus B_R(\mu))\right)\leqslant\mu(\R^m\setminus \supp(\mu))=0.
\end{align*}
This implies that $\nu_t\in\mathscr P(B_R(\mu))\subset\mathscr P_c(\R^m)$.
\end{proof}

\begin{proof}[Proof of Proposition \ref{pro:diff-singleton}]
Given $\alpha\in\partial^+U(\mu)$ and $\beta\in\partial^-U(\mu)$. Then for any $R>0$, $\nu\in\mathscr P(B_R(\mu))$, both \eqref{eq:wass-supdiff} and \eqref{eq:wass-subdiff} hold true. Invoking Lemma \ref{lem:smooth-approx}, for any $f\in C_c^\infty(\R^m)$ and $R>0$, $\nu_t=(\mathrm{id}+t Df)_\#\mu\in\mathscr P(B_R(\mu))$ and
\begin{align*}
	\Gamma_2(\mu,\nu_t)=\left\{\left[\mathrm{id}\times \left(\mathrm{id}+tDf\right)\right]_\#\mu\right\},
\end{align*}
for any $t\in[0,\delta_{f,R}]$. Thus
\begin{align*}
	\int_{\mathbb{R}^m}\left\langle\beta(x),tDf\right\rangle d \mu +o_R\left(t\|Df\|_{L^2(\mu)}\right) \leqslant \int_{\mathbb{R}^m}\left\langle\alpha(x),tDf\right\rangle d \mu+o_R\left(t\|Df\|_{L^2(\mu)}\right)
\end{align*}
and it follows
\begin{align*}
	\int_{\mathbb{R}^m}\left\langle\alpha(x)-\beta(x),Df\right\rangle d \mu\geqslant 0.
\end{align*}
Therefore, $\alpha=\beta$, $\mu$-a.e. since $f\in C_c^\infty(\R^m)$ is arbitrary.
\end{proof}

\begin{proof}[Proof of Proposition \ref{pro:potential-scl-equiv}]
It is easy to check that (3) implies (2) and (2) implies (1). Thus, it is sufficient to prove (1) implies (3).

Given $R>0$, we assume $\nu\in\mathscr P(B(0,R))$ and let $\gamma\in\Gamma_2(\mu,\nu)$. We define $\gamma_{1\to 2}^t$ for $\gamma$ by \eqref{eq:mu-ij}. Theorem \ref{thm:wasser-geo} implies $\Gamma_2(\mu,\gamma_{1\to 2}^t)=\{\gamma_{1,1\to 2}^{t}\}$ (defined in \eqref{eq:mu-kij}) and $\gamma_{1\to 2}^t$ is the constant speed geodesic. Thus $W_2(\mu,\gamma_{1,1\to 2}^{t})=tW_2(\mu,\nu)$ for all $t\in[0,1]$. Now suppose $\alpha\in\partial^+\phi(\mu)$. We obtain
\begin{align*}
	\phi(\gamma_{1\to 2}^t)-\phi(\mu)&\leqslant\int_{\mathbb{R}^{2m}}\langle\alpha(x), y-x\rangle\,d\gamma_{1,1\to 2}^t+o_R(W_2(\mu,\gamma_{1\to 2}^t))\\
			&=\int_{\mathbb{R}^{2m}}\langle\alpha(x), y-x\rangle\,d\gamma_{1,1\to 2}^t+o_R(tW_2(\mu,\nu))\\
			&=\int_{\mathbb{R}^{2m}}\langle\alpha(x), (1-t)x+ty-x\rangle\,d\gamma+o_R(t)\\
			&=\int_{\mathbb{R}^{2m}}\langle\alpha(x), t(y-x)\rangle\,d\gamma+o_R(t).
\end{align*}
Since $\phi(\cdot)$ is locally semiconcave, there exists $C_R>0$ such that
\begin{align*}
	(1-t)\phi(\mu)+t\phi(\nu)-\phi(\gamma_{1\to 2}^t)\leqslant C_Rt(1-t)+W_2^2(\mu,\nu).
\end{align*} 
It follows
\begin{align*}
	\phi(\nu)-\phi(\mu)-C_R(1-t)W_2^2(\mu,\nu)\leqslant \frac{\phi(\gamma_{1\to 2}^t)-\phi(\mu)}{t},
\end{align*}
and 
\begin{align*}
	\phi(\nu)-\phi(\mu)-C_RW_2^2(\mu,\nu)\leqslant \limsup_{t\to 0^+}\frac{\phi(\gamma_{1\to 2}^t)-\phi(\mu)}{t}\leqslant \int_{\R^{2m}}\langle\alpha(x),y-x\rangle\,d\gamma.
\end{align*}
Assertion (3) holds by taking the infimum over all $\gamma\in\Gamma_2(\mu,\nu)$ on the righthand side of the inequality above.
\end{proof}

\begin{proof}[Proof of Theorem \ref{thm:partial+D+}]
Since $\mu\in\mathscr P_c(\R^m)$, there exists $R>0$, such that $\supp(\mu)\subset B(0,R)$. By local semiconcavity of $\phi$, there exists $C_R>0$ such that for any $x,y\in B(0,R)$ and any $t\in[0,1]$,
	\begin{align*}
		(1-t)\phi(x)+t\phi(y)-\phi(tx+(1-t)y)\leqslant C_Rt(1-t)|x-y|^2.
	\end{align*}
	
We begin with necessity and we will prove by contradiction. Suppose there exists $A\subset B(0,R)$ with $\mu(A)>0$ such that $\alpha(x)\notin D^+\phi(x)$ for all $x\in A$. We define a set-valued map
\begin{align*}
	G:A&\rightrightarrows\overline{B(0,R)}\\
		x&\mapsto G(x):=\mathop{\arg\max}_{y\in\overline{B(0,R)}}\left\{\phi(y)-\phi(x)-\langle\alpha(x),y-x\rangle-C_R|y-x|^2\right\}.
\end{align*}
By a standard measurable selection theorem (see, for example \cite[Theorem 18.19]{Aliprantis_Border_book2006}), there is a Borel measurable selection $g:A\to\overline{B(0,R)}$ of $G$, i.e. $x\in A$, $g(x)\in G(x)$. Due to our assumption, $\alpha(x)\notin D^+\phi(x)$ and $g(x)$ is positive on $A$. The map $h:\R^m\to\R^m$,
\begin{align*}
	h(x)=
	\begin{cases}
		x,&x\not\in A\\
		g(x),&x\in A
	\end{cases}
\end{align*} 
is measurable. Set $\nu:=h_\#\mu$. Since $h(\overline{B(0,R)})\subset\overline{B(0,R)}$, we have
\begin{align*}
	h_\#\mu(\R^m\setminus\overline{B(0,R)})=\mu(h^{-1}[\R^m\setminus\overline{B(0,R)}])=\mu(\R^m\setminus[h^{-1}(\overline{B(0,R)}])\leqslant\mu(\R^m\setminus\overline{B(0,R)})=0.
\end{align*}
That is $\nu\in\mathscr P(\overline{B(0,R)})$. By Proposition \ref{pro:potential-scl-equiv} we conclude that if $\alpha\in\partial^+\phi(\mu)$ then
\begin{equation}\label{eq:xi in partial}
	\phi(\nu)-\phi(\mu)-\int_{\R^m}\langle\alpha(x),y-x\rangle\,d(\mathrm{id}\times h)_\#\mu\leqslant C_RW_{2,(\mathrm{id}\times h)_\#\mu}^2(\mu,\nu).
\end{equation}
However, the definition of $h$ and the fact that $g>0$ on $A$ imply
\begin{align*}
	&\,\phi(\nu)-\phi(\mu)-\int_{\R^m}\langle\alpha(x),y-x\rangle\,d(\mathrm{id}\times h)_\#\mu\\
		=&\,\int_{\mathbb{R}^m} \phi(h(x))-\phi(x)-\langle\alpha(x), h(x)-x\rangle\,d\mu \\
    	=&\,\int_A \phi(g(x))-\phi(x)-\left\langle\alpha(x), g(x)-x\right\rangle\, d\mu\\
    	>&\,\int_A C_R|g(x)-x|^2\,d\mu=\int_{\R^m} C_R|h(x)-x|^2\,d\mu\\
    	=&\,C_RW_{2,(\mathrm{id}\times h)_\#\mu}^2(\mu,\nu).
\end{align*}
This leads to a contradiction to \eqref{eq:xi in partial}.
	
Now we turn to prove the sufficiency part. Suppose $\phi\in\Liploc(\R^m)$, $\alpha(x)\in D^+\phi(x)$ for $\mu$-a.e. $x$. Then $\alpha\in L^\infty(\R^m;\mu)$ and
\begin{align*}
	\phi(y)-\phi(x)\leqslant \langle\alpha(x),y-x\rangle+C_R|y-x|^2,
\end{align*}
for every $y\in B(0,R)$ and $\mu$-a.e. $x\in\R^m$. Thus, for any $\nu\in\mathscr P(B(0,R))$ and $\gamma\in\Gamma_2(\mu,\nu)$, integrating both sides of the inequality above, we  get
\begin{align*}
	\phi(\nu)-\phi(\mu)&\leqslant\int_{\R^{2m}}\langle\alpha(x),y-x\rangle+C_R|y-x|^2\,d\gamma\\
		&=\int_{\R^{2m}}\langle\alpha(x),y-x\rangle\,d\gamma+C_RW_2^2(\mu,\nu).
\end{align*}
Taking the infimum over all $\gamma\in\Gamma_2(\mu,\nu)$, we have \eqref{eq:potential-scl-equiv3}. This implies $\alpha\in\partial^+\phi(\mu)$.
\end{proof}

\begin{Lem}[\protect{\cite[Lemma 5.3.2]{Ambrosio_GigliNicola_Savare_book2008}}]\label{lem:gluing lemma}
Let $X_i$ $(i=1,2,3)$ be Polish spaces, $\gamma^{1,2}\in\mathscr P(X_1\times X_2)$ and $\gamma^{1,3}\in\mathscr P(X_1\times X_3)$ with  $(\pi_1)_\#\gamma^{1,2}=(\pi_1)_\#\gamma^{1,3}=\mu^{1}$. Then there exists $\boldsymbol{\mu}\in\mathscr P(X_1\times X_2\times X_3)$, which satisfies
\begin{equation}\label{eq:disintegration}
	(\pi_{1,2})_\#\boldsymbol{\mu}=\gamma^{1,2},\qquad(\pi_{1,3})_\#\boldsymbol{\mu}=\gamma^{1,3}.
\end{equation}
		
Moreover, if $\gamma^{1,2}=\int\gamma_{x_1}^{1,2}d\mu^1$, $\gamma^{1,3}=\int\gamma_{x_1}^{1,3}d\mu^1$ and $\boldsymbol{\mu}=\int\boldsymbol{\mu}_{x_1}d\mu^1$ are the disintegrations of $\gamma^{1,2}$, $\gamma^{1,3}$ and $\boldsymbol{\mu}$ with respect to $\mu^1$ respectively, then \eqref{eq:disintegration} is equivalent to 
\begin{align*}
	\boldsymbol{\mu}_{x_1}\in\Gamma(\gamma_{x_1}^{1,2},\gamma_{x_1}^{1,3})\subset\mathscr P(X_2\times X_3)\qquad \text{for}\ \mu^1-a.e.\ x_1\in X_1.
\end{align*}
\end{Lem}

\begin{Lem}\label{lem:gluing_lambda}
Suppose $\mu_i\in\mathscr P(\R^m)$ with $i=1,2,3$ and $\gamma\in\Gamma(\mu_2,\mu_3)$. For any fixed $\lambda\in[0,1]$ and $t>0$, let $\mu_{\lambda}\in\Gamma_o^{t}(\mu_1,\gamma_{2\to 3}^\lambda)$. Then there exists $\gamma_\lambda\in\Gamma(\mu_1,\gamma)$ such that $(\gamma_\lambda)_{1,2\to 3}^\lambda=\mu_\lambda$. 
\end{Lem}

\begin{proof}
The following is inspired by the proof of \cite[Proposition 7.3.1]{Ambrosio_GigliNicola_Savare_book2008}. Observe that maps
\begin{align*}
	\begin{aligned}
		\Sigma_\lambda:\R^{2m}&\to\R^{2m}\\
			(x_2,x_3)&\mapsto (x_2,(1-\lambda)x_2+\lambda x_3)
	\end{aligned}
	\qquad\text{and}\qquad
	\begin{aligned}
		\Lambda_\lambda:\R^{3m}&\to\R^{3m}\\
			(x_1,x_2,x_3)&\mapsto (x_1,x_2,(1-\lambda)x_2+\lambda x_3)
	\end{aligned}
\end{align*}
are homeomorphisms for any fixed $\lambda\in(0,1]$. Thus, $\gamma_\lambda\in\Gamma(\mu_1,\gamma)$ satisfies $(\gamma_\lambda)_{1,2\to 3}^\lambda=\mu_\lambda$ if and only if for $\nu_\lambda:=(\Lambda_\lambda)_\#\gamma_\lambda$,
\begin{equation}\label{eq:glue-homeo}
	(\pi_{1,3})_\#\nu_\lambda=\mu_\lambda,\qquad (\pi_{2,3})_\#\nu_\lambda=(\Sigma_\lambda)_\#\gamma.
\end{equation}
Then Lemma \ref{lem:gluing lemma} implies the existence of $\nu_\lambda$ and this ensures \eqref{eq:glue-homeo} hold true. This guarantees the existence of $\gamma_\lambda$ since $\Lambda_\lambda$ is invertible. 
		
For the case of $\lambda=0$, we only need to check that the existence of $\gamma_0\in\Gamma(\mu_1,\gamma)$ satisfying $(\pi_{1,2})_\#\gamma_0=\nu$ for each $\nu\in\Gamma^t_o(\mu_1,\mu_2)$. This  can be guaranteed by Lemma \ref{lem:gluing lemma} as well.
\end{proof}

\begin{proof}[Proof of Theorem \ref{thm:dyn-cost}]
Note that for any $v\in\R^m$, $L(x,v)\geqslant\theta_0(|v|)-c_0$, where $\theta_0$ is convex, increasing and superlinear function in $\R^+$. We assume that $\xi\in L_{\mu,\nu}^{t}$ is the dynamical optimal coupling of $C^{t}(\mu,\nu)$, then by Jensen's inequality and the monotonicity of $\theta_0$, for $\mu,\nu\in\mathscr P_c(\R^m)$,
\begin{align*}
	C^{t}(\mu,\nu)
    		&=\E\left(\int_0^tL(\xi(s),\dot\xi(s))\,ds\right)\\
    		&\geqslant \E\left(\int_0^t\theta_0(|\dot\xi(s)|)-c_0\,ds\right)\\
    		&\geqslant t\E\left(\theta_0\left(\left|\frac{1}{t}\int_0^t\dot\xi(s)\,ds\right|\right)-c_0\right)\\
    		&=t\E\left(\theta_0\left(\left|\frac{\xi(t)-\xi(0)}{t}\right|\right)-c_0\right)\\
    		&\geqslant t\theta_0\left(\E\left(\left|\frac{\xi(t)-\xi(0)}{t}\right|\right)\right)-c_0t\\
    		&\geqslant t\theta_0\left(\frac{1}{t}W_1(\mu,\nu)\right)-c_0t.
\end{align*}
Thus, for any $k\in\R$,
\begin{align*}
	C^{t}(\mu,\nu)\geqslant t\theta_0\left(\frac{1}{t}W_1(\mu,\nu)\right)-c_0t\geqslant kW_1(\mu,\nu)-(c_0+\theta_0^*(k))t,
\end{align*} 
where $\theta_0^*$ is Fenchel-Legendre duality of $\theta_0$, which implies $C^{t}(\mu,\cdot)$ is superlinear on $\mathscr P_1(\R^m)$. The proof of superlinearity of $C^t(\cdot,\mu)$ on $\mathscr P_c(\R^m)$ is similar.
    	
Given $\mu,\nu_1,\nu_2\in\mathscr P(K)$ and $t_1,t_2\in[a,b]$, we assume $\gamma_2\in\Gamma_o^{t_2}(\mu,\nu_2)$ and $\gamma_{1,2}\in\Gamma_1(\nu_1,\nu_2)$. Lemma \ref{lem:gluing lemma} provides $\gamma_2*\gamma_{1,2}\in\Gamma(\nu_1,\nu_2,\mu)\in\mathscr P(\R^{3m})$ with 
$(\pi_{3,1})_\#(\gamma_2*\gamma_{1,2})\in\Gamma(\mu,\nu_1)$ and $(\pi_{3,2})_\#(\gamma_2*\gamma_{1,2})=\gamma_2$. Therefore,
\begin{align*}
	&\,C^{t_1}(\mu,\nu_1)-C^{t_2}(\mu,\nu_2)\\
    		\leqslant&\,\int_{\R^{2m}}A_{t_1}(x,y_1)\,d(\pi_{3,1})_\#(\gamma_2*\gamma_{1,2})-\int_{\R^{2m}}A_{t_2}(x,y_2)\,d\gamma_{2}\\
    		=&\,\int_{\R^{3m}}A_{t_1}(x,y_1)-A_{t_2}(x,y_2)\,d(\gamma_2*\gamma_{1,2})\\
    		\leqslant&\,\int_{\R^{3m}}l_{K,a,b}(d(y_1,y_2)+|t_1-t_2|)\,d(\gamma_{2}*\gamma_{1,2})\\
    		=&\,l_{K,a,b}(W_1(\nu_1,\nu_2)+|t_1-t_2|).
\end{align*} 
where $l_{K,a,b}$ denotes the Lipschitz constant of function $(t,y)\mapsto A_t(x,y)$ depending on $K,a$ and $b$.

Suppose $\mu_i\in\mathscr P(K)$, $i=1,2,3$, $\gamma\in\Gamma(\mu_2,\mu_3)$ and $t\in(0,\tau_1]$ fixed. By Lemma \ref{lem:gluing_lambda} we have that for any $\gamma\in\Gamma(\mu_2,\mu_3)$ and $\mu_{\lambda}\in\Gamma_o^{t}(\mu_1,\gamma_{2\to 3}^\lambda)$, there exists $\gamma_\lambda\in\Gamma(\mu_1,\gamma)$ satisfying $(\gamma_\lambda)_{1,2\to 3}^\lambda=\mu_\lambda$. We also have $(\pi_{1,2})_\#\gamma_\lambda\in\Gamma(\mu_1,\mu_2)$ and $(\pi_{1,3})_\#\gamma_\lambda\in\Gamma(\mu_1,\mu_3)$. In addition, for fixed $t\in(0,\tau_1]$,  there exists some $\lambda_K>0$ such that 
    	$K\subset B(x,\lambda_Kt)$ for all $x\in K$. Now we denote $C_{K}:=C_{\lambda_K}$ in Proposition \ref{pro:cvpde}, then
\begin{align}
	&\,(1-\lambda)C^{t}(\mu_1,\mu_2)+\lambda C^{t}(\mu_1,\mu_3)-C^{t}(\mu_1,\gamma_{2\to 3}^\lambda)\notag\\
    	  	=&\,(1-\lambda)\inf_{\boldsymbol{\mu}\in\Gamma(\mu_1,\mu_2)}\int_{\R^{2m}}A_{t}(x,y)\,d\boldsymbol{\mu}+\lambda\inf_{\boldsymbol{\mu}\in\Gamma(\mu_1,\mu_3)}\int_{\R^{2m}}A_{t}(x,y)\,d\boldsymbol{\mu}\\
    	  	&\,\qquad\qquad-\inf_{\boldsymbol{\mu}\in\Gamma(\mu_1,\gamma_{2\to 3}^\lambda)}\int_{\R^{2m}}A_{t}(x,y)\,d\boldsymbol{\mu}\notag\\
    	  	\leqslant&\,\int_{\R^{2m}}(1-\lambda)A_{t}(x,y)\,d(\pi_{1,2})_\#\gamma_\lambda+\int_{\R^{2m}}\lambda A_{t}(x,y)\,d(\pi_{1,3})_\#\gamma_\lambda-\int_{\R^{2m}}A_{t}(x,y)\,d\mu_{\lambda}\notag\\
    	  	=&\,\int_{\R^{3m}}\left[(1-\lambda)A_{t}(x_1,x_2)+\lambda A_{t}(x_1,x_3)-A_{t}(x_1,(1-\lambda)x_2+\lambda x_3)\right]\,d\gamma_\lambda\notag\\
    	  	\leqslant&\,\lambda(1-\lambda)\int_{\R^{3m}}\frac{C_K}{t}|x_2-x_3|^2d\gamma_\lambda\notag\\
    	  	=&\,\lambda(1-\lambda)\int_{\R^{2m}}\frac{C_K}{t}|x_2-x_3|^2d\gamma=\frac{C_K}{t}\lambda(1-\lambda)W_{2,\gamma}^2(\mu_2,\mu_3).\label{eq:meaure-semiconcave}
\end{align}
Consequently, \eqref{eq:meaure-semiconcave} implies $C^t(\mu,\cdot)$ is locally strongly semiconcave on $\mathscr P(K)$ with constant $C_K/t$. 
\end{proof}

\begin{proof}[Proof of Theorem \ref{thm:cost-supdiff}]
For $\nu,\nu^\prime\in\mathscr P_c(\R^m)$, $\gamma^\prime\in\Gamma_2(\nu,\nu')$ and $\gamma\in\Gamma_{o}^{t}(\mu,\nu)$, we know $(\pi_2)_\#\gamma=(\pi_2)_\#\gamma^\prime=\nu$. Following from Lemma \ref{lem:gluing lemma}, $\gamma*\gamma^\prime\in\Gamma(\mu,\nu,\nu^\prime)$ with $(\pi_{2,3})_\#(\gamma*\gamma^\prime)=\gamma^\prime$ and $(\pi_{1,2})_\#(\gamma*\gamma^\prime)=\gamma$. Simultaneously, $(\pi_{1,3})_\#(\gamma*\gamma^\prime)\in\Gamma(\mu,\nu^\prime)$. Here we can get
\begin{align*}
	&\,C^{t}(\mu,\nu^\prime)-C^{t}(\mu,\nu)\\
		=&\,\inf_{\boldsymbol{\mu}\in\Gamma(\mu,\nu^\prime)}\int_{\R^{2m}}A_t(x,y)\,d\boldsymbol{\mu}-\inf_{\boldsymbol{\mu}\in\Gamma(\mu,\nu)}\int_{\R^{2m}}A_t(x,y)\,d\boldsymbol{\mu}\\
		\leqslant&\,\int_{\R^{2m}}A_t(x,y)\,d(\pi_{1,3})_\#(\gamma*\gamma^\prime)-\int_{\R^{2m}}A_t(x,y)\,d(\pi_{1,2})_\#(\gamma*\gamma^\prime)\\
		=&\,\int_{\R^{3m}}A_t(x_1,x_3)-A_t(x_1,x_2)\,d(\gamma*\gamma^\prime)\\
		\leqslant&\,\int_{\R^{3m}}\left(\langle L_v(\eta_{x_1,x_2}(t),\dot\eta_{x_1,x_2}(t)),x_3-x_2\rangle+\frac{C_{\supp(\gamma*\gamma^\prime)}}{2t}|x_3-x_2|^2\right)\,d(\gamma*\gamma^\prime),
\end{align*}
where $C_{\supp(\gamma*\gamma^\prime)}$ is a constant depending on $\supp(\gamma*\gamma^\prime)$. Recall the construction of $\gamma*\gamma^\prime$ mentioned in the proof of Lemma \ref{lem:gluing lemma}, $
	\gamma*\gamma^\prime=\int_{\R^m}\gamma_{x_2}\times\gamma_{x_2}^\prime\,d\nu(x_2),
$ where $\gamma=\int_{\R^m}\gamma_{x_2}\,d\nu(x_2)$ and $\gamma^\prime=\int_{\R^m}\gamma_{x_2}^\prime\,d\nu(x_2)$. It follows 
\begin{align*}
	&\,C^{t}(\mu,\nu^\prime)-C^{t}(\mu,\nu)\\
	\leqslant&\,\int_{\R^{3m}}\left(\langle p_{x_1,x_2}(t),x_3-x_2\rangle+\frac{C_{\supp(\gamma*\gamma^\prime)}}{2t}|x_3-x_2|^2\right)\,d(\gamma*\gamma^\prime)\\
	=&\,\int_{\R^{m}}\int_{\R^{2m}}\left(\langle p_{x_1,x_2}(t),x_3-x_2\rangle+\frac{C_{\supp(\gamma*\gamma^\prime)}}{2t}|x_3-x_2|^2\right)d(\gamma_{x_2}\times\gamma_{x_2}^\prime)d\nu\\
	=&\,\int_{\R^{m}}\int_{\R^{m}}\left(\langle \int_{\R^{m}}p_{x_1,x_2}(t)\,d\gamma_{x_2},x_3-x_2\rangle+\frac{C_{\supp(\gamma*\gamma^\prime)}}{2t}|x_3-x_2|^2\right)d\gamma_{x_2}^\prime d\nu\\
	=&\,\int_{\R^{2m}}\langle\mathbf p_\nu(x_2),x_3-x_2\rangle+\frac{C_{\supp(\gamma*\gamma^\prime)}}{2t}|x_3-x_2|^2\,d\gamma^\prime\\
	=&\,\int_{\R^{2m}}\langle\mathbf p_\nu(x_2),x_3-x_2\rangle\,d\gamma^\prime+\frac{C_{\supp(\gamma*\gamma^\prime)}}{2t}W_2^2(\nu,\nu'),
\end{align*}
where $p_{x,y}(t):=L_v(\eta_{x,y}(t),\dot\eta_{x,y}(t))$. Since our choice of $\gamma^\prime\in\Gamma_2(\mu,\nu)$ is arbitrary, Definition \ref{def:local frechet} shows that $\mathbf p_\nu(x_2)\in \partial_\nu^+C^{t}(\mu,\nu)$. $-\mathbf p_\mu(x_1)\in \partial_\mu^+C^{t}(\mu,\nu)$ can be obtained by a similar argument. For the rest of the proof, we assume that $\gamma\in\Gamma_o^{t}(\mu,\nu)$ and $t^\prime\in(t-\varepsilon,t+\varepsilon)$ with $\varepsilon>0$, 
\begin{align*}
	C^{t^\prime}(\mu,\nu)-C^{t}(\mu,\nu)&\leqslant \int_{\R^{2m}}A_{t^\prime}(x,y)\,d\gamma-\int_{\R^{2m}}A_t(x,y)\,d\gamma\\
	&\leqslant-(t^\prime-t)\int_{\R^{2m}}H(\eta_{x,y}(t),p_{x,y}(t))\,d\gamma+\frac{C_{K,t,\varepsilon}}{2}|t^\prime-t|^2,
\end{align*}
where $K=\text{supp}\,\mu\times\text{supp}\,\nu$ is a compact subset. This estimate completes the proof. 
\end{proof}

\bibliographystyle{plain}
\bibliography{mybib}
\end{document}